\documentclass[aap,preprint]{imsart}

\RequirePackage{amsthm,amsmath,amsfonts,amssymb,mathtools,mathrsfs}
\RequirePackage[authoryear]{natbib}
\RequirePackage{graphicx}
\RequirePackage{tikz}
\usetikzlibrary{calc,arrows.meta}
\RequirePackage{float}
\RequirePackage{anyfontsize}
\usetikzlibrary{matrix,positioning,decorations.pathreplacing}
\RequirePackage[shortlabels]{enumitem}
\RequirePackage{xcolor}
\RequirePackage[colorlinks,citecolor=blue,urlcolor=blue,linkcolor=blue]{hyperref}
\RequirePackage[nameinlink,noabbrev]{cleveref}
\allowdisplaybreaks
\DeclareFontShape{T1}{ptm}{m}{scit}{<->ssub * ptm/m/sc}{}

\startlocaldefs
\numberwithin{equation}{section}

\newcommand{\R}{\mathbb R}
\newcommand{\E}{\mathbb E}
\newcommand{\bP}{\mathbb P}
\newcommand{\bN}{\mathbb N}
\newcommand{\tld}{\widetilde}
\newcommand{\ind}{\mathbf 1}
\newcommand{\Cov}{\operatorname{Cov}}
\newcommand{\Var}{\operatorname{Var}}
\newcommand{\dto}{\overset{D}{\longrightarrow}}
\newcommand{\pto}{\overset{P}{\longrightarrow}}
\DeclareMathOperator{\InjHom}{inj}
\DeclareMathOperator{\supp}{supp}
\DeclareMathOperator{\Bin}{Bin}

\theoremstyle{plain}
\newtheorem{theorem}{Theorem}[section]
\newtheorem{lemma}[theorem]{Lemma}
\newtheorem{corollary}[theorem]{Corollary}
\newtheorem{proposition}[theorem]{Proposition}

\theoremstyle{definition}
\newtheorem{definition}[theorem]{Definition}

\newtheorem{assumption}{}

\newtheorem{remark}[theorem]{Remark}
\newtheorem{example}{Example}

\crefname{assumption}{}{}
\Crefname{assumption}{}{}
\crefname{equation}{equation}{equations}
\Crefname{equation}{Equation}{Equations}
\crefname{definition}{Definition}{Definitions}
\Crefname{definition}{Definition}{Definitions}
\crefname{theorem}{Theorem}{Theorems}
\Crefname{theorem}{Theorem}{Theorems}
\crefname{lemma}{Lemma}{Lemmas}
\Crefname{lemma}{Lemma}{Lemmas}
\crefname{proposition}{Proposition}{Propositions}
\Crefname{proposition}{Proposition}{Propositions}
\crefname{corollary}{Corollary}{Corollaries}
\Crefname{corollary}{Corollary}{Corollaries}
\crefname{section}{Section}{Sections}
\Crefname{section}{Section}{Sections}
\crefname{figure}{Figure}{Figures}
\Crefname{figure}{Figure}{Figures}
\crefname{remark}{Remark}{Remarks}
\Crefname{remark}{Remark}{Remarks}
\crefname{example}{Example}{Examples}
\Crefname{example}{Example}{Examples}

\endlocaldefs

\begin{document}

\begin{frontmatter}

\title{Distributional Limits of Empirical Codegree Processes in Sampled Networks}
\runtitle{Empirical Codegree Processes in Sampled Networks}

\begin{aug}
\author[A]{\fnms{Victor}~\snm{Daniel}\ead[label=e1]{vdd2111@columbia.edu}}
\author[A]{\fnms{Sumit}~\snm{Mukherjee\thanks{Supported by NSF DMS-2515519.}}\ead[label=e2]{sm3949@columbia.edu}}
\address[A]{Department of Statistics, Columbia University
\printead[presep={,\ }]{e1,e2}}
\end{aug}

\begin{abstract}
We study fluctuations of empirical distributions of normalized codegrees under independent vertex sampling from deterministic or random weighted parent networks. We establish a quantitative multivariate normal approximation for sampled biclique counts. Under additional graph-limit and sparsity assumptions, we establish convergence of the finite-dimensional distributions of the empirical codegree process indexed by smooth test functions. We further establish functional convergence on H\"older balls, with a smoothness threshold depending on network sparsity and vertex sampling probabilities. We verify our assumptions for classes of graphon-generated random graphs and give counterexamples showing that the smoothness conditions cannot in general be weakened.
\end{abstract}



\pdfsubject{Preprint}

\end{frontmatter}

\section{Introduction}\label{sec:intro}

When only a subgraph of a network is observed, a central statistical problem is to infer properties of the parent network while accounting for the sampling mechanism. This problem has a long history, including early
work on snowball sampling, graph totals, and connected components \cite{10.1214/aoms/1177705148,MR458659,MR515656}; see also \citet{Frank_2005} for a broader account. The degree distribution is a basic example of a summary feature of the parent graph.
Its estimation from sampled subgraphs already appears in the work of
\citet{Frank1980DegreeFrequencies}, who constructed unbiased estimators
of degree frequencies. More recently,
\citet{ZhangKolaczykSpencer2015} formulated degree-distribution recovery
as an inverse problem and developed a regularized estimation procedure.

In this paper, we consider the subgraph sampling
design: each vertex of a fixed $n$-vertex parent network $G_n$ is chosen independently with
probability $p_n$, and the induced subgraph $\tld G_n$ is observed (\cref{def:induced-vertex-sampling}). Our objective is
to characterize how an empirical distribution of local connectivity
statistics changes under this sampling design. In particular,  we study the empirical distribution of normalized $r$-codegrees.
For an unweighted graph, the $r$-codegree of an $r$-element vertex set is
the number of its common neighbors; we also allow bounded nonnegative
edge weights. The case $r=1$ corresponds to the normalized degree distribution.
For $r=2$, the statistic describes pairwise neighborhood overlap and,
in an unweighted graph, is a normalized version of the dyadwise
shared-partner distribution used in social-network modeling
\cite{HunterGoodreauHandcock2008}.

Vertex sampling in the subgraph sampling model changes the local statistic as well as the collection of sets on which it is evaluated. Thus, the sampled codegrees cannot
be treated as independent observations from the parent codegree distribution. We study how these coupled effects determine the fluctuations of the codegree distribution, viewed through smooth test functions.
Our main result (\cref{thm:empirical-codegree-process}) is a functional central limit theorem for the empirical codegree distribution, centered about the codegree distribution of the parent graph. Write $\nu_n$ and $\widetilde\nu_n$ for the parent and sampled
empirical codegree measures, respectively, with each codegree suitably normalized (see \eqref{eq:codegree-distribution-definition} and \eqref{eq:nutld}).
For a test function $f$, the process of interest is
\begin{align}
    T_n(f):=\sqrt{np_n}\left(
        \int f\,d\widetilde\nu_n-\int f\,d\nu_n
    \right).
    \label{eq:process-def}
\end{align}
Under the structural and sparsity conditions stated below, $T_n$
converges weakly to a centered Gaussian process in
$\ell^\infty(\mathcal F)$, where $\mathcal F$ is a $C^{1,\alpha}$ unit
ball and the required smoothness depends explicitly on sparsity.
The centering is the empirical measure of the finite parent network,
not the distribution induced by a limiting graphon. The theorem
therefore describes sampling uncertainty jointly over a whole class
of smooth summaries.

The connection between codegree moments and biclique counts links our analysis to the sampling theory of subgraph counts. \citet{pmlr-v75-klusowski18a} obtained minimax bounds for motif
estimation under subgraph and neighborhood sampling, and
\citet{KlusowskiWu2020} studied estimation of the number of connected
components from sampled subgraphs. Under independent induced vertex
sampling, \citet{bdm22} established consistency and asymptotic normality
criteria for Horvitz--Thompson motif-count estimators, together with
quantitative normal approximation bounds. Codegree moments are closely
related to biclique counts (\cref{lem:counting-lemma}). We obtain a
quantitative joint Gaussian approximation for normalized sampled
biclique counts, with error of order $(np_n)^{-1/2}$ in a smooth
probability distance, under the conditions of
\cref{thm:stein-error-bound}. This approximation uses the finite-sample
covariance and requires no limiting graphon. Under the additional
conditions of \cref{cor:biclique-clt}, we obtain a joint central limit theorem for sampled biclique densities centered at their parent-network counterparts. These results provide the starting point for establishing the empirical-process limit.

\citet{Janson2026UncoveringGraph} studies the successive uncovering
of vertices of a fixed graph and establishes a functional Gaussian
limit theorem for the evolving edge count. His continuous-time
formulation includes independent vertex sampling. While his process
is indexed by sampling time, ours is indexed by smooth test functions
and captures fluctuations of the empirical codegree distribution
under vertex sampling.

A parallel line of research concerns fluctuations of network
statistics under random graph models. For graphon-generated
graphs, \citet{HladkyPelekisSileikis2021} studied fluctuations
of clique counts, while \citet{BhattacharyaChatterjeeJanson2023}
developed a corresponding theory for general fixed subgraphs,
identifying Gaussian and potentially non-Gaussian limiting
distributions according to graphon regularity conditions.
Building on this theory, \citet{ChatterjeeDanBhattacharya2024}
established joint higher-order fluctuation results and developed
bootstrap procedures for inference on network moments.
Related results for local network statistics include
\citet{DelmasDhersinSciauveau2021}, who studied empirical degree
distributions and partially labeled homomorphism densities,
and \citet{Maugis2024}, who established central limit theorems
for rooted subgraph statistics in inhomogeneous random graphs.
More recently, \citet{Stewart2025EmpiricalDistributions}
obtained concentration and uniform consistency results for
empirical distributions of graph statistics under dependent-edge
models. Unlike these model-based results, our main theorems
describe sampling fluctuations conditional on a fixed parent network, and do not require the true parent network to be generated from a particular model.

Our structural assumptions are formulated in the language of graph
limits \cite{BCLSV08,Lovasz2012,BCCZ2019LpI,BCCZ2018LpII}. We require
convergence of normalized star kernels in a higher-order cut norm and
uniform bounds on normalized common-neighborhood counts (\cref{ass:star-cut-convergence,ass:bounded-codegrees}). Together with
rate conditions relating edge density to expected sample size, they
allow both dense and sparse parent sequences and sampling fractions
tending to zero (\cref{ass:sampling-ratios,ass:sparsity}). We also verify these structural assumptions for classes
of graphon-generated networks, so that the sampling results apply
conditionally on almost every realized parent sequence (\cref{thm:sparse-graphon-samples}).

The limiting process covariance has an explicit representation in terms of
vertex influences (\eqref{eq:limiting-influence-function} and \eqref{eq:limiting-covariance-kernel}). We also give
counterexamples showing that sampling bias and failure of tightness
impose distinct smoothness restrictions (\cref{prop:regularity-counterexample,prop:degree-tightness-counterexample}).

\subsection{Asymptotic notation}
\label{subsec:asymptotic-notation}

For deterministic nonnegative sequences, $a_n=O(b_n)$ and $a_n=o(b_n)$
have their usual meanings, while $a_n\lesssim b_n$ means
$a_n\le Cb_n$ eventually for a constant independent of $n$; we use
$\lesssim_\square$ to display any additional dependence of that constant.
For random variables,
$Z_n=O_P(b_n)$ means that $Z_n/b_n$ is tight and $Z_n=o_P(b_n)$ that
$Z_n/b_n\pto0$. We write $Z_n=O_{L^2}(b_n)$ if
$\E[|Z_n|^2]=O(b_n^2)$.

\subsection{Setup and sampling design}\label{subsec:setup}

A weighted network $G_n$ is a simple graph with vertex set $V(G_n)=\{1,2,\cdots,n\}=[n]$ and real-valued edge weights $G_n(i,j)$.  We set
$G_n(i,j)=0$ when $ij\notin E(G_n)$, and in particular $G_n(i,i)=0$.  For
$S\subseteq [n]$, define
\begin{align*}
    G_n(S,v):=\prod_{u\in S}G_n(u,v),
    \qquad
    d_S:=\sum_{v\in [n]}G_n(S,v).
\end{align*}
For fixed $r\ge1$, the edge density $q_n$, and the normalized empirical $r$-codegree measure $\nu_n$ of $G_n$, are defined by
\begin{align}
    q_n:=\frac1{n^2}\sum_{i,j=1}^nG_n(i,j),
    \qquad
    \nu_n:=\frac{1}{\binom{n}{r}}
    \sum_{S\in\binom{[n]}{r}}
    \delta_{d_S/(nq_n^r)}
    \label{eq:codegree-distribution-definition}
\end{align}
We work with $n\ge r$ and $q_n>0$ when using parent-network
normalizations; the rate conditions below ensure this for all sufficiently
large $n$.

\begin{definition}[Subgraph sampling design]\label[definition]{def:induced-vertex-sampling}
Given $G_n$ as above, and $p_n\in(0,1)$, choose vertices independently with probability $p_n$
and observe the induced weighted subnetwork on the chosen vertices.  The
sampled graph is denoted by $\tld G_n$.
\end{definition}

\begin{figure}
\centering
\begingroup
\definecolor{codegreeblue}{RGB}{0,92,153}

\def\keep#1{\ifcase#1 0%
  \or1\or1\or1\or0\or1\or1\or0\or1\or0\or1%
  \or1\or0\or1\or0\or1\or1\or0\or1\or1\or0\fi}

\def\edges{%
  1/2,2/3,3/4,4/5,5/6,6/7,7/8,8/9,
  9/10,10/11,11/12,12/1,
  13/14,14/15,15/16,16/17,17/18,18/13,
  13/1,13/2,14/3,14/4,15/5,15/6,
  16/7,16/8,17/9,17/10,18/11,18/12,
  19/13,19/17,19/18,20/14,20/15,20/16,19/20%
}

\begin{tikzpicture}[
  x=1cm,y=1cm,font=\small,
  edge/.style={draw=black!65,line width=.55pt},
  common/.style={draw=codegreeblue,line width=1.1pt},
  vertex/.style={
    circle,draw=black,line width=.8pt,
    inner sep=2pt,outer sep=0pt
  },
  retained/.style={vertex,fill=black},
  omitted/.style={vertex,fill=white}
]

\foreach \panel/\offset in {0/0,1/6.8} {
  \begin{scope}[xshift=\offset cm]

    \foreach \i in {1,...,12}
      \coordinate (v\i) at ({90-30*(\i-1)}:2.6);

    \foreach \i in {13,...,18}
      \coordinate (v\i) at ({90-60*(\i-13)}:1.35);

    \coordinate (v19) at (-.48,0);
    \coordinate (v20) at (.48,0);

    \foreach \a/\b in \edges {
      \pgfmathtruncatemacro{\ia}{\a}
      \pgfmathtruncatemacro{\ib}{\b}

      \ifnum\panel=0
        \draw[edge] (v\ia)--(v\ib);
      \else
        \ifnum\keep{\ia}=1
          \ifnum\keep{\ib}=1
            \draw[edge] (v\ia)--(v\ib);
          \fi
        \fi
      \fi
    }

    \foreach \a/\b in {13/2,2/3,13/14,14/3} {
      \ifnum\panel=0
        \draw[common] (v\a)--(v\b);
      \else
        \ifnum\keep{\a}=1
          \ifnum\keep{\b}=1
            \draw[common] (v\a)--(v\b);
          \fi
        \fi
      \fi
    }

    \foreach \i in {1,...,20} {
      \ifnum\keep{\i}=1
        \node[retained] at (v\i) {};
      \else
        \ifnum\panel=0
          \node[omitted] at (v\i) {};
        \fi
      \fi
    }

    \node[above left=3pt,text=codegreeblue]
      at (v13) {$1$};
    \node[right=3pt,text=codegreeblue]
      at (v3) {$2$};
    \node[above right=3pt,text=codegreeblue]
      at (v2) {$3$};

    \ifnum\panel=0
      \node[left=3pt,text=codegreeblue]
        at (v14) {$4$};
    \fi

  \end{scope}
}

\node[align=center] at (0,-3.2) {
  (a) Parent graph $G_n$. Sampled vertices are shown in black.
};

\node[align=center] at (6.8,-3.2) {
  (b) Sampled subgraph $\widetilde G_n$
};
\end{tikzpicture}
\endgroup
\caption{Induced vertex sampling on an unweighted graph. Vertices $1$ and $2$ have two common neighbors in $G_n$, vertices $3$ and $4$, but only vertex $3$ is in $\widetilde G_n$, reducing their codegree from $2$ to $1$.}
\label{fig:sampling}
\end{figure}
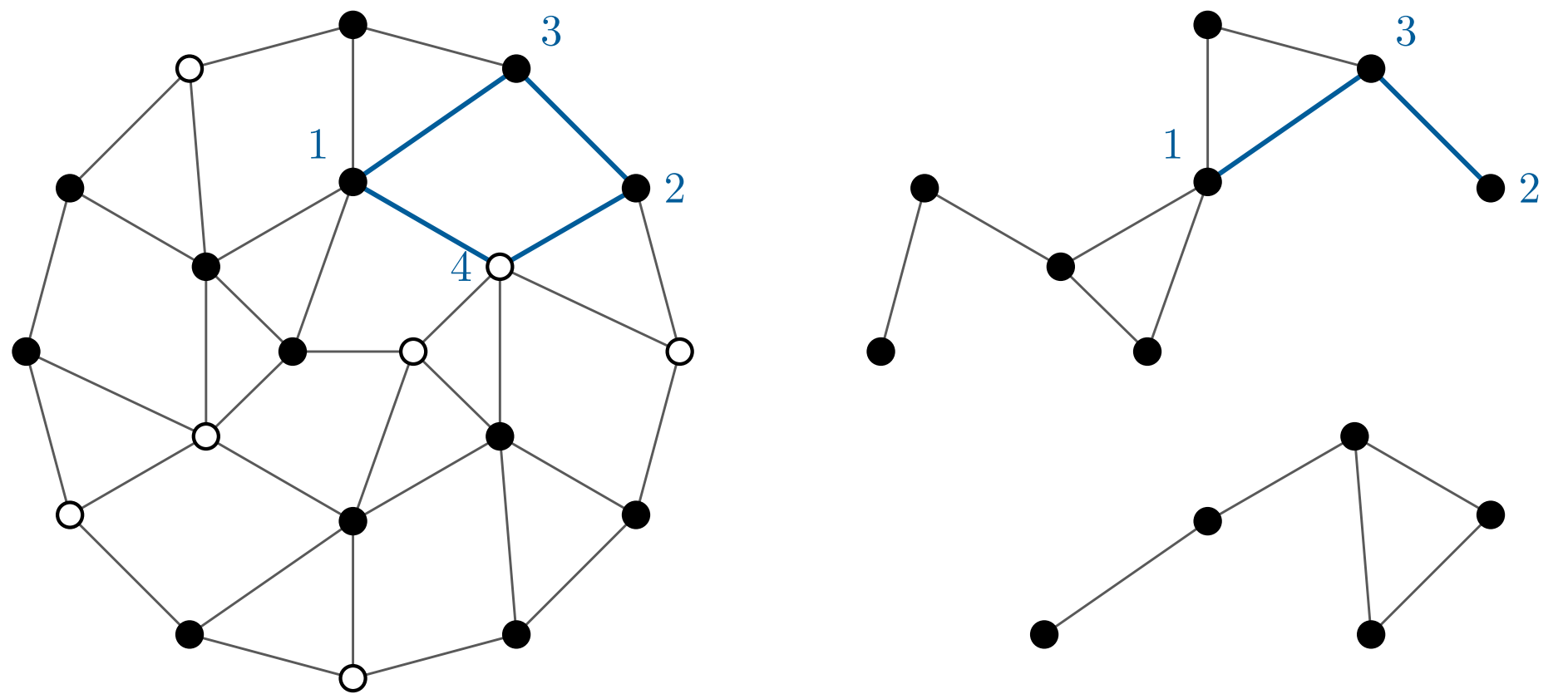

Let $X_1,\ldots,X_n$ be independent $\operatorname{Ber}(p_n)$ variables.  The
sampled graph is $\tld G_n=G_n[\{i:X_i=1\}]$, and we write $N_n:=|\tld V_n|\sim\Bin(n,p_n)$, where $\tld V_n:=V(\tld G_n)$, and denote the edge density and $r$-codegree measure of $\tld G_n$ by
\begin{align}
    \tld q_n:=\frac{1}{N_n^2}\sum_{i,j\in V(\tld G_n)} G_n(i,j),\quad
    \tld \nu_{n}:=\frac{1}{\binom{N_n}r}\sum_{S\in\binom{\tld V_n}r}\delta_{\tld d_S/(N_n \tld q_n^r)},
    \label{eq:nutld}
\end{align}
where $\tld d_S:=\sum_{j=1}^n G_n(S,j)X_j$ for every nonempty $S\subseteq[n]$. When $S\subseteq\tld V_n$, this is the codegree of $S$ in $\tld G_n$. The edge density is defined by the displayed formula when $N_n>0$, and we set $\tld q_n=0$ when $N_n=0$. The displayed formula for $\tld\nu_n$ applies when $N_n\ge r$ and $\tld q_n>0$; otherwise, set $\tld\nu_n:=\delta_0$. For nonempty $S\subseteq[n]$ with $|S|\le r$, write
\begin{align}
    D_n(S):=\frac{d_S}{nq_n^{|S|}},\qquad
    \tld D_n(S):=\frac{\tld d_S}{N_n\tld q_n^{|S|}},
    \label{eq:normalized-discrete-codegrees}
\end{align}
with $\tld D_n(S):=0$ when its denominator vanishes. For singleton sets,
write $d_i:=d_{\{i\}}$, $\tld d_i:=\tld d_{\{i\}}$, $D_n(i):=D_n(\{i\})$, and $\tld D_n(i):=\tld D_n(\{i\})$.
For a bounded measurable function $f$, define
\begin{align*}
    T_n(f):=
    \sqrt{np_n}
    \left(\widetilde{\nu}_n(f)-\nu_n(f)\right),\text{ where }
    \tld\nu_n(f)
:=\int f\ d\tld\nu_{n},\text{ and }
    \nu_n(f):=\int f d\nu_{n}.
\end{align*}

We describe the asymptotic structure of the parent graphs
through their normalized adjacency kernels.

Let $\mathscr W$ be the space of integrable kernels
$W:[0,1]^2\to\mathbb R$ that are symmetric almost everywhere, and let
$\mathscr W_+$ be its subset of nonnegative kernels. Define the normalized
empirical graphon $W_n$ by
\begin{align}
    W_n(x,y):= \frac{G_n(\lceil nx \rceil,\lceil ny \rceil)}{q_n},\quad (x,y)\in(0,1]^2,
    \label{eq:normalized-graphon-definition}
\end{align}
and $W_n\equiv0$ if $q_n=0$. Analogously, if $V(\tld G_n)=\{i\in[n]:X_i=1\}=\{i_1<\dots<i_{N_n}\}$, then we denote by $\tld W_n$ the corresponding normalized empirical graphon of $\tld G_n$, given by
\begin{align}
    \tld W_n(x,y):=\frac{G_n(i_{\lceil N_n x \rceil},i_{\lceil N_n y \rceil})}{\tld q_n},\quad (x,y)\in(0,1]^2,
    \label{eq:normalized-sampled-graphon-definition}
\end{align}
and $\tld W_n\equiv0$ if $\tld q_n=0$.

For $U\in\mathscr W$ and $1\le t\le r$, define the
$t$-star kernel and, wherever the integral is finite, the $t$-codegree function by
\begin{align}
    K_U^{(t)}(\mathbf x;y)
    &:=
    \prod_{a=1}^t U(x_a,y),\quad
    D_U^{(t)}(\mathbf x)
    :=
    \int_0^1 K_U^{(t)}(\mathbf x;y)\,dy,
    \label{eq:codegree-kernel-definition}
\end{align}
for $\mathbf x=(x_1,\ldots,x_t)\in[0,1]^t$
and $y\in[0,1]$. In particular, if $t=1$, $D_W:=D_W^{(1)}$ is the degree function of $W$. Also, for $t\in \bN$ and $F\in L^1([0,1]^t\times[0,1])$, define the \emph{cut norm}
of $F$, viewed as a kernel between $[0,1]^t$ and $[0,1]$, by
\begin{align}
    \|F\|_{\square,t}
    &:=
    \sup_{A\subseteq[0,1]^t,\; B\subseteq[0,1]}
    \left|
        \int_{A\times B}
        F(\mathbf x;y)\,d\mathbf x\,dy
    \right|.
    \label{eq:cut-norm-definition}
\end{align}
For $t=1$, we write $\|F\|_\square:=\|F\|_{\square,1}$,
the cut norm on $[0,1]^2$.

For a finite set $A$, write
\begin{align}
    [n]^A:=\{s:A\to[n]\},
    \qquad
    [n]_A:=\{s\in[n]^A:s\text{ is injective}\},
    \label{eq:function-notation}
\end{align}
and denote the image of $s$ by $\bar s:=s(A)$.  If $A=[m]$, we also write $[n]_m$ and
$(n)_m=n(n-1)\cdots(n-m+1)=m!\binom nm$.

Also, given the intervals $I_i:=((i-1)/n,i/n],\ i\in [n]$, let
\begin{align*}
    \Delta_{n,t}
    &:=
    \bigcup_{s=(s(1),\ldots,s(t))\in[n]_t}
    I_{s(1)}\times\cdots\times I_{s(t)},
\end{align*}
the region corresponding to injective $t$-tuples. In particular, for every $s\in [n]_t$
\begin{align*}
    D_{W_n}^{(t)}(\mathbf x)
    =
    \frac{d_{\bar s}}{nq_n^t} = \frac{1}{nq_n^t}\sum_{j=1}^n \prod_{i=1}^t G_n(s(i),j),
    \qquad
    \mathbf x\in I_{s(1)}\times\cdots\times I_{s(t)}.
\end{align*}

\subsection{Assumptions}\label{subsec:assumptions}

\begin{assumption}\label{ass:weights}
The graph sequence $(G_n)$ has uniformly bounded symmetric weights and no self-loops: for every $n$ and
all $i,j\in[n]$,
\begin{align*}
    0\le G_n(i,j)=G_n(j,i)\le1,\text{ and }G_n(i,i)=0.
\end{align*}
\end{assumption}

\begin{assumption}\label{ass:star-cut-convergence}
There exists $W\in\mathscr W_+$ such that, up to permutations of the vertices of $G_n$,
\begin{align*}
    \max_{1\le t\le r}
    \left\|K_{W_n}^{(t)}-K_W^{(t)}\right\|_{\square,t}
    \overset{n\to\infty}{\longrightarrow} 0.
\end{align*}
\end{assumption}

\begin{assumption}\label{ass:bounded-codegrees}
The graph sequence $(G_n)$ has uniformly bounded codegrees: there
exists $M<\infty$ such that
\begin{align*}
    \sup_n \max_{1\le t\le r}
    \left\|
    D_{W_n}^{(t)}\ind_{\Delta_{n,t}}
    \right\|_\infty=\sup_n \max_{1\le t\le r} \max_{S\in\binom{[n]}{t}} \frac{d_S}{nq_n^t}
    < M.
\end{align*}
In particular, $\nu_n$ is supported in $[0,M]$, and $M> \int_0^1 D_{W_n}(x)\,dx=1$.
\end{assumption}

\begin{assumption}\label{ass:sampling-ratios}
The sampling ratios $(p_n)$ satisfy $p_n\overset{n\to\infty}{\longrightarrow} p_\infty\in[0,1)$.
\end{assumption}

\begin{assumption}\label{ass:sparsity}
We have $np_n\to\infty$ and
\begin{align*}
    0\le\alpha_*<1, \text{ where } \alpha_*:=\frac{r\gamma_*}{1-r\gamma_*},\quad \text{ with }\gamma_*:=
    \limsup_{n\to\infty}
    \frac{\log(1/q_n)}{\log(np_n)}.
\end{align*}
\end{assumption}

Under \cref{ass:star-cut-convergence,ass:bounded-codegrees}, the limiting
kernel satisfies $\int_{[0,1]^2}W=1$ and
$0\le D_W^{(t)}\le M$ almost everywhere for $1\le t\le r$.
Indeed, integrating a star-kernel difference over its last coordinate gives
\begin{align*}
    \|D_{W_n}^{(t)}-D_W^{(t)}\|_1
    \le 2\|K_{W_n}^{(t)}-K_W^{(t)}\|_{\square,t}\longrightarrow0,
\end{align*}
where $\|\cdot\|_1$ denotes the $L^1$ norm. The uniform bound on $\Delta_{n,t}$ and
$|[0,1]^t\setminus\Delta_{n,t}|\to0$ then give the limiting bound.
The normalization follows from the case $t=1$, since $\int W_n=1$.

\subsection{Graphon-generated parent sequences}
\label{subsec:graphon-examples}

The structural assumptions in the previous subsection are formulated for deterministic parent graph sequences. The next result verifies them almost surely for graphon-generated
networks, after which the sampling results apply conditionally on the
realized parent sequence.

\begin{definition}[Random graphon model]\label{def:random-graphon}
Let $W\in\mathscr W_+$ satisfy $\int_{[0,1]^2}W=1$, and let
$0<\rho_n\le1$ be a sequence of target edge densities. Let
$U_1,U_2,\ldots$ be independent uniform random variables on $[0,1]$.
For each $n$, conditionally on these latent variables, let
$(A_{ij}^{(n)})_{1\le i<j\le n}$ be independent Bernoulli variables with
\begin{align*}
    \bP(A_{ij}^{(n)}=1\mid U_1,U_2,\ldots)
    =1\wedge\rho_nW(U_i,U_j).
\end{align*}
Set $A_{ji}^{(n)}=A_{ij}^{(n)}$, $A_{ii}^{(n)}=0$, and let $G_n$ be
the corresponding simple graph. No independence between different values
of $n$ is required. Let $\sigma_n$ order the latent variables increasingly.
In this example, use the corresponding relabeling in the definition of
$W_n$, so that
\begin{align}
    q_n=\frac1{n^2}\sum_{i,j=1}^n A_{ij}^{(n)},\qquad
    W_n(x,y)=\frac{A_{\sigma_n(i)\sigma_n(j)}^{(n)}}{q_n},
    \quad (x,y)\in I_i\times I_j,
    \label{eq:W_n-q-definition}
\end{align}
with $W_n\equiv0$ if $q_n=0$.
\end{definition}

\begin{theorem}
\label[theorem]{thm:sparse-graphon-samples}
Fix $r\ge1$ and consider the random graphon model defined above. Suppose that $n\rho_n^r/\log n\to\infty$, $\max_{1\le t\le r}\|D_W^{(t)}\|_\infty<\infty$, and that either
\begin{enumerate}[(i)]
    \item $\|W\|_\infty<\infty$ and
    $(\limsup_{n\to\infty}\rho_n)\|W\|_\infty\le1$; or
    \item $\rho_n\to0$.
\end{enumerate}
Then, almost surely,
\begin{align*}
    q_n=\left(1+o(1)\right)\rho_n,
    \quad
    \max_{1\le t\le r}
    \|K_{W_n}^{(t)}-K_W^{(t)}\|_{\square,t}\overset{n\to\infty}{\longrightarrow} 0.
\end{align*}
Moreover, almost surely, for every $1\le t\le r$ and sufficiently
large $n$,
\begin{align*}
    \|D_{W_n}^{(t)}\ind_{\Delta_{n,t}}\|_\infty
    \le4\bigl(1+\|D_W^{(t)}\|_\infty\bigr).
\end{align*}
\end{theorem}

\begin{remark}\label[remark]{rem:check}
Thus \cref{ass:weights,ass:star-cut-convergence,ass:bounded-codegrees}
hold almost surely (provided $G_n$ is not empty), 
for a suitable $M$ larger than the displayed bound $4\bigl(1+\|D_W^{(t)}\|_\infty\bigr)$. If, in addition,
$p_n\to p_\infty\in[0,1)$, $np_n\to\infty$, and
\begin{align*}
    \limsup_{n\to\infty}
    \frac{\log(1/\rho_n)}{\log(np_n)}
    <\frac{1}{2r},
\end{align*}
then \cref{ass:sampling-ratios,ass:sparsity} also hold almost surely.
Indeed, $q_n/\rho_n\to1$ implies that the limsup above equals
$\gamma_*$ almost surely, and $\gamma_*<1/(2r)$ is equivalent to $\alpha_*<1$.
\end{remark}

\subsection{Distributional results}
\label{subsec:distributional-results}
For a positive integer $k$ and $\alpha\in(0,1]$, let
$C^{k,\alpha}([0,M])$ denote the set of $k$ times continuously differentiable
functions on $[0,M]$ whose $k$th derivative is $\alpha$-H\"older, equipped with
the norm
\begin{align}
    \|f\|_{C^{k,\alpha}([0,M])}
    :=
    \sum_{i=0}^k \|f^{(i)}\|_\infty
    +
    [f^{(k)}]_\alpha,\qquad [f^{(k)}]_\alpha:=\sup_{x\ne y} \frac{|f^{(k)}(x)-f^{(k)}(y)|}{|x-y|^\alpha}.
    \label{eq:holder-norms}
\end{align}
We use the convention
\begin{align*}
    C^{k,0}([0,M]):=C^k([0,M]),
    \qquad
    \|f\|_{C^{k,0}}:=\sum_{i=0}^k\|f^{(i)}\|_\infty .
\end{align*}
Write
\begin{align*}
    \mathcal B_{k,\alpha}
    :=
    \left\{
        f\in C^{k,\alpha}([0,M]):
        \|f\|_{C^{k,\alpha}([0,M])}\le1
    \right\}.
\end{align*}

When a test function $f\in C^1([0,M])$ is evaluated at a sampled
codegree outside $[0,M]$, we use the bounded extension
$f(x):=f(x\wedge M)$, $x\ge0$. In formulas involving $f'$, extend
$f'$ separately by $f'(x):=f'(x\wedge M)$; this is not the derivative
of the bounded extension beyond $M$.

For $k\ge1$, let $K_{r,k}$ be the complete bipartite graph with
vertex set $[r+k]$ and parts
$L(K_{r,k})=[r]$ and $R(K_{r,k})=\{r+1,\ldots,r+k\}$.
Write $K_1$ for a singleton and $K_{1,1}=K_2$ for an edge.
For a finite simple graph $H$ and a weighted graph $G$, let
\begin{align*}
    \InjHom(H,G)
    :=\sum_{s\in V(G)_{V(H)}}
    \prod_{uv\in E(H)}G(s(u),s(v)),
\end{align*}
where $V(G)_{V(H)}$ is the set of injective maps from $V(H)$ to $V(G)$.
Define the normalized injective densities by
\begin{align}
    t_{\mathrm{inj}}(H,G_n)
    &:=\frac{\InjHom(H,G_n)}{(n)_{|V(H)|}q_n^{|E(H)|}},&
    t_{\mathrm{inj}}(H,\tld G_n)
    &:=\frac{\InjHom(H,\tld G_n)}{(N_n)_{|V(H)|}\tld q_n^{|E(H)|}}.
    \label{eq:inj-density-definition}
\end{align}
A normalized injective density is defined to be zero when its denominator
vanishes.

The quantitative approximation uses centered counts with deterministic
normalization:
\begin{align}
    Y_n(H):=
    \frac{\InjHom(H,\tld G_n)-p_n^{|V(H)|}\InjHom(H,G_n)}
    {(np_n)^{|V(H)|-1/2}q_n^{|E(H)|}}.
    \label{eq:centered-counts}
\end{align}
For fixed $m\ge1$, let
$\mathcal H=\mathcal H_m=(K_{r,1},\ldots,K_{r,m},K_1,K_{1,1})$,
and set $\mathbf Y_n=(Y_n(H))_{H\in\mathcal H}$ and
$\Sigma_n=\Var(\mathbf Y_n)$. Here $\mathcal H$ is an ordered family:
when $r=1$, the two occurrences of $K_{1,1}$ are kept as separate,
identical coordinates. Singular covariance matrices are therefore allowed. The coordinates corresponding to $K_1$ and $K_{1,1}$
account for the sampled vertex count and edge density
in the subsequent random normalization.

For random vectors $X,Y\in\mathbb R^d$ with finite first moments, define
\begin{align}
    d_3(X,Y):=
    \sup_{\varphi\in\mathcal F_3(\mathbb R^d)}
    |\E\varphi(X)-\E\varphi(Y)|,
    \label{eq:smooth-distance-definition}
\end{align}
where $\mathcal F_3(\mathbb R^d)$ consists of the functions
$\varphi\in C^3(\mathbb R^d)$ satisfying
\begin{align*}
    \sup_x\|\nabla\varphi(x)\|_2\le1,\qquad
    \sup_x\|\nabla^2\varphi(x)\|_{\mathrm{op}}\le1,\qquad
    \operatorname{Lip}(\nabla^2\varphi)\le1.
\end{align*}
The Lipschitz seminorm uses the operator norm on matrices and the Euclidean
norm on $\mathbb R^d$. We also write $d_3(X,\mathcal N(0,\Sigma))$
when the second argument is a Gaussian law.

The first distributional result is a normal approximation for $\mathbf Y_n$ with its finite-sample covariance. It requires neither a limiting graphon nor convergence of the sampling ratios.

\begin{theorem}
\label[theorem]{thm:stein-error-bound}
Fix $m\ge1$. Suppose that $(G_n)$ satisfies
\cref{ass:weights,ass:bounded-codegrees}, that $np_n\to\infty$, and that
$np_nq_n^r\ge c_0$ for all sufficiently large $n$, for some constant $c_0>0$.
Then
\begin{align}
    d_3(\mathbf Y_n,\mathcal N(0,\Sigma_n))
    \lesssim_{r,m,M,c_0}(np_n)^{-1/2}.
    \label{eq:stein-error-rate}
\end{align}
\end{theorem}

To specify the limiting covariance, let $\mathbf U=(U_1,\ldots,U_r)$ have independent uniform coordinates on $[0,1]$. For $f\in C^1([0,M])$, define the vertex influence function
$\phi_f:[0,1]\to\mathbb R$ by
\begin{align}
\phi_f(x)
&:=r\left\{
\E[f(D_W^{(r)}(\mathbf U))\mid U_1=x]
-\E[f(D_W^{(r)}(\mathbf U))]
\right\}
\notag\\
&\quad+
\E\left[
 f'(D_W^{(r)}(\mathbf U))
 \left\{
 \prod_{a=1}^rW(U_a,x)-D_W^{(r)}(\mathbf U)
 \right\}
\right]
\notag\\
&\quad-2r\,\E[D_W^{(r)}(\mathbf U)f'(D_W^{(r)}(\mathbf U))]
\{D_W(x)-1\}.
\label{eq:limiting-influence-function}
\end{align}
Each term is centered in $x$, so $\int_0^1\phi_f(x)\,dx=0$. Also, define
\begin{align}
\Gamma(f,g):=(1-p_\infty)
\int_0^1\phi_f(x)\phi_g(x)\,dx,
\qquad f,g\in C^1([0,M]).
\label{eq:limiting-covariance-kernel}
\end{align}
These quantities are finite. Indeed,
\begin{align*}
\int_{[0,1]^r}\prod_{a=1}^rW(u_a,x)\,d\mathbf u=D_W(x)^r,
\end{align*}
so the bounds following the assumptions imply that $\phi_f$ is bounded.

Under the graph-limit assumptions, covariance convergence and adjustment
for the random vertex and edge-density normalizations yield the following
corollary.

\begin{corollary}
\label[corollary]{cor:biclique-clt}
Fix $m\ge1$. Suppose that $(G_n)$ satisfies
\cref{ass:weights,ass:star-cut-convergence,ass:bounded-codegrees}, that
$(p_n)$ satisfies \cref{ass:sampling-ratios}, and that $np_nq_n^r\to\infty$.
Then
\begin{align*}
    \sqrt{np_n}\left(
        t_{\mathrm{inj}}(K_{r,k},\tld G_n)
        -t_{\mathrm{inj}}(K_{r,k},G_n)
    \right)_{1\le k\le m}
    \dto\mathcal N(0,\Gamma_m),
\end{align*}
where $\Gamma_m=(\Gamma(x^k,x^\ell))_{1\le k,\ell\le m}$.
\end{corollary}

{
\begin{definition}\label[definition]{def:Gf_brownian_int}
Let $B$ be a standard Brownian motion on $[0,1]$. Define the centered
Gaussian process $G=\{G(f):f\in C^1([0,M])\}$ by
\begin{align*}
G(f):=\sqrt{1-p_\infty}\int_0^1\phi_f(x)\,dB(x).
\end{align*}
By the It\^o isometry, $G$ has covariance kernel $\Gamma$
defined in \eqref{eq:limiting-covariance-kernel}.
\end{definition}
}
The joint convergence of sampled biclique densities extends
to smooth functions of the empirical codegree distribution, as stated in the following theorem.

\begin{theorem}\label[theorem]{thm:empirical-codegree-process}
Suppose that the graph sequence $(G_n)$ satisfies
\cref{ass:weights,ass:star-cut-convergence,ass:bounded-codegrees}, the sampling
ratios $(p_n)$ satisfy \cref{ass:sampling-ratios}, and the two sequences jointly satisfy \cref{ass:sparsity}. Then, for every
$\alpha_*<\alpha\le1$,
\begin{align*}
    T_n\dto G,
\end{align*}
in the sense of finite-dimensional distributions on
$C^{1,\alpha}([0,M])$, where $G$ is the process from \cref{def:Gf_brownian_int}.
Moreover, the convergence holds in $\ell^\infty(\mathcal B_{1,\alpha})$ for every $\alpha>\max\{\alpha_*,1/2\}$.
If $\liminf_{n\to\infty}q_n>0$ (in which case $\gamma_*=\alpha_*=0$), then the finite-dimensional convergence
holds more generally on $C^1([0,M])$.
\end{theorem}

\subsection{Examples}\label{subsec:worked-examples}

We illustrate the limiting process through the following examples.
In each example, generate the parent graphs as in
Definition \ref{def:random-graphon}, and assume $p_n\to0$,
$np_n\to\infty$, and
$\limsup_{n\to\infty}\frac{\log(1/\rho_n)}{\log(np_n)}<\frac1{2r}$.
Together with the bounds and restrictions specified below,
\cref{thm:sparse-graphon-samples} and \cref{rem:check}
ensure that \cref{thm:empirical-codegree-process} applies
conditionally on almost every realized parent sequence.

\begin{example}[Degeneracy of the degree process]
For $r=1$, the limiting process satisfies
\begin{align*}
G\equiv0 \iff D_W\equiv1.
\end{align*}
By \eqref{eq:limiting-influence-function},
$D_W\equiv1$ implies $\phi_f\equiv0$ for every $f$. Conversely, writing $d=D_W$ and using the symmetry
of $W$, we obtain
\begin{align*}
\int_0^1 \left(3d(x)^2\phi_{z^2}(x)
-2d(x)\phi_{z^3}(x)\right)\ dx
=\int_0^1(d(x)-1)^4\ dx
\quad+3\left(\int_0^1(d(x)-1)^2\,dx\right)^2.
\end{align*}
If $G\equiv0$, the left-hand side vanishes, forcing $d\equiv1$.
\end{example}

\begin{example}[Stochastic block model]
\label{ex:two-block-models}
Consider a balanced two-block model with
$A_1=[0,1/2]$, $A_2=(1/2,1]$, and
\begin{align*}
W(x,y):=\frac1q
\begin{cases}
a,&(x,y)\in A_1^2,\\
b,&(x,y)\in A_2^2,\\
c,&\text{otherwise},
\end{cases}
\qquad q:=\frac{a+b+2c}{4},
\end{align*}
where $a,b,c\in[0,1]$ and $a+b+c>0$. Assume
$(\limsup_n\rho_n)\max\{a,b,c\}\le q$.

For a tuple $\mathbf x\in[0,1]^r$ with
$t$ coordinates in $A_1$, splitting the defining
integral over $A_1$ and $A_2$ gives
\begin{align}
D_W^{(r)}(\mathbf x)
=\frac{a^tc^{r-t}+c^tb^{r-t}}{2q^r}.
\label{eq:two-block-codegree}
\end{align}

If $a=b$, we show that $\phi_f\equiv0$ for every
$f\in C^1([0,M])$.
Let $\sigma(x):=1-x$, $x\in[0,1]$, which interchanges
$A_1$ and $A_2$ except at $1/2$.
Then $W(\sigma(x),\sigma(y))=W(x,y)$, and \eqref{eq:two-block-codegree} gives
$D_W^{(r)}(\sigma\mathbf x)=D_W^{(r)}(\mathbf x)$
, where
$\sigma\mathbf x:=(\sigma(x_1),\dots,\sigma(x_r))$.
Since $\sigma$ preserves Lebesgue measure,
$\sigma\mathbf U\overset{d}=\mathbf U$ for
$\mathbf U=(U_1,\dots,U_r)\overset{\text{i.i.d.}}{\sim}\operatorname{Unif}(0,1)$.
Hence,
\begin{align*}
\E[f(D_W^{(r)}(\mathbf U))\mid U_1=\sigma(x)]
&=\E\left[f\left(D_W^{(r)}(\sigma(x),\sigma(U_2),\ldots,\sigma(U_r))\right)\right]\\
&=\E[f(D_W^{(r)}(\mathbf U))\mid U_1=x],
\end{align*}
and
\begin{align*}
\E\left[f'(D_W^{(r)}(\mathbf U))
\prod_{j=1}^r W(U_j,\sigma(x))\right]
&=\E\left[f'\left(D_W^{(r)}(\sigma\mathbf U)\right)
\prod_{j=1}^r W(\sigma(U_j),\sigma(x))\right]\\
&=\E\left[f'\left(D_W^{(r)}(\mathbf U)\right)
\prod_{j=1}^r W(U_j,x)\right].
\end{align*}
Since $D_W\equiv1$, \eqref{eq:limiting-influence-function} gives $\phi_f(\sigma(x))=\phi_f(x)$ for $x\ne1/2$.
The same formula shows that $\phi_f$ is constant on each $A_i$, since $W(x,\cdot)$ is.
Thus $\phi_f$ is constant on $[0,1]$, and
$\int_0^1\phi_f(x)\,dx=0$ yields $\phi_f\equiv0$.

If $a>b$, \eqref{eq:two-block-codegree} attains
its maximum $(a^r+c^r)/(2q^r)$ precisely when
$t=r$, that is, on $A_1^r$.
Choose a smooth $f$ equal to $1$ near this maximum
and $0$ near all other codegree values.
Since $f'$ vanishes at every codegree value,
\eqref{eq:limiting-influence-function} yields
\begin{align*}
\phi_f
=\frac{r}{2^r}(\ind_{A_1}-\ind_{A_2}),
\qquad
\Gamma(f,f)=(1-p_\infty)\frac{r^2}{4^r}>0.
\end{align*}
The case $b>a$ is symmetric. Therefore,
\begin{align*}
G\equiv0
&\iff a=b
&& (r\ge1),\\
D_W^{(r)}\text{ is constant}
&\iff a=b=c
&& (r\ge2).
\end{align*}
The second equivalence follows from
\eqref{eq:two-block-codegree}: equality of the values at $t=0,r$ gives $a=b$, while equality at $t=1,r$ then gives $a=c$.
In particular, when $a=b\ne c$, the codegree
function is nonconstant for every $r\ge2$,
although the limiting process is degenerate.
\end{example}

\begin{example}[Rank-one graphon]
\label[example]{ex:rank-one-graphon}
Let $h\in L^\infty([0,1])$ be nonnegative with
$\int_0^1h(u)\,du=1$, and define
\begin{align*}
W(x,y):=h(x)h(y).
\end{align*}
Then $\int W=1$, and by \eqref{eq:codegree-kernel-definition}
\begin{align}
D_W(x)=h(x),
\qquad
D_W^{(t)}(x_1,\ldots,x_t)
=m_t\prod_{j=1}^t h(x_j),
\qquad
m_t:=\int_0^1 h(u)^t\ du.
\label{eq:one-rank-moments}
\end{align}
In particular,
$\|D_W^{(t)}\|_\infty\le m_t\|h\|_\infty^t<\infty$
for every fixed integer $t\ge1$.
Assume that $(\limsup_n\rho_n)\|h\|_\infty^2\le1$.

For the degree process ($r=1$),
\eqref{eq:limiting-influence-function} simplifies to
\begin{align*}
\phi_f(x)
=f(h(x))-\int_0^1f(h(u))\ du
-(h(x)-1)\int_0^1h(u)f'(h(u))\ du.
\end{align*}
A rank-one graphon can give rise to a limiting covariance kernel
of infinite rank.
For example, take $h(x)=2x$, so that $W(x,y)=4xy$.
For integers $k\ge1$,
\begin{align*}
\phi_{z^k}(x)
=(2x)^k-\frac{2^k}{k+1}\{1+k(2x-1)\}.
\end{align*}
Thus $\phi_z\equiv0$, while
$\phi_{z^2},\ldots,\phi_{z^m}$ are polynomials of distinct
degrees and are linearly independent in $L^2([0,1])$
for every $m\ge2$.
By \eqref{eq:limiting-covariance-kernel}, the matrix
$\bigl(\Gamma(z^i,z^j)\bigr)_{2\le i,j\le m}$
is positive definite for every $m\ge2$.
Hence the limiting degree process has an infinite-rank
covariance kernel.

For general $h$ as above and $r=2$,
\eqref{eq:one-rank-moments} gives
$D_W^{(2)}(x,y)=m_2h(x)h(y)$.
Substitution in \eqref{eq:limiting-influence-function} yields
\begin{align*}
\phi_z(x)
&=h(x)^2-2m_2h(x)+m_2,\\
\phi_{z^2}(x)
&=4m_2^3\{h(x)^2-m_2\}
-8m_2^4\{h(x)-1\}
=4m_2^3\phi_z(x).
\end{align*}
Thus the test function
\begin{align*}
f_*(z)=z^2-4m_2^3z
\end{align*}
satisfies $\phi_{f_*}\equiv0$ for $r=2$.
For $r=1$, however, $\phi_z\equiv0$, so using linearity of $f\mapsto \phi_f$ gives
$\phi_{f_*}=\phi_{z^2}$.
In particular, for $h(x)=2x$,
\begin{align*}
\Gamma(f_*,f_*)>0 \quad\text{for }r=1,
\quad\text{but}\quad
\Gamma(f_*,f_*)=0 \quad\text{for }r=2.
\end{align*}
Thus positive limiting variance for a fixed test function
at $r=1$ need not persist at $r=2$.
\end{example}

\begin{example}[Unbounded graphon]
\label[example]{ex:power-graphon}
Let $0<\theta<1/r$ and, for $x\ne y$, define
\begin{align*}
W(x,y):=c_\theta|x-y|^{-\theta},
\quad\text{where }
c_\theta:=\frac{(1-\theta)(2-\theta)}2.
\end{align*}
Then $\int W=1$.
Assume $\rho_n\to0$.
Although $W$ is unbounded, its degree function
\begin{align*}
D_W(x)=\frac{2-\theta}{2}
\left[x^{1-\theta}+(1-x)^{1-\theta}\right]
\end{align*}
is bounded.
More generally, H\"older's inequality gives
\begin{align*}
\|D_W^{(t)}\|_\infty
\le \sup_{x\in[0,1]}\int_0^1W(x,y)^t\,dy
=\frac{2^{t\theta}c_\theta^t}{1-t\theta}<\infty,
\qquad 1\le t\le r.
\end{align*}

Put $m_r=\int_0^1D_W(u)^r\,du$.
For the identity test function, Fubini's theorem and
\eqref{eq:limiting-influence-function} give
\begin{align*}
\phi_z(x)
=D_W(x)^r
+rc_\theta\int_0^1|x-u|^{-\theta}D_W(u)^{r-1}\,du
+\{r-1-2rD_W(x)\}m_r.
\end{align*}
For $r=1$, we have $m_1=1$, and the displayed formula gives
$\phi_z\equiv0$.
At the same order, substitution of $f(z)=z^2$ gives
\begin{align*}
\phi_{z^2}(x)
=D_W(x)^2
+2c_\theta\int_0^1|x-u|^{-\theta}D_W(u)\,du
+\{1-4D_W(x)\}m_2,
\end{align*}
where $m_2=\int_0^1D_W(u)^2\,du$.
Thus, $\phi_{z^2}$ is not identically zero for $r=1$, and $\phi_z$ is not identically zero for $r\ge2$, and \eqref{eq:limiting-covariance-kernel} gives
\begin{align*}
\begin{aligned}
\Gamma(z,z)=0,\quad
\Gamma(z^2,z^2)=\int_0^1\phi_{z^2}(x)^2\,dx>0
&\qquad \text{if }r=1,\\
\Gamma(z,z)=\int_0^1\phi_z(x)^2\,dx>0
&\qquad \text{if }2\le r<1/\theta.
\end{aligned}
\end{align*}
Thus the limiting process is nondegenerate for every $r$ with
$r\theta<1$.
\end{example}

\subsection{Sharpness of the smoothness conditions}
\label{subsec:sharpness-results}

The two smoothness thresholds in
\cref{thm:empirical-codegree-process} reflect distinct phenomena. The condition
$\alpha>\alpha_*$ cannot in general be relaxed below $\alpha_*$, since sampling
bias may dominate the $\sqrt{np_n}$ scale and prevent finite-dimensional
convergence to the centered Gaussian limit; this is shown in
\cref{prop:regularity-counterexample}. By contrast, \cref{prop:degree-tightness-counterexample} shows that when
$\alpha<1/2$, convergence of finite-dimensional distributions may still hold while
tightness over the H\"older unit ball fails.
\begin{proposition}
\label[proposition]{prop:regularity-counterexample}
Let $r=1$ and $p_n\equiv1/2$. Fix $0<\gamma<1/2$. 
Then, there exists a sequence of simple unweighted graphs $G_n$ with edge density $q_n$ satisfying $n^{\gamma}q_n\to 1$  and limiting graphon $W\equiv1$, so that  \cref{ass:weights}-\cref{ass:sparsity} hold with $\gamma_*=\gamma, \alpha_*=\frac{\gamma}{1-\gamma}<1$. Moreover, for the choice $f(x):=|x-1|^{1+\alpha}$ with $0<\alpha<\alpha_*$ we have $T_n(f)\pto+\infty$. In particular, $T_n(f)$ does not converge in distribution to $G(f)$.
\end{proposition}

\begin{proposition}
\label[proposition]{prop:degree-tightness-counterexample}
Let $r=1$ and $p_n\equiv1/2$. Fix  $ \gamma\in(0,1/3)$ and $\alpha\in (\gamma/(1-\gamma),1/2)$. Then there exists a sequence of simple unweighted graphs $G_n$ with edge density $q_n$  satisfying $n^\gamma q_n\to1$  and a limiting graphon with a non-constant  degree function, so that
 \cref{ass:weights}-\cref{ass:sparsity} hold with $\gamma_*=\gamma, \alpha_*=\frac{\gamma}{1-\gamma}<\alpha$.
Moreover, the process
 $T_n$ converges in finite-dimensional distributions on
$C^{1,\alpha}([0,M])$ to $G$, but
\begin{align*}
    \sup_{f\in\mathcal B_{1,\alpha}}|T_n(f)|\pto+\infty.
\end{align*}
In particular, $T_n$ is not tight in
$\ell^\infty(\mathcal B_{1,\alpha})$, although finite-dimensional
convergence holds. The same conclusion holds for $T_n-\E T_n$.
\end{proposition}

\subsection{Overview of the proofs}
\label{subsec:proof-strategy}

The moments of the empirical codegree distribution are closely
related to biclique counts. For each fixed $k\ge1$,
\cref{lem:counting-lemma}\emph{(i)} shows that the $k$th
moment of the empirical $r$-codegree distribution can be
approximated by the normalized injective count of $K_{r,k}$.
When $r=1$, this reduces to the familiar relation between
degree moments and star counts. We therefore begin by
studying the joint fluctuations of sampled biclique counts.

Recall that the vertex sampling random variables $X_i=\ind\{i\in V(\tld G_n)\}$, $i\in[n]$,
are independent Bernoulli variables. The contribution of
a copy of $H$, represented by an injective map $s$, to its
sampled count contains the factor
$\prod_{i\in s(V(H))}X_i$. Hence contributions from
vertex-disjoint copies are independent. Using this
dependency structure and the multivariate Stein bound
of \citet{Raic2004Decomposable}, we prove the quantitative
normal approximation in \cref{thm:stein-error-bound}.
This follows the approach to motif estimation under vertex
sampling in \citet{bdm22}, but requires explicit bounds
on the error terms under our sparsity assumptions.
These terms involve mixed third moments, which reduce
to weighted counts of connected collections of overlapping
copies of bicliques, vertices, and edges.

One might try to control these counts using the
$L^p$ graphon counting lemma of
\citet[Theorem~2.20]{BCCZ2019LpI}.
However, this result requires uniform $L^p$ bounds
for some $p>\Delta(H)$, where $\Delta(H)$ denotes
the maximum degree of $H$. For an unweighted graph, the normalized adjacency kernel satisfies
$\|W_n\|_{L^p}=q_n^{1/p-1}$, where $q_n$ is the
edge density. Thus, this condition fails whenever
$q_n\to0$.
Other sparse counting lemmas, such as those of
\citet{ConlonFoxZhao2014}, rely on additional
pseudorandomness conditions, including jumbledness and relative regularity, which are not imposed here.

Instead, we use boundedness of normalized codegrees
through order $r$, a condition satisfied by a broad
class of deterministic and random graphs. For example,
in sparse Erd\H{o}s--R\'enyi graphs, the degrees
normalized by $nq_n$ are uniformly bounded with high
probability whenever $nq_n\gg\log n$. Analogous
bounds for codegrees through order $r$ hold when
$nq_n^r\gg\log n$. More generally, our
bounded-codegree assumption gives, for any fixed
$r$-degenerate graph $H$ (see
\cref{def:graph-degeneracy,lem:degen-count-bound}),
\begin{align*}
    \InjHom(H,G_n)
    \lesssim_H n^{|V(H)|}q_n^{|E(H)|}.
\end{align*}

The graphs formed by overlapping bicliques need not themselves be $r$-degenerate. By a delicate
combinatorial argument, in
\cref{lem:biclique-spanning-subgraph,lem:union-spanning-subgraph},
we show that they admit $r$-degenerate spanning
subgraphs retaining sufficiently many edges while preserving the connectedness of the join.
The construction extends the spanning-tree argument by choosing suitable degeneracy orderings for successive copies. This gives the estimates in \cref{lem:multijoin_bound}, which are strong enough to control the Stein error at rate $O((np_n)^{-1/2})$ whenever $np_nq_n^r$
is bounded away from zero.
For $r=1$, the counting estimates use only bounded
normalized degrees and give, for any fixed tree $H$,
the bound
\begin{align*}
    \InjHom(H,G_n)
    \lesssim_H n^{|V(H)|}q_n^{|V(H)|-1}.
\end{align*}
For $r=2$, the additional control of
pairwise codegrees gives, for example,
\begin{align*}
\InjHom(K_{2,k},G_n)\lesssim n^{k+2}q_n^{2k},
\end{align*}
whereas a spanning-tree argument using the $r=1$ bound above
 gives only $O(n^{k+2}q_n^{k+1})$.
The additional powers of $q_n$ are important for
controlling overlaps in sparse graphs.

The normal approximation in \cref{thm:stein-error-bound}
uses a finite-sample covariance matrix. To obtain the limiting Gaussian law, we identify its limit through
pairs of copies sharing exactly one vertex. Contributions
from further overlaps are negligible by
\cref{lem:multijoin_bound}, while convergence of the
corresponding one-vertex-join densities follows from
\cref{lem:counting-lemma}\emph{(iii)}. Accounting for
the random vertex count and edge density then gives the CLT in
\cref{cor:biclique-clt}. Combining this corollary along with \cref{lem:counting-lemma}, the convergence of finite-dimensional process marginals follows for polynomial test functions, using the moment comparison at the $\sqrt{np_n}$ scale. For CLT with general $C^{1,\alpha}$ functions (addressed in \cref{lem:fdd_convergence}), we use polynomial approximation in a weaker
H\"older norm, together with uniform variance and bias bounds
from \cref{lem:second_moment_bias_bound} for the approximation errors. The sparsity-dependent
restriction $\alpha>\alpha_*$ arises from the latter
estimates.

To establish functional convergence in
\cref{thm:empirical-codegree-process}, we prove tightness in \cref{prop:tightness_holder}. Flipping a single inclusion variable $X_i$ affects the vertex sets contributing to the empirical distribution, their codegrees, and the random normalizers. \Cref{lem:holder-process-hamming-stability} controls these effects on a high-probability event by bounding process increments in terms of the Hamming distance between sampling configurations.
Combining this estimate with Warnke's variance-sensitive
bounded-differences inequality, based on Freedman's martingale inequality, yields sub-Gaussian concentration for the expectation-centered process, including when the sampling rate $p_n$ vanishes (\cref{lem:holder-process-concentration}). A chaining argument using standard entropy bounds for
$\mathcal B_{1,\alpha}$ then establishes tightness of the expectation centered process when $\alpha>1/2$ (\cref{lem:centered-process-tightness}). Finally, uniform convergence of the bias to zero
shows that the expectation-centered process is
asymptotically equivalent to $T_n$, completing
the proof of \cref{thm:empirical-codegree-process}.

In \cref{thm:sparse-graphon-samples}, we demonstrate that our conditions \cref{ass:weights,ass:star-cut-convergence,ass:bounded-codegrees} hold when the parent graph $G_n$ is a sparse inhomogeneous random graph, produced by latent uniform random variables. The proof separates
the approximation of the conditional expectations of the
normalized star kernels, using the latent variables,
from their fluctuations due to edge sampling. The latter
are controlled in cut norm. Uniform bounds on normalized
codegrees require an additional argument, using concentration
inequalities simultaneously over all vertex tuples of
sizes up to $r$.

Finally, we give two counterexamples to justify the assumption $\alpha>\max(\alpha_*,1/2)$ in \cref{thm:empirical-codegree-process}, by identifying two distinct obstacles to functional
convergence. For $\alpha$ below the sparsity-dependent threshold $\alpha_*$, sampling bias can prevent Gaussian convergence even for a fixed test function (see Proposition \ref{prop:regularity-counterexample}). On the other hand, finite-dimensional convergence may hold
while process tightness breaks in the regime $\alpha_*<\alpha<1/2$ (see Proposition \ref{prop:degree-tightness-counterexample}).

\subsection{Discussion and future directions}
\label{sec:future-directions}

We discuss three directions for future work.

{\bf Simultaneous inference.}
The covariance formula \eqref{eq:limiting-covariance-kernel} suggests
estimating each vertex's contribution from the sampled graph. These
estimates could be used to construct a Gaussian multiplier bootstrap
for $\sup_{f\in\mathcal B_{1,\alpha}}|T_n(f)|$ and obtain simultaneous
confidence bounds for $\nu_n(f)$, as in related work on $U$-processes
\citep{ChenKato2020}.

{\bf Higher-order fluctuations.}
The examples in \cref{ex:two-block-models,ex:rank-one-graphon}
show that the limiting process, or individual marginals,
may be degenerate even when the codegree function is nonconstant.
For subgraph counts in graphon-generated random graphs,
\citet{BhattacharyaChatterjeeJanson2023,ChatterjeeDanBhattacharya2024}
obtained higher-order limits with potentially non-Gaussian
components. An analogous investigation of the empirical
codegree process under network sampling would require analyzing higher-order sampling terms, possible residual first-order
fluctuations, and sampling bias. Which terms contribute at
leading order may depend on the sampling rate and the
finite parent graph sequence.

{\bf Unequal sampling probabilities.}
\citet{HsiehHsuKoKovarikLogan2024} use sampling weights to estimate
the average degree and the global clustering coefficient when sampling
rates vary across vertex types. A natural
question is whether our functional limit theorem extends to a weighted
codegree distribution under independent sampling with unequal known
probabilities. \citet{BoistardLopuhaaRuizGazen2017} prove related
limits for weighted empirical distributions of directly observed
values. Here both the codegrees and the edge density used to normalize
them must be estimated from the sampled graph.

\section{Biclique counts and vertex joins}
\label{sec:biclique_counts}

This section proves \cref{thm:stein-error-bound,cor:biclique-clt}.
We first establish the deterministic kernel and overlap estimates needed
for the normal approximation and covariance identification. The counting
limits use $nq_n^r\to\infty$, whereas the connected-copy bound itself does
not require this convergence.

\subsection{Bicliques and one-vertex joins}

Recall the graphs $K_{r,k}$, $K_1$, and $K_{1,1}$ and the parts
$L(K_{r,k})$, $R(K_{r,k})$ from \cref{subsec:distributional-results}.
When the biclique is clear from context, write simply $L$ and $R$.

If $H$ and $H'$ are finite simple graphs and $u\in V(H)$, $v\in V(H')$, their
\emph{one-vertex join} is denoted by $H\oplus_{u,v}H'$, which is the graph obtained from $H$ and $H'$ by identifying
$u$ with $v$.  Thus $|V(H\oplus_{u,v}H')|=|V(H)|+|V(H')|-1$.
For example, if
$a\in L(K_{r,k})$ and $b\in R(K_{r,\ell})$, then
$K_{r,k}\oplus_{a,b}K_{r,\ell}$ is the pair of bicliques obtained by
identifying the left vertex $a$ of the first copy with the right vertex $b$ of
the second copy.

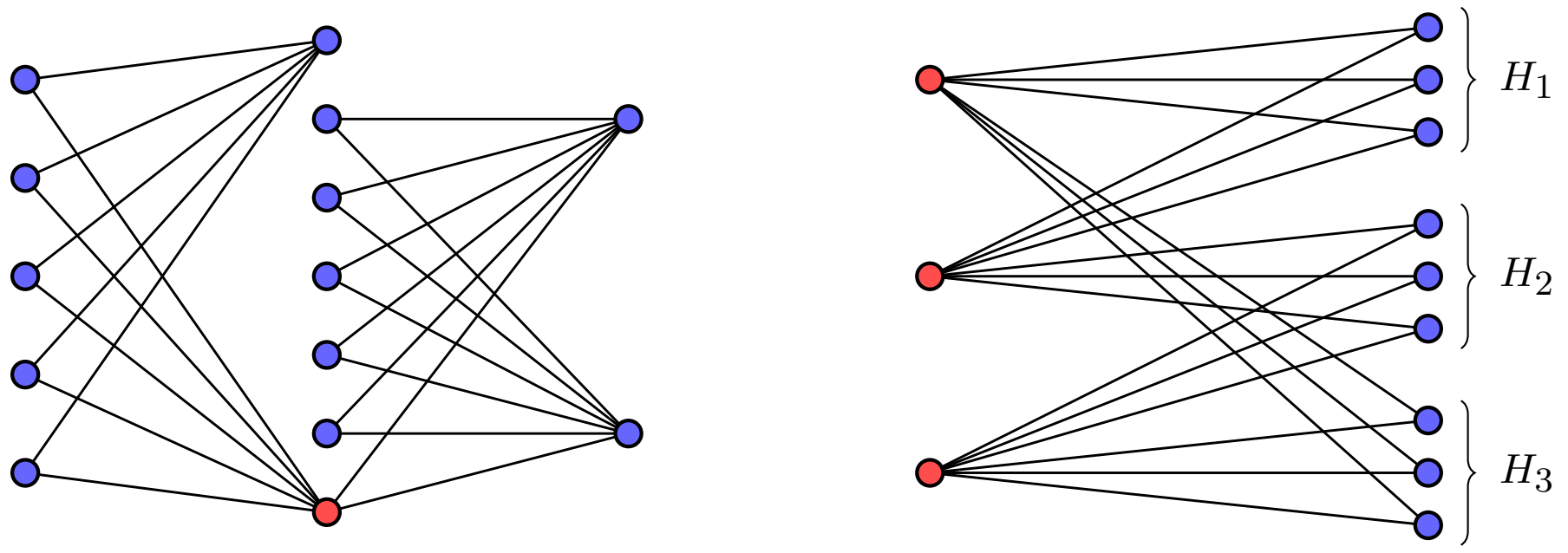
\begin{figure}[H]
\centering
\begin{tikzpicture}[
 x=1cm,y=1cm,font=\small,
 edge/.style={draw=black,line width=.55pt},
 vertex/.style={circle,draw=black,thick,inner sep=2pt,
   outer sep=0pt,fill=blue!60},
 shared/.style={vertex,fill=red!70},
 group/.style={decorate,decoration={brace,amplitude=3pt},thin}
]

\begin{scope}[shift={(0,0)}]
 \coordinate (s) at (0,-1.8);
 \coordinate (a) at (0,1.8);
 \coordinate (b) at (2.3,1.2);
 \coordinate (c) at (2.3,-1.2);
 \coordinate (l0) at (-2.3,1.5);
 \coordinate (l1) at (-2.3,0.75);
 \coordinate (l2) at (-2.3,0);
 \coordinate (l3) at (-2.3,-0.75);
 \coordinate (l4) at (-2.3,-1.5);
 \coordinate (r0) at (0,1.2);
 \coordinate (r1) at (0,0.6);
 \coordinate (r2) at (0,0);
 \coordinate (r3) at (0,-0.6);
 \coordinate (r4) at (0,-1.2);
 \draw[edge] (a)--(l0);
 \draw[edge] (a)--(l1);
 \draw[edge] (a)--(l2);
 \draw[edge] (a)--(l3);
 \draw[edge] (a)--(l4);
 \draw[edge] (b)--(r0);
 \draw[edge] (b)--(r1);
 \draw[edge] (b)--(r2);
 \draw[edge] (b)--(r3);
 \draw[edge] (b)--(r4);
 \draw[edge] (b)--(s);
 \draw[edge] (c)--(r0);
 \draw[edge] (c)--(r1);
 \draw[edge] (c)--(r2);
 \draw[edge] (c)--(r3);
 \draw[edge] (c)--(r4);
 \draw[edge] (c)--(s);
 \draw[edge] (l0)--(s);
 \draw[edge] (l1)--(s);
 \draw[edge] (l2)--(s);
 \draw[edge] (l3)--(s);
 \draw[edge] (l4)--(s);
 \node[shared] at (s) {};
 \node[vertex] at (a) {};
 \node[vertex] at (b) {};
 \node[vertex] at (c) {};
 \node[vertex] at (l0) {};
 \node[vertex] at (l1) {};
 \node[vertex] at (l2) {};
 \node[vertex] at (l3) {};
 \node[vertex] at (l4) {};
 \node[vertex] at (r0) {};
 \node[vertex] at (r1) {};
 \node[vertex] at (r2) {};
 \node[vertex] at (r3) {};
 \node[vertex] at (r4) {};
 \node[align=center,text width=5.8cm,anchor=north] at (0,-2.4) {(a) $K_{2,5}$ and $K_{2,6}$\\one shared left--right vertex};
\end{scope}
\begin{scope}[shift={(6.7,0)}]
 \coordinate (a) at (-2.1,1.5);
 \coordinate (b) at (-2.1,0);
 \coordinate (c) at (-2.1,-1.5);
 \coordinate (r00) at (1.7,1.9);
 \coordinate (r01) at (1.7,1.5);
 \coordinate (r02) at (1.7,1.1);
 \coordinate (r10) at (1.7,0.4);
 \coordinate (r11) at (1.7,0);
 \coordinate (r12) at (1.7,-0.4);
 \coordinate (r20) at (1.7,-1.1);
 \coordinate (r21) at (1.7,-1.5);
 \coordinate (r22) at (1.7,-1.9);
 \draw[edge] (a)--(r00);
 \draw[edge] (a)--(r01);
 \draw[edge] (a)--(r02);
 \draw[edge] (a)--(r20);
 \draw[edge] (a)--(r21);
 \draw[edge] (a)--(r22);
 \draw[edge] (b)--(r00);
 \draw[edge] (b)--(r01);
 \draw[edge] (b)--(r02);
 \draw[edge] (b)--(r10);
 \draw[edge] (b)--(r11);
 \draw[edge] (b)--(r12);
 \draw[edge] (c)--(r10);
 \draw[edge] (c)--(r11);
 \draw[edge] (c)--(r12);
 \draw[edge] (c)--(r20);
 \draw[edge] (c)--(r21);
 \draw[edge] (c)--(r22);
 \node[shared] at (a) {};
 \node[shared] at (b) {};
 \node[shared] at (c) {};
 \node[vertex] at (r00) {};
 \node[vertex] at (r01) {};
 \node[vertex] at (r02) {};
 \node[vertex] at (r10) {};
 \node[vertex] at (r11) {};
 \node[vertex] at (r12) {};
 \node[vertex] at (r20) {};
 \node[vertex] at (r21) {};
 \node[vertex] at (r22) {};
 \draw[group] (1.95,2.05)--(1.95,0.95) node[midway,right=5pt,font=\footnotesize] {$H_1$};
 \draw[group] (1.95,0.55)--(1.95,-0.55) node[midway,right=5pt,font=\footnotesize] {$H_2$};
 \draw[group] (1.95,-0.95)--(1.95,-2.05) node[midway,right=5pt,font=\footnotesize] {$H_3$};
 \node[align=center,text width=5.8cm,anchor=north] at (0,-2.4) {(b) Three copies of $K_{2,3}$\\pairwise one-vertex overlaps};
\end{scope}
\end{tikzpicture}
\caption{Biclique joins for $r=2$, with shared vertices
in red. Panel (a) shows a one-vertex join, and (b) shows three copies of $K_{2,3}$
with distinct pairwise shared vertices.}
\label{fig:gluing}
\end{figure}

For a finite simple graph $H$, let $(U_a)_{a\in V(H)}$ be independent random
variables uniformly distributed on $[0,1]$. Define the \emph{homomorphism density} of $H$ by
\begin{align*}
    t(H,W):= {\mathbb E}\left[
        \prod_{\{a,b\}\in E(H)}W(U_a,U_b)
    \right]
    =\int_{[0,1]^{V(H)}}
    \prod_{uv\in E(H)} W(x_u,x_v)\,d\mathbf x_{V(H)}.
\end{align*}
Also, for $u\in V(H)$, define its \emph{one-point conditional density function}
\begin{align}
    t_u(H,W)(y)
    &:={\mathbb E}\left[
        \prod_{\{a,b\}\in E(H)}W(U_a,U_b)
        \,\middle|\,U_u=y
    \right]
    \notag\\
    &=
    \int_{[0,1]^{V(H)\setminus\{u\}}}
    \prod_{\{a,b\}\in E(H)}W(x_a,x_b)
    \prod_{a\in V(H)\setminus\{u\}}dx_a,
    \qquad x_u=y.
    \label{eq:conditional-homomorphism-density}
\end{align}
This gives
\begin{align}
    t(H\oplus_{u,v}H',W)
    &=
    \int_0^1 t_u(H,W)(y)t_v(H',W)(y)\,dy.
    \label{eq:vertex-join-density-factorization}
\end{align}

For a finite simple graph $H$ and $s\in[n]^{V(H)}$, set
\begin{align}
    G_n(H,s)
    :=
    \prod_{uv\in E(H)}
    G_n\bigl(s(u),s(v)\bigr),
    \label{eq:homo-weights}
\end{align}
so that (see \eqref{eq:function-notation})
\begin{align}
    \InjHom(H,G_n)
    =
    \sum_{s\in[n]_{V(H)}} G_n(H,s),\quad \InjHom(H,\tld G_n)
    =
    \sum_{s\in[n]_{V(H)}} G_n(H,s)\prod_{a\in V(H)} X_{s(a)}.
    \label{eq:inj-polynomial}
\end{align}
In particular, we have $G_n(K_{r,k},s)=\prod_{a=1}^r \prod_{b=r+1}^{r+k} G_n(s(a), s(b))$ for all $r,k\in \bN$ and $s\in [n]^{r+k}$.

\subsection{Kernel estimates}

We use the following standard comparison between the cut norm and the
associated $L^\infty-L^1$ operator norm; see
\citet[Section~4, in particular (4.2)--(4.4) and
Remarks~4.2--4.3]{Janson2013Graphons}. We include the short proof for
completeness.

\begin{lemma}\label[lemma]{lem:cut-norm-L1-Linfty-ineq}
Let $t\in \bN$ and $F\in L^1([0,1]^t\times[0,1])$. For
$\phi\in L^\infty([0,1]^t)$ and $\psi\in L^\infty([0,1])$,
\begin{align*}
    \left|
    \int_{[0,1]^t\times[0,1]}
    \phi(\mathbf x)\psi(y)F(\mathbf x;y)\,d\mathbf x\,dy
    \right|
    \le 4\|\phi\|_\infty\|\psi\|_\infty\|F\|_{\square,t}.
\end{align*}
\end{lemma}

\begin{proof}
First suppose $0\le\phi,\psi\le1$. Then
\begin{align*}
    \int_{[0,1]^t\times[0,1]} \phi(\mathbf x)\psi(y)F(\mathbf x;y)\,d\mathbf x\,dy
    =
    \int_0^1\int_0^1
    \int_{\{\phi>u\}\times\{\psi>v\}}
    F(\mathbf x;y)\,d\mathbf x\,dy\,du\,dv .
\end{align*}
Taking absolute values and applying the definition of $\|F\|_{\square,t}$ (see \eqref{eq:cut-norm-definition}) gives the desired bound with constant $1$ in this nonnegative case. For general bounded $\phi$ and $\psi$, rescale by their suprema and apply the preceding estimate to the four pairs of positive and negative parts.
\end{proof}

\begin{remark}
Suppose that \cref{ass:weights,ass:star-cut-convergence,ass:bounded-codegrees}
hold, and fix $1\le t\le r$. Apply
\cref{lem:cut-norm-L1-Linfty-ineq} with
$F=K_{W_n}^{(t)}-K_W^{(t)}$,
$\phi(\mathbf x)=\operatorname{sgn}(D_{W_n}^{(t)}(\mathbf x)-D_W^{(t)}(\mathbf x))$,
and $\psi\equiv1$. This gives
\begin{align*}
    \|D_{W_n}^{(t)}-D_W^{(t)}\|_1
    \le4\|K_{W_n}^{(t)}-K_W^{(t)}\|_{\square,t}\overset{n\to\infty}{\longrightarrow} 0.
\end{align*}
As noted after the assumptions, $0\le D_W^{(t)}\le M$ almost everywhere.
For $k\ge1$, use $|x^k-y^k|\le kM^{k-1}|x-y|$ on $[0,M]$,
applying the finite-graph bound only on $\Delta_{n,t}$, to obtain
\begin{align*}
    &\|(D_{W_n}^{(t)})^k\ind_{\Delta_{n,t}}-(D_W^{(t)})^k\|_1\\
    &\qquad\le kM^{k-1}\|D_{W_n}^{(t)}-D_W^{(t)}\|_1
             +M^k|[0,1]^t\setminus\Delta_{n,t}|\\
    &\qquad\lesssim_{M,k,t}
       \|K_{W_n}^{(t)}-K_W^{(t)}\|_{\square,t}+n^{-1}.
\end{align*}
For $k=0$, the difference is just
$\ind_{\Delta_{n,t}}-1$. Therefore, for every integer $k\ge0$,
\begin{align}
    (D_{W_n}^{(t)})^k\ind_{\Delta_{n,t}}
    &\overset{n\to\infty}{\longrightarrow}(D_W^{(t)})^k\quad\text{in }L^1([0,1]^t),
    &\|D_W^{(t)}\|_\infty&\le M,
    \quad1\le t\le r.
    \label{eq:codegree-L1conv-bdd}
\end{align}
\end{remark}

\subsection{Bounds for connected collections of copies}

The following estimates control the connected overlaps that appear in
the covariance calculations and in the Stein bound.

\begin{definition}\label[definition]{def:connected-copies}
Let $m\in \bN$, let $H_1,\ldots,H_m$ be finite simple graphs, and $s_j\in[n]_{V(H_j)}$, $1\le j\le m$. Define
$\mathcal G(s_1,\ldots,s_m)$ on the vertex set $[m]=\{1,\dots,m\}$, with an edge between
$i$ and $j$ whenever $\bar s_i\cap \bar s_j\neq\emptyset$. We call $(s_1,\ldots,s_m)$ connected if
$\mathcal G(s_1,\ldots,s_m)$ is connected, and denote the set of all
such tuples by $\mathcal K_{n,m}(H_1,\ldots,H_m)$. Also, for $h\in\bN\cup\{0\}$, let $\mathcal K_{n,m}^{(h)}(H_1,\ldots,H_m)$ be the set of the tuples $(s_1,\ldots,s_m)\in\mathcal K_{n,m}(H_1,\ldots,H_m)$ such that
\begin{align}\label{eq:v}
    \left|\bigcup_{j=1}^m\bar s_j\right|
    =
    \sum_{j=1}^m|V(H_j)|-(m-1)-h.
\end{align}
\end{definition}

Note that, to ensure connectedness of $\mathcal G(s_1,\ldots,s_m)$, necessarily we must have $h\ge 0$, as
\begin{align*}
    \left|\bigcup_{j=1}^m\bar s_j\right|\le \sum_{j=1}^m |V(H_j)| - (m-1),\quad (s_1,\ldots,s_m)\in \mathcal K_{n,m}.
\end{align*}

\begin{lemma}\label[lemma]{lem:multijoin_bound}
Suppose that $(G_n)$ satisfies
\cref{ass:weights,ass:bounded-codegrees}, let $m\in \bN$
and $H_1,\ldots,H_m\in \{K_{r,k}:k\ge1\}\cup\{K_1,K_{1,1}\}$ be fixed, and set $v:=\sum_{j=1}^m |V(H_j)|,\ e:=\sum_{j=1}^m |E(H_j)|$. Then
\begin{align*}
    \hspace{-0.4in}
    \sum_{(s_1,\ldots,s_m)\in
        \mathcal K_{n,m}^{(h)}(H_1,\ldots,H_m)}
    \prod_{j=1}^mG_n(H_j,s_j)\le
    2^{\binom{v}{2}}(vM)^v
    n^{v-(m-1)}q_n^e\left(nq_n^r\right)^{-h}.
\end{align*}
\end{lemma}

The proof of \cref{lem:multijoin_bound} is given in Appendix \ref{sec:join_bounds}.

\subsection{Biclique and one-vertex-join counts}

The covariance structure of the empirical codegree process involves bicliques and one-vertex joins of two bicliques.

\begin{lemma}
\label[lemma]{lem:counting-lemma}
Suppose that the graph sequence $(G_n)$ satisfies
\cref{ass:weights,ass:star-cut-convergence,ass:bounded-codegrees}, and that
$nq_n^r\to\infty$.
Then the following hold.

\emph{(i)} For every fixed $k\ge1$,
\begin{align*}
    \left|
        \int x^k\,d\nu_n(x)
        -t_{\mathrm{inj}}(K_{r,k},G_n)
    \right|
    &\lesssim_{r,k,M}(nq_n^r)^{-1}.
\end{align*}

\emph{(ii)} For every fixed $k\ge1$,
\begin{align*}
    t_{\mathrm{inj}}(K_{r,k},G_n)
    &\longrightarrow t(K_{r,k},W),
\end{align*}
\emph{(iii)} For every fixed
$H,H'\in\{K_{r,k}:k\ge1\}\cup\{K_1,K_{1,1}\}$,
$u\in V(H)$, and $v\in V(H')$,
\begin{align*}
    t_{\mathrm{inj}}(H\oplus_{u,v}H',G_n)
    &\longrightarrow
    t(H\oplus_{u,v}H',W).
\end{align*}
\end{lemma}

\begin{proof}[Proof of \cref{lem:counting-lemma}]
By \eqref{eq:inj-density-definition} and \eqref{eq:homo-weights}, we have
\begin{align*}
    t_{\mathrm{inj}}(K_{r,k},G_n)=\frac{1}{(n)_{r+k}q_n^{rk}}\sum_{s\in [n]_{V(K_{r,k})}} G_n(K_{r,k},s).
\end{align*}
Also, by \eqref{eq:codegree-distribution-definition}
\begin{align}\label{eq:moment}
   \int x^k\ d\nu_n(x)
    &=\frac{1}{\binom{n}{r}}\sum_{\substack{S\subseteq [n]\\|S|=r}} \left(\frac{d_S}{nq_n^r}\right)^k
    &=\frac{1}{(n)_r n^kq_n^{rk}}\sum_{\substack{s\in[n]^{V(K_{r,k})}\\
        s|_{L(K_{r,k})}\ \mathrm{injective}}}
        G_n(K_{r,k},s),
\end{align}
where $V(K_{r,k})$ and $L(K_{r,k})$ are as defined in \cref{subsec:distributional-results}.

Moreover, since $G_n(i,i)=0$, every nonzero term $G_n(K_{r,k},s)$
corresponding to a noninjective map $s$ must identify two vertices of
$R(K_{r,k})$. Hence
\begin{align*}
    \left|
        \int x^k\,d\nu_n(x)
        -
        t_{\mathrm{inj}}(K_{r,k},G_n)
    \right|
    &\le
    \left|
        1-\frac{(n)_r n^k}{(n)_{r+k}}
    \right|
    \int x^k\,d\nu_n(x)\\
    &+
    \frac{1}{(n)_{r+k}q_n^{rk}}\sum_{1\le a<b\le k}
    \sum_{\substack{s\in[n]^{V(K_{r,k})}\\
        s|_{L(K_{r,k})}\ \mathrm{injective}\\
        s(r+a)=s(r+b)}}
        G_n(K_{r,k},s).
\end{align*}
By \cref{ass:bounded-codegrees},
\begin{align*}
    &\left|
        1-\frac{(n)_r n^k}{(n)_{r+k}}
    \right|
    \int x^k\,d\nu_n(x)
    \le
    M^k
    \left|
        1-\frac{(n)_r n^k}{(n)_{r+k}}
    \right|
    \lesssim_{r,k}\frac{M^k}{n}
    \le \frac{M^k}{nq_n^r}.
\end{align*}

If $k=1$, every noninjective $s\in[n]^{V(K_{r,1})}$
with $s|_{L(K_{r,1})}$ injective identifies $r+1$ with some vertex in $[r]$ and therefore has zero weight, as $G_n(i,i)=0$ for all $i\in [n]$. Hence the second term above vanishes.

It remains to consider $k\ge2$. For such $k$,
\begin{align}
    &\frac{1}{(n)_{r+k}q_n^{rk}}
    \sum_{1\le a<b\le k}
    \sum_{\substack{s\in[n]^{V(K_{r,k})}\\
        s|_{L(K_{r,k})}\ \mathrm{injective}\\
        s(r+a)=s(r+b)}}
    G_n(K_{r,k},s)\notag\\
    &=
    \frac{1}{(n)_{r+k}q_n^{rk}}
    \sum_{1\le a<b\le k}
    \sum_{\substack{s\in[n]^{V(K_{r,k})}\\
        s|_{L(K_{r,k})}\ \mathrm{injective}\\
        s(r+a)=s(r+b)}}
    \prod_{u=1}^r\prod_{v=r+1}^{r+k}
    G_n(s(u),s(v))\notag\\
    &\overset{\ref{ass:weights}}{\le}
    \binom{k}{2}
    \frac{1}{(n)_{r+k}q_n^{rk}}
    \sum_{\substack{s'\in[n]^{V(K_{r,k-1})}\\
        s'|_{L(K_{r,k-1})}\ \mathrm{injective}}}
    \prod_{u=1}^r\prod_{v=r+1}^{r+k-1}
    G_n(s'(u),s'(v))\notag\\
    &=
    \binom{k}{2}
    \frac{(n)_r n^{k-1}q_n^{r(k-1)}}
         {(n)_{r+k}q_n^{rk}}
    \int x^{k-1}\,d\nu_n(x)
    \lesssim_{r,k}
    M^{k-1}(nq_n^r)^{-1}.
    \label{eq:biclique-repeated-right-bound}
\end{align}
where the last equality uses arguments similar to \eqref{eq:moment}.
This proves \emph{(i)}.

In order to prove \emph{(ii)}, note that since each cell in $\Delta_{n,r}$ has Lebesgue measure $n^{-r}$, then
\begin{align*}
    \int x^k d\nu_n(x) = \frac{n^r}{(n)_r} \int_{\Delta_{n,r}}
    \left(D_{W_n}^{(r)}(\mathbf x)\right)^k\ d\mathbf x
\end{align*}
Moreover, Assumption~\ref{ass:bounded-codegrees}, and \eqref{eq:codegree-L1conv-bdd} imply
\begin{align*}
    & \int_{\Delta_{n,r}}
        \bigl(D_{W_n}^{(r)}(\mathbf x)\bigr)^k\,d\mathbf x
        \overset{n\to \infty}{\longrightarrow}
        \int_{[0,1]^r}
        \bigl(D_W^{(r)}(\mathbf x)\bigr)^k\,d\mathbf x
    =t(K_{r,k},W)
\end{align*}
for every $k\in \bN$. This and \emph{(i)} imply \emph{(ii)}.

Now, we prove \emph{(iii)}. By
\eqref{eq:conditional-homomorphism-density}, for every $k\in \bN$ and $u\in V(K_{r,k})$
\begin{align}
    t_u(K_{r,k},W)(y)
    &=
    \begin{cases}
    \displaystyle
    \int_{[0,1]^{r-1}}
    \bigl(D_W^{(r)}(y,x_2,\ldots,x_r)\bigr)^k
    \,dx_2\cdots dx_r,
    &u\in L(K_{r,k}),
    \\[0.15in]
    \displaystyle
    \int_{[0,1]^r}
    K_W^{(r)}(\mathbf x;y)
    \bigl(D_W^{(r)}(\mathbf x)\bigr)^{k-1}
    \,d\mathbf x,
    &u\in R(K_{r,k}),
    \end{cases}
    \label{eq:cond-biclique-density-identity}
\end{align}
where $K_W^{(r)}$ and $D_W^{(r)}$ are as in \eqref{eq:codegree-kernel-definition}. For a left vertex $u_L\in L(K_{r,k})$ and a right vertex
$u_R\in R(K_{r,k})$, define the semi-injective conditional densities by
\begin{align*}
    t_{u_L}^{\mathrm{si}}(K_{r,k},G_n)(y)
    &:=\int_{[0,1]^{r-1}}\ind_{\Delta_{n,r}}(y,\mathbf x)
       \bigl(D_{W_n}^{(r)}(y,\mathbf x)\bigr)^k\,d\mathbf x,\\
    t_{u_R}^{\mathrm{si}}(K_{r,k},G_n)(y)
    &:=\int_{\Delta_{n,r}}K_{W_n}^{(r)}(\mathbf x;y)
       \bigl(D_{W_n}^{(r)}(\mathbf x)\bigr)^{k-1}\,d\mathbf x.
\end{align*}
Now fix $k,\ell\in \bN$ and $u\in V(K_{r,k}),\ v\in V(K_{r,\ell})$, and set $J:=K_{r,k}\oplus_{u,v}K_{r,\ell}$, and note that $|V(J)|=2r+k+\ell-1,\ |E(J)|=r(k+\ell)$. We establish the following facts, from which \emph{(iii)} follows.

\textbf{Claim 1.} We have
\begin{align*}
    &\int_0^1t_u^{\mathrm{si}}(K_{r,k},G_n)(y)
              t_v^{\mathrm{si}}(K_{r,\ell},G_n)(y)\,dy\\
    &\qquad\longrightarrow
      \int_0^1t_u(K_{r,k},W)(y)t_v(K_{r,\ell},W)(y)\,dy
      =t(J,W),
\end{align*}
where the equality is \eqref{eq:vertex-join-density-factorization}.

\textbf{Claim 2.} We have
\begin{align*}
    &\left|t_{\mathrm{inj}}(J,G_n)
       -\int_0^1t_u^{\mathrm{si}}(K_{r,k},G_n)(y)
                t_v^{\mathrm{si}}(K_{r,\ell},G_n)(y)\,dy\right|\\
    &\qquad\lesssim_{r,k,\ell,M}(nq_n^r)^{-1}.
\end{align*}

We first prove Claim 1. It suffices to show that, for every
$k\in\bN$ and $u\in V(K_{r,k})$,
\begin{align}
    t_u^{\mathrm{si}}(K_{r,k},G_n)\overset{n\to \infty}{\longrightarrow}
    t_u(K_{r,k},W)\quad \text{in }L^1([0,1]),
    \notag\\
    \sup_n\|t_u^{\mathrm{si}}(K_{r,k},G_n)\|_\infty
    \vee \|t_u(K_{r,k},W)\|_\infty
    &\le M^{r+k-1}.
    \label{eq:conditional-si-convergence}
\end{align}

Suppose first that $u\in L(K_{r,k})$. By Fubini's theorem,
$|x^k-y^k|\le kM^{k-1}|x-y|$ for $x,y\in[0,M]$, and
\eqref{eq:codegree-L1conv-bdd},
\begin{align*}
    \|t_u^{\mathrm{si}}(K_{r,k},G_n)
        -t_u(K_{r,k},W)\|_1
    &\le
    kM^{k-1}
    \|D_{W_n}^{(r)}\ind_{\Delta_{n,r}}
        -D_W^{(r)}\|_1
    \longrightarrow0.
\end{align*}
The defining formulas and
\eqref{eq:codegree-L1conv-bdd} also give
\begin{align*}
    \|t_u^{\mathrm{si}}(K_{r,k},G_n)\|_\infty
    \vee
    \|t_u(K_{r,k},W)\|_\infty
    \le M^k
    \le M^{r+k-1}.
\end{align*}
Now suppose that $u\in R(K_{r,k})$. Write
\begin{align*}
    F_n&:=K_{W_n}^{(r)}-K_W^{(r)},\\
    \phi_n(\mathbf x)&:=\bigl(D_{W_n}^{(r)}(\mathbf x)\bigr)^{k-1}
                         \ind_{\Delta_{n,r}}(\mathbf x),
    &\phi(\mathbf x)&:=\bigl(D_W^{(r)}(\mathbf x)\bigr)^{k-1}.
\end{align*}
Then the difference of the conditional densities equals
\begin{align*}
    &t_u^{\mathrm{si}}(K_{r,k},G_n)(y)-t_u(K_{r,k},W)(y)\\
    &\qquad=\int_{[0,1]^r}F_n(\mathbf x;y)\phi_n(\mathbf x)\,d\mathbf x
       +\int_{[0,1]^r}K_W^{(r)}(\mathbf x;y)
                       (\phi_n-\phi)(\mathbf x)\,d\mathbf x.
\end{align*}
Since $\|\phi_n\|_\infty\le M^{k-1}$, apply
\cref{lem:cut-norm-L1-Linfty-ineq} with
$\psi(y)=\operatorname{sgn}(\int F_n(\mathbf x;y)\phi_n(\mathbf x)\,d\mathbf x)$
to obtain
\begin{align*}
    \left\|\int_{[0,1]^r}F_n(\mathbf x;\cdot)\phi_n(\mathbf x)\,d\mathbf x\right\|_1
    &\le4M^{k-1}\|F_n\|_{\square,r}\longrightarrow0.
\end{align*}
For the second term, Tonelli's theorem gives
\begin{align*}
    &\left\|\int_{[0,1]^r}K_W^{(r)}(\mathbf x;\cdot)
                         (\phi_n-\phi)(\mathbf x)\,d\mathbf x\right\|_1\\
    &\qquad\le\int_{[0,1]^r}D_W^{(r)}(\mathbf x)
                         |\phi_n-\phi|(\mathbf x)\,d\mathbf x
       \le M\|\phi_n-\phi\|_1\longrightarrow0,
\end{align*}
by \eqref{eq:codegree-L1conv-bdd}.
Finally, for every $y\in[0,1]$,
\begin{align*}
    t_u^{\mathrm{si}}(K_{r,k},G_n)(y)
    &\le
    M^{k-1}
    \int_{[0,1]^r}K_{W_n}^{(r)}(\mathbf x;y)\,d\mathbf x\\
    &=
    M^{k-1}\bigl(D_{W_n}(y)\bigr)^r
    \le M^{r+k-1}.
\end{align*}
The same calculation with $W$ in place of $W_n$ gives $\|t_u(K_{r,k},W)\|_\infty\le M^{r+k-1}$. This proves \eqref{eq:conditional-si-convergence}, and Claim 1 follows.

Now, we prove Claim 2. Let $\mathcal A_n$ be the set of pairs
$(s_1,s_2)$ with $s_1\in[n]^{V(K_{r,k})}$ and
$s_2\in[n]^{V(K_{r,\ell})}$
such that $s_1|_{L(K_{r,k})}$ and $s_2|_{L(K_{r,\ell})}$ are
injective and $s_1(u)=s_2(v)$. For $\mathcal C\subseteq\mathcal A_n$, set
\begin{align*}
    w_n(\mathcal C):=\sum_{(s_1,s_2)\in \mathcal C} G_n(K_{r,k},s_1)G_n(K_{r,\ell},s_2)
\end{align*}
and note that for $i\in [n]$ and $y\in I_i=\left(\frac{i-1}{n},\frac{i}{n}\right]$, we have the identity (see \eqref{eq:cond-biclique-density-identity})
\begin{align}
    t_u^{\mathrm{si}}(K_{r,k},G_n)(y)
    &=
    \sum_{\substack{s\in[n]^{V(K_{r,k})}\\
        s(u)=i,\;s|_{L(K_{r,k})}\ \mathrm{injective}}}
    \frac{G_n(K_{r,k},s)}{n^{r+k-1}q_n^{rk}}.
    \label{eq:cond-si-discrete-representation}
\end{align}
Expanding both conditional densities and grouping by the common image
of the identified vertex gives
\begin{align*}
    \int_0^1t_u^{\mathrm{si}}(K_{r,k},G_n)(y)
             t_v^{\mathrm{si}}(K_{r,\ell},G_n)(y)\,dy
    &=\frac{w_n(\mathcal A_n)}{n^{2r+k+\ell-1}q_n^{r(k+\ell)}}.
\end{align*}
Recall that for $s\in [n]^A$, we write $\bar s:=s(A)=\{s(a):a\in A\}$, and partition $\mathcal A_n$ into the following subsets:
\begin{align*}
    \mathcal A_n^{(0)}
    &:=
    \{(s_1,s_2)\in\mathcal A_n:
      s_1,s_2\text{ are injective and }|\bar s_1\cap\bar s_2|=1\},\\
    \mathcal A_n^{(1)}
    &:=
    \{(s_1,s_2)\in\mathcal A_n:
      s_1,s_2\text{ are injective and }
      |\bar s_1\cap\bar s_2|\ge2\},\\
    \mathcal A_n^{(2)}
    &:=
    \{(s_1,s_2)\in\mathcal A_n:
      s_1\text{ is not injective}\},\\
    \mathcal A_n^{(3)}
    &:=
    \{(s_1,s_2)\in\mathcal A_n:
      s_1\text{ is injective and }s_2\text{ is not injective}\}.
\end{align*}

The class $\mathcal A_n^{(0)}$ is in bijection with the injective
maps from $J$ to $[n]$. Moreover,
\begin{align*}
    w_n(\mathcal A_n^{(0)})
    &\le w_n(\mathcal A_n)
    = n^{2r+k+\ell-1}q_n^{r(k+\ell)}
    \int_0^1
    t_u^{\mathrm{si}}(K_{r,k},G_n)(y)
    t_v^{\mathrm{si}}(K_{r,\ell},G_n)(y)\,dy\\
    & \overset{\eqref{eq:conditional-si-convergence}}{\le} M^{2r+k+\ell-2}
    n^{2r+k+\ell-1}q_n^{r(k+\ell)}.
\end{align*}

For $(s_1,s_2)\in\mathcal A_n^{(1)}$, write $|\bar s_1\cap\bar s_2|=h+1$. Then $h\ge1$, and the union of the two copies has $2r+k+\ell-1-h$
vertices. By \cref{lem:multijoin_bound} with $m=2,\ H_1=K_{r,k}$, and $H_2=K_{r,\ell}$,
\begin{align*}
    w_n(\mathcal A_n^{(1)})
    &\lesssim_{r,k,\ell,M}
      \sum_{h=1}^{r+\min\{k,\ell\}-1}
      n^{2r+k+\ell-1-h}q_n^{r(k+\ell)-rh}\\
    &\lesssim_{r,k,\ell,M}
      n^{2r+k+\ell-2}q_n^{r(k+\ell-1)},
\end{align*}
where the last inequality uses $nq_n^r\ge 1$.

For each $i\in[n]$, fix $y_i\in I_i=\left(\frac{i-1}{n},\frac{i}{n}\right]$. Also, if $(s_1,s_2)\in\mathcal A_n^{(2)}$, every nonzero term $G_n(K_{r,k},s_1)G_n(K_{r,\ell},s_2)$ satisfies $s_1(r+a)=s_1(r+b)$ for some $1\le a<b\le k$. Therefore,
\begin{align*}
    w_n(\mathcal A_n^{(2)})
    &\le
    \sum_{1\le a<b\le k}
    \sum_{\substack{
        s_1\in[n]^{V(K_{r,k})}\\
        s_1|_{L(K_{r,k})}\ \mathrm{injective}\\
        s_1(r+a)=s_1(r+b)}}
    G_n(K_{r,k},s_1)
    \sum_{\substack{
        s_2\in[n]^{V(K_{r,\ell})}\\
        s_2(v)=s_1(u)\\
        s_2|_{L(K_{r,\ell})}\ \mathrm{injective}}}
    G_n(K_{r,\ell},s_2)\\
    &=
    n^{r+\ell-1}q_n^{r\ell}
    \sum_{1\le a<b\le k}
    \sum_{\substack{
        s_1\in[n]^{V(K_{r,k})}\\
        s_1|_{L(K_{r,k})}\ \mathrm{injective}\\
        s_1(r+a)=s_1(r+b)}}
    G_n(K_{r,k},s_1)
    t_v^{\mathrm{si}}(K_{r,\ell},G_n)
        (y_{s_1(u)})\\
        &\le
    n^{r+\ell-1}q_n^{r\ell}
    \bigl\|
        t_v^{\mathrm{si}}(K_{r,\ell},G_n)
    \bigr\|_\infty
    \sum_{1\le a<b\le k}
    \sum_{\substack{
        s_1\in[n]^{V(K_{r,k})}\\
        s_1|_{L(K_{r,k})}\ \mathrm{injective}\\
        s_1(r+a)=s_1(r+b)}}
    G_n(K_{r,k},s_1)\\
    &\lesssim_{r,k,\ell,M}
    n^{r+\ell-1}q_n^{r\ell}\cdot
    n^{r+k} q_n^{rk}\cdot
    (nq_n^r)^{-1}
    =
    n^{2r+k+\ell-2}q_n^{r(k+\ell-1)},
\end{align*}
where the last inequality above follows from the same bounds in \eqref{eq:biclique-repeated-right-bound}. Similarly, we have $w_n(\mathcal A_n^{(3)})\lesssim_{r,k,\ell,M} n^{2r+k+\ell-2}q_n^{r(k+\ell-1)}$, then
\begin{align}
    &\left|
        t_{\mathrm{inj}}(J,G_n)
        -
        \int_0^1
        t_u^{\mathrm{si}}(K_{r,k},G_n)(y)
        t_v^{\mathrm{si}}(K_{r,\ell},G_n)(y)\,dy
    \right|\label{eq:tinj_tsi_comparison}\\
    &\quad\le
    \left|
        \frac{n^{2r+k+\ell-1}}{(n)_{2r+k+\ell-1}}-1
    \right|
    \frac{w_n(\mathcal A_n^{(0)})}
         {n^{2r+k+\ell-1}q_n^{r(k+\ell)}}
    +
    \frac{
        \sum_{a=1}^3 w_n(\mathcal A_n^{(a)})
    }{n^{2r+k+\ell-1}q_n^{r(k+\ell)}}\notag\\
    &\quad\lesssim_{r,k,\ell,M}
    n^{-1}
    +
    \frac{n^{2r+k+\ell-2}q_n^{r(k+\ell-1)}}
         {n^{2r+k+\ell-1}q_n^{r(k+\ell)}}
    \lesssim_{r,k,\ell,M}
    (nq_n^r)^{-1}.\notag
\end{align}
The preceding argument proves \emph{(iii)} for $H=K_{r,k}$ and
$H'=K_{r,\ell}$. We next treat explicitly the joins involving
$K_{1,1}$ and $K_1$. Fix $k\ge1$,
$u\in V(K_{r,k})$, and $v\in V(K_{1,1})$. The argument used to prove
Claim~2, applied with $K_{1,1}$ in place of $K_{r,\ell}$ and noting
that its one-point conditional density is $D_{W_n}$, gives
\begin{align*}
    t_{\mathrm{inj}}
    (K_{r,k}\oplus_{u,v}K_{1,1},G_n)
    &=
    \int_0^1
    t_u^{\mathrm{si}}(K_{r,k},G_n)(y)D_{W_n}(y)\,dy + O\left((nq_n^r)^{-1}\right)\\
    &\overset{n\to \infty}{\longrightarrow}
    \int_0^1t_u(K_{r,k},W)(y)D_W(y)\,dy = t(K_{r,k}\oplus_{u,v}K_{1,1},W),
\end{align*}
by \eqref{eq:conditional-si-convergence},
\eqref{eq:codegree-L1conv-bdd}, and
\eqref{eq:vertex-join-density-factorization}. Similarly, for $u,v\in V(K_{1,1})$,
\begin{align*}
    t_{\mathrm{inj}}
    (K_{1,1}\oplus_{u,v}K_{1,1},G_n)\overset{n\to\infty}{\longrightarrow}
    \int_0^1D_W(y)^2\,dy =t(K_{1,1}\oplus_{u,v}K_{1,1},W),
\end{align*}
by \eqref{eq:codegree-L1conv-bdd} and
\eqref{eq:vertex-join-density-factorization}. If at least one of $H,H'$, say $H'$, is $K_1$, then $H\oplus_{u,v}H'$ is isomorphic to $H$. Thus the claim follows from \emph{(ii)}, while the case $K_1$ and $K_{1,1}$ is immediate. This completes the proof of \emph{(iii)}.
\end{proof}

\begin{remark}\label[remark]{rem:deterministic-counting-comparisons}
The estimate in \cref{lem:counting-lemma}\emph{(i)}
is deterministic and uses only
\cref{ass:weights,ass:bounded-codegrees}.
Consequently, it also applies to $\tld G_n$
on the event $\mathcal E_n$ defined in
\cref{lem:normalizer_codegree_estimates},
with $(G_n,n,q_n)$ replaced by
$(\tld G_n,N_n,\tld q_n)$,
provided $N_n\tld q_n^r\ge1$.
\end{remark}

\subsection{Normal approximation for sampled counts}
\label{subsec:sampled-count-normal-approximation}

We now use the distance $d_3$, centered counts $Y_n(H)$, ordered family
$\mathcal H$, and covariance matrix $\Sigma_n$ introduced in
\cref{subsec:distributional-results}.

To obtain the multivariate normal approximation for the sampled biclique
counts, we use the following form of
\cite[Theorem~2]{Raic2004Decomposable} for a general positive
semidefinite covariance matrix. For completeness, a proof using Stein's
method is given in \cref{app:stein-bounds}.

\begin{theorem}\label[theorem]{thm:raic-decomposable-bound}
Let $\mathcal H$ be a finite indexed family of graphs and let
$X_1,\ldots,X_n$ be independent random variables. Let $\left\{\xi_n(H,s):H\in\mathcal H,\ s\in[n]_{V(H)}\right\}$ be a centered, square-integrable $\mathbb R^m$-valued stochastic process such that $\xi_n(H,s)$ is measurable with respect to $\sigma(X_i:i\in\bar s)$. Set
\begin{align*}
    S&:=\sum_{H\in\mathcal H}
    \sum_{s\in[n]_{V(H)}}\xi_n(H,s),
    \qquad
    \Sigma:=\Var(S),
\end{align*}
and let $\mathcal K_{n,3}(H_1,H_2,H_3)$ be as in
Definition~\ref{def:connected-copies}. Then
\begin{align*}
    &d_3(S,\mathcal N(0,\Sigma))\\
    &\quad\le\frac13\sum_{H_1,H_2,H_3\in\mathcal H}
       \sum_{(s_1,s_2,s_3)\in\mathcal K_{n,3}(H_1,H_2,H_3)}
       \Biggl\{\E\prod_{a=1}^3\|\xi_n(H_a,s_a)\|\\
    &\hspace{4cm}
       +\E\left[\prod_{a=1}^2\|\xi_n(H_a,s_a)\|\right]
          \E\|\xi_n(H_3,s_3)\|\Biggr\}.
\end{align*}
\end{theorem}

For $H\in\mathcal H$ and $s\in[n]_{V(H)}$, set
\begin{align*}
    Z_{H,s}:=\prod_{u\in V(H)}X_{s(u)}-p_n^{|V(H)|},
\end{align*}
so that
\begin{align}
    Y_n(H)
    =
    \frac{1}{(np_n)^{|V(H)|-1/2}q_n^{|E(H)|}}
    \sum_{s\in[n]_{V(H)}}G_n(H,s)Z_{H,s}.
    \label{eq:copy-indexed-expansion}
\end{align}

Define the matrix of one-vertex-join contributions by
\begin{align}
\Sigma_n^\oplus(H,H')
:=(1-p_n)\sum_{u\in V(H)}\sum_{v\in V(H')}
t_{\mathrm{inj}}(H\oplus_{u,v}H',G_n),
\qquad H,H'\in\mathcal H.
\label{eq:join-covariance-definition}
\end{align}
For matrices indexed by $\mathcal H\times\mathcal H$, write
$\|A\|_{\mathrm F}:=
(\sum_{H,H'\in\mathcal H}|A(H,H')|^2)^{1/2}$
for the Frobenius norm.

\begin{proof}[Proof of \cref{thm:stein-error-bound}]
For $H\in\mathcal H$ and $s\in[n]_{V(H)}$, let $\mathbf e_H$ be the
corresponding coordinate vector in $\mathbb R^{|\mathcal H|}$ and set
\begin{align*}
    \xi_n(H,s)
    :=
    \frac{G_n(H,s)Z_{H,s}}
    {(np_n)^{|V(H)|-1/2}q_n^{|E(H)|}}\,\mathbf e_H.
\end{align*}
Then $\mathbf Y_n=\sum_{H\in\mathcal H}\sum_s\xi_n(H,s)$ by
\eqref{eq:copy-indexed-expansion}.

Fix $H_1,H_2,H_3\in\mathcal H$ and
$s_i\in[n]_{V(H_i)}$, $1\le i\le3$. Since
$|\bar s_i|=|V(H_i)|$ and
$|Z_{H_i,s_i}|=\bigg|\prod_{u\in V(H_i)} X_{s_i(u)}-p_n^{|V(H_i)|}\bigg|
\le \prod_{u\in V(H_i)} X_{s_i(u)}+p_n^{|V(H_i)|}$,
expanding the product and using independence gives
\begin{align*}
    \E\left|
        Z_{H_1,s_1}Z_{H_2,s_2}Z_{H_3,s_3}
    \right|
    &\le
    \sum_{B\subseteq[3]}
    p_n^{
        \sum_{i\notin B}|V(H_i)|
        +
        \left|\bigcup_{i\in B}\bar s_i\right|
    }\le
    8p_n^{|\bar s_1\cup\bar s_2\cup\bar s_3|}.
\end{align*}
Indeed, for every $B\subseteq[3]$, $|\bar s_1\cup\bar s_2\cup\bar s_3|\le\left|\bigcup_{i\in B}\bar s_i\right|+\sum_{i\notin B}|\bar s_i|$,
and $0<p_n\le1$. The same argument gives
\begin{align*}
    \E\left|Z_{H_1,s_1}Z_{H_2,s_2}\right|
    \le
    4p_n^{|\bar s_1\cup\bar s_2|},\qquad
    \E\left|Z_{H_3,s_3}\right|
    \le
    2p_n^{|\bar s_3|}.
\end{align*}
Since
$|\bar s_1\cup\bar s_2\cup\bar s_3|
\le|\bar s_1\cup\bar s_2|+|\bar s_3|$, it follows that
\begin{align*}
    &\E\left|
        Z_{H_1,s_1}Z_{H_2,s_2}Z_{H_3,s_3}
    \right|
    +
    \E\left|Z_{H_1,s_1}Z_{H_2,s_2}\right|
    \E\left|Z_{H_3,s_3}\right|
    \le 16p_n^{|\bar s_1\cup\bar s_2\cup\bar s_3|}.
\end{align*}

The triples in $\mathcal K_{n,3}(H_1,H_2,H_3)$ are connected and hence satisfy $|\bar s_1\cup\bar s_2\cup\bar s_3|=\sum_{i=1}^3|V(H_i)|-2-h$ for some $h\ge0$. Therefore, \cref{thm:raic-decomposable-bound} yields
\begin{align*}
    d_3(\mathbf Y_n,\mathcal N(0,\Sigma_n))
    &\lesssim_{r,m}
    \sum_{H_1,H_2,H_3\in\mathcal H}
    (np_n)^{3/2-\sum_{i=1}^3|V(H_i)|}
    q_n^{-\sum_{i=1}^3|E(H_i)|}\\
    &\qquad\qquad\times
    \sum_{h\ge0}
    p_n^{\sum_{i=1}^3|V(H_i)|-2-h}
    \sum_{\substack{\mathbf s\in\mathcal K_{n,3}^{(h)}(H_1,H_2,H_3)\\ \mathbf s=(s_1,s_2,s_3)}}
    \prod_{i=1}^3G_n(H_i,s_i)\\
    &\overset{\cref{lem:multijoin_bound}}{\lesssim_{r,m,M}}
    (np_n)^{-1/2}
    \sum_{H_1,H_2,H_3\in\mathcal H}
    \sum_{h\ge0}
    (np_nq_n^r)^{-h}\lesssim_{r,m,M,c_0}(np_n)^{-1/2},
\end{align*}
Here only finitely many values of $h$, bounded in terms of $r,m$,
can occur, and $np_nq_n^r\ge c_0$ eventually. This proves
\eqref{eq:stein-error-rate}.
\end{proof}

The next lemma identifies the leading contribution to $\Sigma_n$.

\begin{lemma}\label[lemma]{lem:count-covariance-comparison}
Suppose that the graph sequence $G_n$ satisfies
\cref{ass:weights,ass:bounded-codegrees} and that
$np_nq_n^r\to\infty$. Then
\begin{align*}
\|\Sigma_n^\oplus-\Sigma_n\|_{\mathrm F}
\lesssim_{r,m,M}(np_nq_n^r)^{-1}.
\end{align*}
\end{lemma}

\begin{proof}
Fix $H,H'\in\mathcal H$.
If $s_1\in[n]_{V(H)}$, $s_2\in[n]_{V(H')}$, and
$|\bar s_1\cap\bar s_2|=\ell$, then
\begin{align*}
\Cov(Z_{H,s_1},Z_{H',s_2})
=
p_n^{|V(H)|+|V(H')|-\ell}
-p_n^{|V(H)|+|V(H')|}.
\end{align*}
Thus \eqref{eq:copy-indexed-expansion} gives
\begin{align}
\Sigma_n(H,H')
&=
\sum_{\ell=1}^{\min\{|V(H)|,|V(H')|\}}
\left(
p_n^{|V(H)|+|V(H')|-\ell}
-p_n^{|V(H)|+|V(H')|}
\right)
\notag\\
&\quad\times
\sum_{\substack{
s_1\in[n]_{V(H)},\ s_2\in[n]_{V(H')}\\
|\bar s_1\cap\bar s_2|=\ell}}
\frac{G_n(H,s_1)G_n(H',s_2)}
{(np_n)^{|V(H)|+|V(H')|-1}
q_n^{|E(H)|+|E(H')|}}.
\label{eq:cov-expansion}
\end{align}

Pairs sharing exactly one vertex correspond to injective maps
on a one-vertex join, so
\begin{align*}
\sum_{\substack{
s_1\in[n]_{V(H)},\ s_2\in[n]_{V(H')}\\
|\bar s_1\cap\bar s_2|=1}}
G_n(H,s_1)G_n(H',s_2)
=
\sum_{u\in V(H)}\sum_{v\in V(H')}
\InjHom(H\oplus_{u,v}H',G_n).
\end{align*}
Since each join has $|V(H)|+|V(H')|-1$ vertices and
$|E(H)|+|E(H')|$ edges, the definition of
$t_{\mathrm{inj}}$ shows that the $\ell=1$ term in
\eqref{eq:cov-expansion} is
\begin{align*}
\frac{(n)_{|V(H)|+|V(H')|-1}}
     {n^{|V(H)|+|V(H')|-1}}
\Sigma_n^\oplus(H,H').
\end{align*}

For $\ell\ge2$, pairs sharing $\ell$ vertices belong to
$\mathcal K_{n,2}^{(\ell-1)}(H,H')$.
Applying \cref{lem:multijoin_bound} gives
\begin{align*}
&\sum_{\substack{
s_1\in[n]_{V(H)},\ s_2\in[n]_{V(H')}\\
|\bar s_1\cap\bar s_2|=\ell}}
G_n(H,s_1)G_n(H',s_2)\\
&\qquad\lesssim_{r,m,M}
n^{|V(H)|+|V(H')|-\ell}
q_n^{|E(H)|+|E(H')|-r(\ell-1)}.
\end{align*}
Substituting into \eqref{eq:cov-expansion}, we obtain
\begin{align*}
\left|
\Sigma_n(H,H')
-
\frac{(n)_{|V(H)|+|V(H')|-1}}
     {n^{|V(H)|+|V(H')|-1}}
\Sigma_n^\oplus(H,H')
\right|
&\lesssim_{r,m,M}
\sum_{\ell=2}^{\min\{|V(H)|,|V(H')|\}}
(np_nq_n^r)^{-(\ell-1)}\\
&\lesssim_{r,m,M}(np_nq_n^r)^{-1}.
\end{align*}
The entries of $\Sigma_n^\oplus$ are uniformly bounded by
the same counting bound with $\ell=1$, so replacing
$(n)_{|V(H)|+|V(H')|-1}/n^{|V(H)|+|V(H')|-1}$ by $1$
introduces an error bounded by
$\lesssim_{r,m,M}n^{-1}\le(np_nq_n^r)^{-1}$. Since $m$ is fixed, the resulting entrywise bound completes the proof.
\end{proof}

\begin{proof}[Proof of \cref{cor:biclique-clt}]
By \cref{lem:counting-lemma}\emph{(iii)} and
$p_n\to p_\infty$, the entries of $\Sigma_n^\oplus$ converge.
Together with \cref{lem:count-covariance-comparison}, this gives
\begin{align}
\Sigma_n(H,H')
\longrightarrow
\Sigma_\infty(H,H')
:=(1-p_\infty)
\sum_{u\in V(H)}\sum_{v\in V(H')}
t(H\oplus_{u,v}H',W).
\label{eq:limiting-sampled-count-covariance}
\end{align}
Writing
$\Sigma_\infty:=
(\Sigma_\infty(H,H'))_{H,H'\in\mathcal H}$,
\cref{thm:stein-error-bound} now yields
\begin{align*}
\mathbf Y_n\dto\mathcal N(0,\Sigma_\infty).
\end{align*}
We now apply the delta method. For $1\le k\le m$, let
\begin{align*}
    \tld \theta_{n,k}:=
    \frac{\InjHom(K_{r,k},\tld G_n)}{(np_n)^{r+k}q_n^{rk}},
    \qquad
    \theta_{n,k}
    :=
    \frac{\InjHom(K_{r,k},G_n)}{n^{r+k}q_n^{rk}},
    \qquad
    \mu_k:=t(K_{r,k},W).
\end{align*}
With $N_n=|V(\tld G_n)|$ as in \cref{sec:intro}, using \eqref{eq:centered-counts} and $\InjHom(K_{1,1},\tld G_n)=N_n^2\tld q_n$, we get
\begin{align*}
    Y_n(K_{r,k})=&\sqrt{np_n}\left(\tld \theta_{n,k}-\theta_{n,k}\right),\\
Y_n(K_1)=&\sqrt{np_n}\left(\frac{N_n}{np_n}-1\right),\quad
Y_n(K_{1,1})=\sqrt{np_n}\left(\frac{N_n^2\widetilde q_n}{(np_n)^2q_n}-1\right),
\end{align*}
where $K_1$ and $K_{1,1}$ are as defined in \cref{subsec:distributional-results}. For each $1\le k\le m$, we have $$N_n^{r+k}/(N_n)_{r+k}=1+O_{P}((np_n)^{-1}),\text{ and } n^{r+k}/(n)_{r+k}=1+O(n^{-1}),$$ so these factors are negligible on the $\sqrt{np_n}$ scale. Moreover, using \eqref{eq:inj-density-definition} twice gives
\begin{align*}
    t_{\mathrm{inj}}(K_{r,k},\widetilde G_n)
    =&
    \frac{N_n^{r+k}}{(N_n)_{r+k}}
    \ \tld \theta_{n,k}
    \left(\frac{N_n}{np_n}\right)^{2rk-r-k}
    \left(
        \frac{N_n^2\, \tld q_n}{(np_n)^2 q_n}
    \right)^{-rk}\\
    =&\frac{N_n^{r+k}}{(N_n)_{r+k}}\Phi_k\left((\tld\theta_{n,k})_{k=1}^m,\frac{N_n}{np_n},\frac{N_n^2\, \tld q_n}{(np_n)^2 q_n}\right),\\
    t_{\mathrm{inj}}(K_{r,k},G_n)
    &=\frac{n^{r+k}}{(n)_{r+k}}\theta_{n,k}=\frac{n^{r+k}}{(n)_{r+k}}\Phi_k\left((\theta_{n,k})_{k=1}^m,1,1\right),
\end{align*}
where $\Phi:\R^m\times(0,\infty)^2\longrightarrow\R^m$ is given by
\begin{align}
    \Phi_k(x_1,\ldots,x_m,y,z):=x_k y^{2rk-r-k}z^{-rk}.
    \label{eq:delta_method_Phi}
\end{align}
Also,  \cref{lem:counting-lemma}\emph{(ii)} gives that $\theta_{n,k}\to\mu_k$ for $1\le k\le m$.
Therefore
\begin{align*}
    \sqrt{np_n}
    \left(
        t_{\mathrm{inj}}(K_{r,k},\widetilde G_n)
        -t_{\mathrm{inj}}(K_{r,k},G_n)
    \right)_{k=1}^m
    \dto\mathcal N(0,D_\Phi(\boldsymbol \mu)\Sigma_\infty D_\Phi(\boldsymbol \mu)^T),
\end{align*}
where $\boldsymbol\mu=(\mu_1,\ldots,\mu_m,1,1)$. The Jacobian is
\begin{align*}
    D_\Phi(\boldsymbol\mu)
    =\left(I_m\ \middle|\ ((2rk-r-k)\mu_k)_{k=1}^m
                   \ \middle|\ (-rk\mu_k)_{k=1}^m\right).
\end{align*}
It remains to identify the limiting covariance.
By
\eqref{eq:vertex-join-density-factorization}, and \eqref{eq:limiting-sampled-count-covariance} we have
\begin{align*}
    \Sigma_\infty(H,H')
    &=(1-p_\infty)
    \int_0^1
    \left(\sum_{u\in V(H)}t_u(H,W)(x)\right)
    \left(\sum_{w\in V(H')}t_w(H',W)(x)\right)\,dx.
\end{align*}
Moreover,
\begin{align*}
    \sum_{u\in V(K_1)}t_u(K_1,W)(x)
    =1,\text{ and }
    \sum_{u\in V(K_{1,1})}t_u(K_{1,1},W)(x)
    =2D_W(x),
\end{align*}
while \eqref{eq:cond-biclique-density-identity} gives
\begin{align*}
    \sum_{u\in V(K_{r,k})}t_u(K_{r,k},W)(x)
    &=r\,\E\left[(D_W^{(r)}(U_1,\dots,U_r))^k\mid U_1=x\right]
    \\ &\quad+k\,\E\left[
        \bigl(D_W^{(r)}(U_1,\dots,U_r)\bigr)^{k-1}\prod_{a=1}^rW(U_a,x)
    \right].
\end{align*}
For $1\le k\le m$, set
\begin{align*}
    b_k(z):=\sum_{u\in V(K_{r,k})}t_u(K_{r,k},W)(z)
            +(2rk-r-k)\mu_k-2rk\mu_kD_W(z).
\end{align*}
The preceding identities give, for $1\le k,\ell\le m$,
\begin{align*}
    \bigl(D_\Phi(\boldsymbol\mu)\Sigma_\infty
               D_\Phi(\boldsymbol\mu)^{\mathsf T}\bigr)_{k\ell}
    &=(1-p_\infty)\int_0^1b_k(z)b_\ell(z)\,dz.
\end{align*}
By \eqref{eq:limiting-influence-function},
\begin{align*}
    b_k(z)=\phi_{x^k}(z).
\end{align*}
Therefore
\begin{align}
    \bigl(D_\Phi(\boldsymbol\mu)\Sigma_\infty
               D_\Phi(\boldsymbol\mu)^{\mathsf T}\bigr)_{k\ell}
    =\Gamma(x^k,x^\ell),
    \label{eq:delta-method-process-covariance}
\end{align}
which completes the proof.
\end{proof}

\section{Convergence of the empirical codegree process}
\label{sec:process_convergence}

This section proves \cref{thm:empirical-codegree-process}. We first obtain
uniform variance and bias bounds, then prove finite-dimensional convergence
by reducing polynomial test functions to sampled biclique counts and using
polynomial approximation, and finally prove tightness by
variance-sensitive concentration, chaining, and uniform bias control.
Throughout, use the normalized codegrees in
\eqref{eq:normalized-discrete-codegrees} and the bounded extensions from
\cref{subsec:distributional-results}; thus, when $N_n\ge r$ and
$\tld q_n>0$,
\begin{align*}
\nu_n(f)&=\binom nr^{-1}\sum_{S\in\binom{[n]}r}f(D_n(S)),&
\tld\nu_n(f)&=\binom{N_n}r^{-1}\sum_{S\in\binom{\tld V_n}r}f(\tld D_n(S)),
\end{align*}
while otherwise $\tld\nu_n(f)=f(0)$ by convention. Hence
$T_n(f)=\sqrt{np_n}(\widetilde\nu_n(f)-\nu_n(f))$ (see \eqref{eq:process-def}).

To separate the fluctuations in the sampled vertex count, edge density,
and codegrees, write
\begin{align*}
N_n &= np_n(1+\delta),
& \tld q_n &= q_n(1+\omega),\\
\tld d_S &= np_nq_n^r\{D_n(S)+\Delta_S\},
& & S\in\binom{[n]}r.
\end{align*}

\begin{lemma}\label[lemma]{lem:normalizer_codegree_estimates}
Define the \emph{good event}
\begin{align*}
    \mathcal E_n:=\left\{|N_n/(np_n)-1|\le1/2,\ |\tld q_n/q_n-1|\le1/2,\max_{1\le t\le r}
    \|D_{\widetilde W_n}^{(t)}\ind_{\Delta_{N_n,t}}\|_\infty\le M\right\},
\end{align*}
and suppose that $G_n$ satisfies
\cref{ass:weights,ass:bounded-codegrees}, $np_n\to\infty$, and $\inf_n np_nq_n>0$. Then
\begin{align*}
    N_n=\left(1+O_{L^2}((np_n)^{-1/2})\right)np_n,\quad \tld q_n=\left(1+O_{L^2}((np_n)^{-1/2})\right)q_n.
\end{align*}
If, in addition, $np_nq_n^r/\log(np_n)\to\infty$, then $\bP(\mathcal E_n^c)\lesssim(np_n)^{-1}$.
\end{lemma}

\begin{proof}
It follows directly from $N_n\sim \operatorname{Bin}(n,p_n)$ that $N_n=\left(1+O_{L^2}((np_n)^{-1/2})\right)np_n$. Setting
\begin{align*}
    \widetilde{e}_n:=\frac{1}{n^2p_n^2q_n}\sum_{i,j=1}^nG_n(i,j)X_iX_j, \text{ note that }\E[\widetilde{e}_n]=\frac{1}{n^2q_n}\sum_{i,j=1}^nG_n(i,j).
\end{align*}
Since terms corresponding to disjoint edges are independent, a
covariance calculation gives
\begin{align*}
    \E[(\tld e_n-1)^2]=\Var(\tld e_n)
    &\le\frac{2n^2p_n^2q_n+4p_n^3\sum_{i=1}^nd_i^2}
               {(n^2p_n^2q_n)^2}\\
    &\le\frac{2n^2p_n^2q_n+4M^2n^3p_n^3q_n^2}{(n^2p_n^2q_n)^2}
       \le\frac4{np_n}\left(\frac1{np_nq_n}+M^2\right),
\end{align*}
where $d_i\le Mnq_n$. The normalizations satisfy
\begin{align*}
\frac{\tld q_n}{q_n}=\left(\frac{np_n}{N_n}\right)^2\tld e_n=(1+\delta)^{-2}\tld e_n,
\end{align*}
provided $\{N_n>0\}$. On $\{|N_n/(np_n)-1|\le1/2\}$ we have $|\tld q_n/q_n-1|\le4|\tld e_n-1|+6|N_n/(np_n)-1|$.
On its complement, $0\le\widetilde q_n\le1$ gives
$|\tld q_n/q_n-1|\le q_n^{-1}$. The binomial Chernoff bound gives
$\bP(|N_n/(np_n)-1|>1/2)\le2e^{-np_n/12}$. Consequently,
\begin{align*}
    \E[(\tld q_n/q_n-1)^2]
    &\lesssim\E[(\tld e_n-1)^2]+\E[(N_n/(np_n)-1)^2]+q_n^{-2}e^{-np_n/12}
    \\ &\lesssim(np_n)^{-1}+(np_n)^2e^{-np_n/12}
    \lesssim(np_n)^{-1}.
\end{align*}
Here we used $\inf_n np_nq_n>0$.
This proves $\tld q_n=\left(1+O_{L^2}((np_n)^{-1/2})\right)q_n$.

Suppose now that $np_nq_n^r/\log(np_n)\to\infty$. Set $M_0:=\sup_n\max_{1\le t\le r}\|D_{W_n}^{(t)}\ind_{\Delta_{n,t}}\|_\infty<M$, and choose
$\varepsilon_0\in(0,1/2]$ sufficiently small that
$(M_0+\varepsilon_0)/(1-\varepsilon_0)^{r+1}<M$.
For $1\le t\le r$ and $S\in\binom{[n]}t$,
$\Var(\tld d_S)\le\E[\tld d_S]\le M_0np_nq_n^t$.
Bernstein's inequality therefore gives, for some $c>0$,
\begin{align*}
    \max_{S\in\binom{[n]}t}
    \bP\bigl(\tld d_S>np_nq_n^t(M_0+\varepsilon_0)\bigr)
    \le e^{-c np_nq_n^t}.
\end{align*}
The absence of self-loops makes $X_S=\prod_{j\in S}X_j$ independent of $\tld d_S=\sum_{j\notin S}G_n(S,j)X_j$.
Since $q_n^t\ge q_n^r$ for $t\le r$, a union bound gives
\begin{align*}
    \bP\left(
    \max_{1\le t\le r}
    \max_{S\in\binom{[n]}t:\,X_S=1}
    \frac{\tld d_S}{np_nq_n^t}>M_0+\varepsilon_0\right)
    \le\sum_{t=1}^r\binom nt p_n^te^{-c np_nq_n^t}
    \lesssim(np_n)^{-1}.
\end{align*}
Moreover, the $L^2$ estimates just proved and Chebyshev's inequality give
\begin{align*}
    \bP(|N_n/(np_n)-1|>\varepsilon_0)
    +\bP(|\tld q_n/q_n-1|>\varepsilon_0)\lesssim(np_n)^{-1}.
\end{align*}
Thus, with probability $1-O((np_n)^{-1})$, we have
$|N_n/(np_n)-1|,|\tld q_n/q_n-1|\le\varepsilon_0$ and
$\tld d_S\le np_nq_n^t(M_0+\varepsilon_0)$ for every $1\le t\le r$ and
$S\in\binom{[n]}t$ with $X_S=1$.
On this event,
\begin{align*}
    \tld D_n(S)
    \le\frac{M_0+\varepsilon_0}{(1-\varepsilon_0)^{t+1}}
    \le\frac{M_0+\varepsilon_0}{(1-\varepsilon_0)^{r+1}}<M,
\end{align*}
which completes the proof.
\end{proof}

\begin{lemma}\label[lemma]{lem:second_moment_bias_bound}
Suppose that $G_n$ satisfies
\cref{ass:weights,ass:bounded-codegrees},
$np_n\to\infty$, $\inf_n np_nq_n>0$, and let $T_n$ be as defined in \cref{subsec:setup}.
Then, uniformly over $f\in C^1([0,M])$,
\begin{align*}
    \Var[T_n(f)]\lesssim\|f\|_{C^1([0,M])}^2.
\end{align*}
If, in addition, $G_n$ and $p_n$ satisfy
\cref{ass:sparsity}, then, for every fixed
$\beta\in(\alpha_*,1]$, uniformly over
$f\in C^{1,\beta}([0,M])$,
\begin{align*}
    |\E[T_n(f)]|
    \lesssim_\beta
    \|f\|_{C^{1,\beta}([0,M])}.
\end{align*}
If $\liminf_{n\to\infty} q_n>0$,
the same bound holds with $\beta=0$, i.e., with respect to the $C^1$ norm.
\end{lemma}

\begin{proof}
All implicit constants below are independent of $n$ and $f$.
Define the proxy
\begin{align}
    \overline T_n(f)
    :=\sqrt{np_n}\left\{
    \frac1{p_n^r\binom nr}
    \sum_{S\in\binom{[n]}r}X_S
    f\left(\frac{\tld d_S}{np_nq_n^r}\right)-\nu_n(f)
    \right\}.
    \label{eq:proxy_Tn}
\end{align}
Adding and subtracting the sampled average with deterministic
codegree normalization gives, when $N_n\ge r$ and $\tld q_n>0$,
\begin{align*}
T_n(f)-\overline T_n(f)
={}&
\sqrt{np_n}
\left(1-\frac{\binom{N_n}{r}}{p_n^r\binom nr}\right)
\frac1{\binom{N_n}{r}}
\sum_{S\in\binom{\tld V_n}r}
f\left(\frac{\tld d_S}{np_nq_n^r}\right)
\\
&+\frac{\sqrt{np_n}}{\binom{N_n}{r}}
\sum_{S\in\binom{\tld V_n}r}
\left[
f\bigl(\tld D_n(S)\bigr)
-f\left(\frac{\tld d_S}{np_nq_n^r}\right)
\right].
\end{align*}
Set $a:=\frac{N_n}{np_n}(\frac{\tld q_n}{q_n})^r$, so that
$\tld d_S/(np_nq_n^r)=a\tld D_n(S)$.
Since $f$ is extended by the convention $f(x)=f(x\wedge M)$,
boundedness and Lipschitz continuity give
\begin{align*}
\sup_{x\ge0}|f(ax)-f(x)|
\lesssim_M
\|f\|_{C^1([0,M])}\min\{1,|a-1|\}.
\end{align*}
Applying this bound to each summand in the second term
of the preceding identity, and bounding the average in
the first term by $\|f\|_\infty$, yields
\begin{align*}
|T_n(f)-\overline T_n(f)|
\lesssim_M
\sqrt{np_n}\,\|f\|_{C^1([0,M])}
\left[
\left|1-\frac{\binom{N_n}{r}}{p_n^r\binom nr}\right|
+\min\{1,|a-1|\}
\right].
\end{align*}
Finally, the elementary inequality
$\min\{1,|uv^r-1|\}\lesssim_r|u-1|+|v-1|$,
applied with $u=N_n/(np_n)$ and $v=\tld q_n/q_n$, gives
\begin{align*}
|T_n(f)-\overline T_n(f)|
\lesssim_{r,M}
\sqrt{np_n}\,\|f\|_{C^1([0,M])}
\left[
\left|1-\frac{\binom{N_n}{r}}{p_n^r\binom nr}\right|
+\left|\frac{N_n}{np_n}-1\right|
+\left|\frac{\tld q_n}{q_n}-1\right|
\right].
\end{align*}
Since $\binom{N_n}{r}=\sum_{S\in\binom{[n]}r}X_S$, and $\E[X_S]=p_n^r,\ \Cov(X_S,X_T)=p_n^{2r-|S\cap T|}-p_n^{2r}$ for all $S,T\in\binom{[n]}{r}$, expanding the variance and grouping pairs
by their intersection size gives
\begin{align*}
    \E\left[\left(\frac{\binom{N_n}{r}}{p_n^r\binom nr}-1\right)^2\right]
    &=\frac1{p_n^{2r}\binom nr^2}
       \sum_{S,T\in\binom{[n]}r}
           (p_n^{2r-|S\cap T|}-p_n^{2r})\\
    &=\sum_{k=1}^r
       \frac{\binom rk\binom{n-r}{r-k}}{\binom nr}(p_n^{-k}-1)\\
    &\lesssim_r\sum_{k=1}^r(np_n)^{-k}\lesssim_r(np_n)^{-1}.
\end{align*}
Combining this estimate with the $L^2$ normalizer bounds in
\cref{lem:normalizer_codegree_estimates} proves
\begin{align}
\E|T_n(f)-\overline T_n(f)|^2
\lesssim \|f\|_{C^1([0,M])}^2.
\label{eq:deterministic-normalizer-comparison}
\end{align}
Hence, it suffices to establish the variance and bias bounds
for $\overline T_n(f)$.

\medskip
\noindent\emph{Variance bound.}
Let $(X_1',\ldots,X_n')$ be an independent copy of
$(X_1,\ldots,X_n)$, and let $\overline T_n^{(i)}(f)$ be obtained by
replacing $X_i$ by $X_i'$. Write $X_S^{(i)}$ and $\tld d_S^{(i)}$ for the
corresponding versions of $X_S$ and $\tld d_S$. Then
\begin{align}
    \overline T_n(f)-\overline T_n^{(i)}(f)
    =\frac{\sqrt{np_n}}{p_n^r\binom nr}(H_i+L_i),
    \label{eq:single-site-two-contributions}
\end{align}
where
\begin{align*}
    H_i:=\sum_{S\ni i}(X_S-X_S^{(i)})
    f\left(\frac{\tld d_S}{np_nq_n^r}\right),
    \qquad
    L_i:=\sum_{S\not\ni i}X_S
    \left\{f\left(\frac{\tld d_S}{np_nq_n^r}\right)
    -f\left(\frac{\tld d_S^{(i)}}{np_nq_n^r}\right)\right\}.
\end{align*}
Indeed, the absence of self-loops gives $\tld d_S=\tld d_S^{(i)}$ when $i\in S$,
whereas $X_S=X_S^{(i)}$ when $i\notin S$. Since
\begin{align*}
\sum_{\substack{S\in\binom{[n]}r\\S\ni i}}
|X_S-X_S^{(i)}|
=|X_i-X_i'|
\sum_{T\in\binom{[n]\setminus\{i\}}{r-1}}
\prod_{j\in T}X_j
=|X_i-X_i'|
\binom{\sum_{j\ne i}X_j}{r-1},
\end{align*}
independence and binomial moment bounds give
\begin{align*}
    \E H_i^2
    \le\|f\|_\infty^2\E|X_i-X_i'|
    \E\left[\binom{\sum_{j\ne i}X_j}{r-1}^2\right]
    \lesssim \|f\|_\infty^2p_n(np_n)^{2(r-1)}.
\end{align*}
For $S\not\ni i$, Lipschitz continuity gives
\begin{align*}
\left|f\left(\frac{\tld d_S}{np_nq_n^r}\right)
    -f\left(\frac{\tld d_S^{(i)}}{np_nq_n^r}\right)\right|
&\le |\widetilde{d}_S-\widetilde{d}_S^{(i)}|
\frac{\|f'\|_\infty}{np_nq_n^r}\\
&\le \frac{\|f'\|_\infty}{np_nq_n^r}
|X_i-X_i'| G_n(S,i).
\end{align*}
Since $\sum_{S\not\ni i}X_SG_n(S,i)\le \widetilde{d}_i^{r}$, we have

\begin{align*}
\E L_i^2
&\le
\frac{\|f'\|_\infty^2}{(np_nq_n^r)^2}
\E[(X_i-X_i')^2]\,
\E\Big[\Big(\sum_{S\not\ni i}X_SG_n(S,i)\Big)^2\Big]\\
&\le \frac{2p_n\|f'\|_\infty^2}{(np_nq_n^r)^2} \E[(\tld d_{\{i\}})^{2r}]
\lesssim
\frac{\|f'\|_\infty^2}{(np_nq_n^r)^2}p_n(np_nq_n)^{2r}
=
\|f'\|_\infty^2p_n(np_n)^{2(r-1)},
\end{align*}
where the first and second inequalities in the last line above use the fact that $\tld d_{\{i\}}=\sum_{j\ne i}G_n(j,i)X_j$ is independent of $X_i$ and $X_i'$ and is a sum of independent $[0,1]$-valued random variables with mean at most $Mnp_nq_n$, respectively. Substituting the bound of $\E H_i^2$ and $\E L_i^2$ into
\eqref{eq:single-site-two-contributions} and applying the Efron--Stein
inequality yields
\begin{align*}
    \Var[\overline T_n(f)]
    \le\frac12\sum_{i=1}^n
    \E|\overline T_n(f)-\overline T_n^{(i)}(f)|^2
    \lesssim\|f\|_{C^1([0,M])}^2.
\end{align*}
Together with \eqref{eq:deterministic-normalizer-comparison}, this proves
the variance bound.

\medskip
\noindent\emph{Bias bound.}
Suppose now that \cref{ass:sparsity} holds.
Let $\beta\in (\alpha_*,1]$, and let
$f\in C^{1,\beta}([0,M])$. For
$S\in\binom{[n]}{r}$, set $\Delta_S:=(\tld d_S-p_nd_S)/(np_nq_n^r)$.
Then $\E[\Delta_S]=0$ and
$\Var(\tld d_S)\le p_n d_S\le Mnp_nq_n^r$,
where the last inequality follows from \cref{ass:bounded-codegrees}.
Thus, by Lyapunov's inequality,
\begin{align*}
    \E|\Delta_S|^{1+\beta}
    \le \frac{(\Var(\tld d_S))^{(1+\beta)/2}}{(np_nq_n^r)^{1+\beta}}
    \le M^{(1+\beta)/2}(np_nq_n^r)^{-(1+\beta)/2}.
\end{align*}
Let $M_0$ be as in the proof of \cref{lem:normalizer_codegree_estimates}. Then $D_n(S)=d_S/(nq_n^r)\le M_0<M$ for all $S$, and uniformly over $x\in[0,M_0]$ and $z\ge-x$, the extension of $f(x)=f(x\wedge M)$ satisfies
\begin{align}
    |f(x+z)-f(x)-f'(x)z|
    \lesssim_\beta
    \|f\|_{C^{1,\beta}([0,M])}|z|^{1+\beta}.
    \label{eq:bias-taylor-bound}
\end{align}
Indeed, this is the usual Taylor estimate when $x+z\le M$.
Otherwise, $z\ge M-M_0$, and the Lipschitz bound gives the same estimate, with a constant depending on $M-M_0>0$. Since $\tld d_S$ is independent of $X_S$ and
$\tld d_S/(np_nq_n^r)=D_n(S)+\Delta_S$, by \eqref{eq:proxy_Tn}, we have
\begin{align*}
    \E[\overline T_n(f)]
    =\frac{\sqrt{np_n}}{\binom nr}
    \sum_{S\in\binom{[n]}r}
    \E[f(D_n(S)+\Delta_S)-f(D_n(S))].
\end{align*}

Applying the Taylor remainder bound to this identity and taking
expectations, the linear terms vanish since $\E[\Delta_S]=0$
for every $S$. The preceding bound on $\E|\Delta_S|^{1+\beta}$
then gives
\begin{align*}
    |\E[\overline T_n(f)]|
    \lesssim_\beta
    \sqrt{np_n}\,
    \|f\|_{C^{1,\beta}([0,M])}
    (np_nq_n^r)^{-(1+\beta)/2}.
\end{align*}
Together with
\eqref{eq:deterministic-normalizer-comparison},
this gives
\begin{align*}
    |\E[T_n(f)]|
    \lesssim_\beta
    \|f\|_{C^{1,\beta}([0,M])}
    \left[
        1+(np_n)^{-\beta/2}q_n^{-r(1+\beta)/2}
    \right].
\end{align*}
Since $$\beta>\alpha_*=\frac{r\gamma_*}{1-r\gamma_*}\Rightarrow \frac{\beta}{1+\beta}>{r\gamma_*},$$ we can
choose $\varepsilon>0$
sufficiently small that
$r(\gamma_*+\varepsilon)(1+\beta)<\beta$.
By the definition of $\gamma_*$ (see \cref{ass:sparsity}),
for all sufficiently large $n$,
\begin{align}
    \frac{\log(1/q_n)}{\log(np_n)}
    \le \gamma_*+\varepsilon,
    \label{eq:exponent}
\end{align}
and hence $q_n\ge (np_n)^{-(\gamma_*+\varepsilon)}$.
Therefore,
\begin{align*}
    (np_n)^{-\beta/2}q_n^{-r(1+\beta)/2}
    &\le
    (np_n)^{-\{\beta-r(\gamma_*+\varepsilon)(1+\beta)\}/2}
    \longrightarrow 0.
\end{align*}
If $\liminf_{n\to\infty}q_n>0$, the bound in
\eqref{eq:bias-taylor-bound} also holds with
$\beta=0$, and the term in brackets above is bounded.
\end{proof}

\subsection{Finite-dimensional convergence}

We next establish finite-dimensional convergence of $T_n$. For polynomial
test functions, the fluctuations reduce to those of sampled biclique counts
and the normalizing vertex and edge counts. The extension to general
$C^{1,\alpha}$ test functions uses polynomial approximation in a weaker
H\"older norm. The following lemma is a standard consequence of Bernstein
polynomial approximation; see, e.g., \citet[Ch.~10]{DeVoreLorentz1993}.

\begin{lemma}\label[lemma]{lem:poly_approx_weaker_holder}
Let $0\le\beta<\alpha\le1$. For every
$f\in C^{1,\alpha}([0,M])$, there exist polynomials $P_\ell$ such that
\(
    \|P_\ell-f\|_{C^{1,\beta}([0,M])}
    \overset{\ell\to\infty}{\longrightarrow}0.
\)
\end{lemma}

We will also use the following auxiliary lemma.

\begin{lemma}
\label[lemma]{lem:limiting-influence-bound}
For every $f\in C^1([0,M])$,
\begin{align*}
    \|\phi_f\|_\infty
    \lesssim_{r,M}
    \|f\|_{C^1([0,M])}.
\end{align*}
\end{lemma}

\begin{proof}
By \eqref{eq:limiting-influence-function} and
$0\le D_W^{(t)}\le M$ almost everywhere,
\begin{align*}
|\phi_f(x)|
\le 2r\|f\|_\infty+
\{M^r+M+2rM(M+1)\}\|f'\|_\infty
\end{align*}
for almost every $x$. Taking the essential supremum proves the claim.
\end{proof}

\begin{lemma}\label[lemma]{lem:fdd_convergence}
Suppose that the graph sequence $G_n$ satisfies \cref{ass:weights,ass:star-cut-convergence,ass:bounded-codegrees}, the sampling
ratios $(p_n)$ satisfy \cref{ass:sampling-ratios}, and the two sequences jointly satisfy \cref{ass:sparsity}. Let $\left\{G(f):f\in C^1([0,M])\right\}$ be as in \cref{def:Gf_brownian_int}, and fix $\alpha_*<\alpha\le1$. If $f_1,\ldots,f_m\in C^{1,\alpha}([0,M])$, then
\begin{align*}
    \bigl(T_n(f_1),\ldots,T_n(f_m)\bigr)
    \dto
    \bigl(G(f_1),\ldots,G(f_m)\bigr),
\end{align*}
If $\liminf_n q_n>0$, the same conclusion holds for arbitrary
$f_1,\ldots,f_m\in C^1([0,M])$.
\end{lemma}

\begin{proof}
By linearity of $T_n$, bilinearity of $\Gamma$, and the
Cram\'er--Wold device, it is enough to show, for each fixed $f\in C^{1,\alpha}([0,M])$, that
\begin{align*}
    T_n(f)\dto G(f).
\end{align*}
\medskip
\noindent\emph{Polynomial test functions.}
We first consider the monomial $P_k(x):=x^k$. If $k=0$, then
$T_n(P_0)=0$. Moreover, \eqref{eq:limiting-influence-function} gives $\phi_{P_0}\equiv0$, so
$\Gamma(P_0,P_0)=0$. Suppose now that $k\ge1$. Fix $\varepsilon\in(0,1/(2r)-\gamma_*)$, which exists since $\alpha_*<1$ is equivalent to $\gamma_*<1/(2r)$. Let
$\mathcal E_n$ be as in \cref{lem:normalizer_codegree_estimates}. Since
\cref{ass:sparsity} implies $np_nq_n^r/\log(np_n)\to\infty$, that lemma gives $\bP(\mathcal E_n)\to1$.

On $\mathcal E_n$, we have $N_n>r$ for all sufficiently large $n$, and so $\widetilde \nu_n$ is well defined.  Moreover, by the definition of $\mathcal E_n$ we have $\max_{\substack{S\in\binom{[n]}r\\X_S=1}}\tld D_n(S)\le M$. Invoking
\cref{lem:counting-lemma}\emph{(i)} with $G_n$ replaced by $\widetilde G_n$ on the event $\mathcal E_n$ (see \cref{rem:deterministic-counting-comparisons}) gives
\begin{align*}
    \left|
        \widetilde\nu_n(P_k)-t_{\mathrm{inj}}(K_{r,k},\widetilde G_n)
    \right|
    \lesssim_{r,k,M}\frac{1}{N_n\widetilde q_n^r},
\end{align*}
Multiplying both sides by $\sqrt{np_n}$ and noting that $N_n=\Theta(np_n), \widetilde q_n=\Theta(q_n)$ on $\mathcal E_n$ gives
\begin{align*}
    &\sqrt{np_n}\left|
        \widetilde\nu_n(P_k)-t_{\mathrm{inj}}(K_{r,k},\widetilde G_n)
    \right|
    \lesssim_{r,k,M}\left(\sqrt{np_n}\,q_n^r\right)^{-1}
    \le (np_n)^{-1/2+r(\gamma_*+\varepsilon)}
    \underset{\text{\cref{ass:sparsity}}}{\overset{n\to\infty}{\longrightarrow}}0,
\end{align*}
where we use the fact that $q_n\ge (np_n)^{-\gamma_*-\varepsilon}$ for all large $n$, by \eqref{eq:exponent}.
Together with $\bP(\mathcal E_n)\to1$, this proves
\begin{align*}
    \widetilde\nu_n(P_k)=t_{\mathrm{inj}}(K_{r,k},\widetilde G_n)
    +o_P((np_n)^{-1/2}).
\end{align*}
For the parent graph $G_n$, \cref{lem:counting-lemma}\emph{(i)} and the same argument give
\begin{align*}
    \nu_n(P_k)=t_{\mathrm{inj}}(K_{r,k},G_n)+o((np_n)^{-1/2}).
\end{align*}
Combining these two estimates yields
\begin{align}
    T_n(P_k)=\sqrt{np_n}\left(
    t_{\mathrm{inj}}(K_{r,k},\widetilde G_n)
    -t_{\mathrm{inj}}(K_{r,k},G_n)\right)+o_P(1).
    \label{eq:monomial-biclique-reduction}
\end{align}
Therefore, \cref{cor:biclique-clt} and
\eqref{eq:monomial-biclique-reduction} give joint convergence for any finite
collection of monomials. Since $T_n$ and the map $f\mapsto\phi_f$ are linear,
the same conclusion holds for every finite collection of polynomials. In
particular, for every polynomial $P$, we have $T_n(P)\dto G(P)$.

\medskip
\noindent\emph{Polynomial approximation.}

Fix $\beta\in(\alpha_*,\alpha)$.
Since $f\in C^{1,\alpha}([0,M])$,
\cref{lem:poly_approx_weaker_holder} gives
polynomials $(P_\ell)$ such that
\begin{align*}
    \|P_\ell-f\|_{C^{1,\beta}([0,M])}
    \overset{\ell\to\infty}{\longrightarrow}0.
\end{align*}
Set $R_\ell:=f-P_\ell$. By
\cref{lem:second_moment_bias_bound},
\begin{align*}
    &\limsup_{n\to\infty}\E\left[|T_n(R_\ell)|^2\right]
    \lesssim
    \|R_\ell\|_{C^1([0,M])}^2
    +\|R_\ell\|_{C^{1,\beta}([0,M])}^2
    \overset{\ell\to\infty}{\longrightarrow}0.
\end{align*}
Also,  noting that $G(f)-G(P_\ell)=G(f-P_\ell)$ (see \cref{def:Gf_brownian_int}),
we get
\begin{align*}
    \E\left[|G(f)-G(P_\ell)|^2\right]
    =\Gamma(R_\ell,R_\ell)
    =(1-p_\infty)\int_0^1\phi_{R_\ell}(x)^2\,dx
    \lesssim\|R_\ell\|_{C^1([0,M])}^2
    \overset{\ell\to\infty}{\longrightarrow}0,
\end{align*}
where the last inequality uses
\cref{lem:limiting-influence-bound} with $f=R_\ell$.
By the approximation argument, it follows that $T_n(f)\dto G(f)$. 
If $\liminf_n q_n>0$, then
$\gamma_*=\alpha_*=0$, and the endpoint case of
\cref{lem:second_moment_bias_bound} gives the
bias bound $|\E[T_n(f)]|\lesssim \|f\|_{C^1([0,M])}$.
The same approximation argument then applies to arbitrary $f\in C^1([0,M])$, since polynomials
are dense in $C^1([0,M])$ by the Weierstrass
approximation theorem applied to $f'$ and
integration. This completes the proof.
\end{proof}

\subsection{Tightness}

\begin{proposition}
\label[proposition]{prop:tightness_holder}
Suppose that the graph sequence $G_n$ satisfies
\cref{ass:weights,ass:star-cut-convergence,ass:bounded-codegrees},
the sampling ratios $(p_n)$ satisfy \cref{ass:sampling-ratios},
and the two sequences jointly satisfy \cref{ass:sparsity}.
If $\max\{\alpha_*,1/2\}<\alpha\le1$, then $T_n$ is tight in $\ell^\infty(\mathcal B_{1,\alpha})$.
\end{proposition}

We first control the change in $T_n(f)$ under resampling and
deduce a sub-Gaussian concentration bound for its centered version
on $\mathcal E_n$. A chaining argument based on Dudley's entropy
inequality then gives asymptotic $C^1$-equicontinuity.
Finally, uniform control of the bias yields tightness of $T_n$.

For $x\in\{0,1\}^n$, let $T_n(f;x)$ denote the value of
$T_n(f)$ under the sampling configuration $x$, and write
\begin{align*}
\tld V(x)&:=\{i:x_i=1\},
&
N_n(x)&:=|\tld V(x)|=\sum_{i=1}^n x_i,\\
\tld d_S(x)&:=\sum_{v\in \tld V(x)}G_n(S,v),
&
\widetilde q_n(x)&:=
\frac{1}{\{1\vee N_n(x)\}^2}
\sum_{i,j\in \tld V(x)}G_n(i,j).
\end{align*}
We equip $\{0,1\}^n$ with the Hamming distance
$d_H(x,x'):=|\tld V(x)\triangle \tld V(x')|=\sum_{i=1}^n |x_i-x_i'|$.

\begin{lemma}\label[lemma]{lem:holder-process-hamming-stability}
Suppose that $G_n$ satisfies \cref{ass:weights}.
Let $\mathcal E_n\subseteq\{0,1\}^n$ be as in \cref{lem:normalizer_codegree_estimates}. If $np_n\ge8r$, then, for every $f\in C^1([0,M])$ and
$x,x'\in\mathcal E_n$,
\begin{align*}
|T_n(f;x)-T_n(f;x')|
\lesssim_{r,M}
(np_n)^{-1/2}\|f\|_{C^1([0,M])}d_H(x,x').
\end{align*}
\end{lemma}

\begin{proof}
Fix $x,x'\in\mathcal E_n$, and set
$d:=d_H(x,x')$ and $I:=\tld V(x)\cap \tld V(x')$.
For $z\in\{x,x'\}$, write
\begin{align}
    \widetilde D_n(S;z)
    &:=\frac{\tld d_S(z)}{N_n(z)\widetilde q_n(z)^r},
       \quad S\in\binom{\tld V(z)}r,\label{eq:tldDn}\\
    \widetilde\nu_n(f;z)
    &:=\frac1{\binom{N_n(z)}r}
       \sum_{S\in\binom{\tld V(z)}r}f(\widetilde D_n(S;z)).\notag
\end{align}
Then
\begin{align}
T_n(f;x)-T_n(f;x')
=
\sqrt{np_n}\{\widetilde\nu_n(f;x)
-\widetilde\nu_n(f;x')\}.
\label{eq:hamming-stability-basic-identity}
\end{align}
If $d\ge np_n/4$, the result follows from
$|\widetilde\nu_n(f;x)-\widetilde\nu_n(f;x')|
\le2\|f\|_\infty$.
Suppose that $d<np_n/4$.

Since $x,x'\in\mathcal E_n$, for $z\in\{x,x'\}$,
\begin{align}
\frac{np_n}{2}\le N_n(z)\le\frac{3np_n}{2},
\qquad
\frac{q_n}{2}\le\widetilde q_n(z)\le\frac{3q_n}{2},
\qquad
\tld d_S(z)\le M N_n(z)\widetilde q_n(z)^t
\label{eq:good-event-uniform-comparisons}
\end{align}
for $S\in\binom{\tld V(z)}t$ and $1\le t\le r$.
In particular, \begin{align}\label{eq:double_dip}
|I|\ge N_n(x)-d\ge np_n/4.
\end{align}

Define the averages over the common tuples by
\begin{align*}
\widetilde\nu_n^I(f;z)
:=
\binom{|I|}r^{-1}
\sum_{S\in\binom Ir}f(\widetilde D_n(S;z)),
\qquad z\in\{x,x'\}.
\end{align*}
We use the decomposition
\begin{align}
|\widetilde\nu_n(f;x)-\widetilde\nu_n(f;x')|
\le{}&
|\widetilde\nu_n(f;x)-\widetilde\nu_n^I(f;x)|+
|\widetilde\nu_n^I(f;x)-\widetilde\nu_n^I(f;x')|
\notag\\
&+
|\widetilde\nu_n^I(f;x')-\widetilde\nu_n(f;x')|.
\label{eq:hamming-stability-decomposition}
\end{align}
We first bound the middle term in
\eqref{eq:hamming-stability-decomposition}.
Canceling the edges within $I$ and using symmetry gives
\begin{align*}
N_n(x)^2\widetilde q_n(x)-N_n(x')^2\widetilde q_n(x')
&=\sum_{i,j\in\tld V(x)} G_n(i,j) - \sum_{i,j\in\tld V(x')} G_n(i,j)\\
&=
2\sum_{i\in\tld V(x)\setminus I}\sum_{j\in I}G_n(i,j)
+\sum_{i,j\in\tld V(x)\setminus I}G_n(i,j)\\
&\quad-
2\sum_{i\in\tld V(x')\setminus I}\sum_{j\in I}G_n(i,j)
-\sum_{i,j\in\tld V(x')\setminus I}G_n(i,j).
\end{align*}
By nonnegativity and the $t=1$ case of
\eqref{eq:good-event-uniform-comparisons},
\begin{align*}
|N_n(x)^2\widetilde q_n(x)-N_n(x')^2\widetilde q_n(x')|
&\le
2\sum_{i\in\tld V(x)\setminus I}\tld d_{\{i\}}(x)
+2\sum_{i\in\tld V(x')\setminus I}\tld d_{\{i\}}(x')\\
&\le
2Md\left(N_n(x)\widetilde q_n(x)
+N_n(x')\widetilde q_n(x')
\right)\lesssim_M dnp_nq_n,
\end{align*}
since $|\tld V(x)\setminus I|+|\tld V(x')\setminus I|=d$.
To pass from edge counts to edge densities, write
\begin{align*}
\widetilde q_n(x)-\widetilde q_n(x')
&=
\frac{N_n(x)^2\widetilde q_n(x)-N_n(x')^2\widetilde q_n(x')}
     {N_n(x)^2}
+\widetilde q_n(x')
\frac{N_n(x')^2-N_n(x)^2}{N_n(x)^2}.
\end{align*}
Using $|N_n(x)-N_n(x')|\le d$ and
\eqref{eq:good-event-uniform-comparisons}, we obtain
\begin{align}
|\widetilde q_n(x)-\widetilde q_n(x')|
\lesssim_M q_n\frac d{np_n},\quad
\left|
\frac{N_n(x')\widetilde q_n(x')^r}
     {N_n(x)\widetilde q_n(x)^r}-1
\right|
\lesssim_{r,M}\frac d{np_n}.
\label{eq:good-event-denominator-stability}
\end{align}

For $S\in\binom Ir$, we have
\begin{align*}
\tld d_S(x)-\tld d_S(x')
&=
\sum_{i\in\tld V(x)\setminus I}G_n(S,i)
-\sum_{i\in\tld V(x')\setminus I}G_n(S,i).
\end{align*}
Also, by \eqref{eq:tldDn},
\begin{align*}
\widetilde D_n(S;x)-\widetilde D_n(S;x')
&=
\frac{\tld d_S(x)-\tld d_S(x')}
     {N_n(x)\widetilde q_n(x)^r}
     +
     \widetilde D_n(S;x')
\left\{
\frac{N_n(x')\widetilde q_n(x')^r}
     {N_n(x)\widetilde q_n(x)^r}-1
\right\}.
\end{align*}
Since $\widetilde D_n(S;x')\le M$,
\eqref{eq:good-event-denominator-stability} gives
\begin{align*}
|\widetilde D_n(S;x)-\widetilde D_n(S;x')|
\lesssim_{r,M}
\frac{1}{np_nq_n^r}
\sum_{i\in\tld V(x)\triangle\tld V(x')}G_n(S,i)
+\frac d{np_n}.
\end{align*}
For each $i\in\tld V(x)\triangle\tld V(x')$, let
$z\in\{x,x'\}$ be the configuration containing $i$.
Since $I\subseteq\tld V(z)$,
\begin{align*}
\sum_{S\in\binom Ir}G_n(S,i)
&\le
\left(\sum_{v\in I}G_n(i,v)\right)^r
\le \bigl(\tld d_{\{i\}}(z)\bigr)^r
\lesssim_{r,M}(np_nq_n)^r,
\end{align*}
by \eqref{eq:good-event-uniform-comparisons}.
Moreover, $|I|\ge np_n/4$ by \eqref{eq:double_dip} and
$np_n\ge8r$ by hypothesis, so
$\binom{|I|}{r}\gtrsim_r(np_n)^r$.
Averaging over $S\in\binom Ir$ and interchanging the sums,
we therefore obtain
\begin{align*}
\frac{1}{\binom{|I|}{r}}
\sum_{S\in\binom Ir}
|\widetilde D_n(S;x)-\widetilde D_n(S;x')|\lesssim_{r,M}
\frac{d(np_nq_n)^r}{np_nq_n^r\binom{|I|}{r}}
+\frac d{np_n}\lesssim_{r,M}\frac d{np_n}.
\end{align*}
The mean value theorem now gives
\begin{align}
|\widetilde\nu_n^I(f;x)-\widetilde\nu_n^I(f;x')|
\lesssim_{r,M}\frac d{np_n}\|f'\|_\infty,
\label{eq:hamming-stability-common-tuples}
\end{align}
which bounds the middle term in
\eqref{eq:hamming-stability-decomposition}. Proceeding to control the first and third terms in the RHS of \eqref{eq:hamming-stability-decomposition}, for $z\in\{x,x'\}$, we have $N_n(z)-|I|\le d$ (see \eqref{eq:double_dip}) and
$N_n(z)-r\ge \frac{np_n}{2}-\frac{np_n}{8}=3np_n/8$. Hence
\begin{align*}
1-\frac{\binom{|I|}{r}}{\binom{N_n(z)}r}
=
1-\prod_{j=0}^{r-1}\frac{|I|-j}{N_n(z)-j}
\le
\sum_{j=0}^{r-1}\frac{N_n(z)-|I|}{N_n(z)-j}
\lesssim_r\frac d{np_n},
\end{align*}
where we used
$1-\prod_j(1-a_j)\le\sum_j a_j$ for $a_j\in[0,1]$.
Therefore,
\begin{align}
|\widetilde\nu_n(f;z)-\widetilde\nu_n^I(f;z)|
\le
2\left\{1-\frac{\binom{|I|}{r}}{\binom{N_n(z)}r}\right\}
\|f\|_\infty
\lesssim_r\frac d{np_n}\|f\|_\infty.
\label{eq:hamming-stability-index-sets}
\end{align}
Combining \eqref{eq:hamming-stability-decomposition},
\eqref{eq:hamming-stability-common-tuples}, and
\eqref{eq:hamming-stability-index-sets} along with
\eqref{eq:hamming-stability-basic-identity} proves the result.
\end{proof}

\begin{lemma}\label[lemma]{lem:holder-process-concentration}
Set
$S_n(f):=T_n(f)-\E[T_n(f)]$, for $f\in C^1([0,M])$. 
Suppose that the graph $G_n$ satisfies
\cref{ass:weights,ass:bounded-codegrees}, and further we have $np_n\to\infty$ and
$np_nq_n^r/\log(np_n)\to\infty$.
Then there exist $C<\infty$ and $n_0\in\bN$ such that,
for every $n\ge n_0$, $f\in C^1([0,M])$, and $s\ge0$,
\begin{align}
\bP\left(
\mathcal E_n\cap
\left\{|S_n(f)|>C\sqrt s\,\|f\|_{C^1}\right\}
\right)
\le2e^{-s}.
\label{eq:holder-process-subgaussian-tail}
\end{align}
\end{lemma}

\begin{proof}
The claim is immediate if $f=0$.
Suppose that $f\ne0$, and take $n$ sufficiently large that
$np_n\ge8r$.
Lemma \ref{lem:holder-process-hamming-stability} gives the existence of a finite positive constant $c=c_{r,M}$, such that for all $f\in C^1([0,M])$ and $x,x'\in \mathcal{E}_n$ we have
$$|T_n(f;x)-T_n(f;x')|\le c (np_n)^{-1/2} \|f\|_{C^1([0,M])} d_H(x,x')=L d_H(x,x'),$$
where 
$L:=c(np_n)^{-1/2}\|f\|_{C^1}$. By McShane's extension theorem (followed by truncation)
there is an $L$-Lipschitz function
$F:\{0,1\}^n\to\mathbb R$
equal to $T_n(f;\cdot)$ on $\mathcal E_n$ and satisfying
$\|F\|_\infty\le2\sqrt{np_n}\|f\|_\infty$.
By the variance-sensitive bounded-differences inequality
of \citet[Corollary~6]{Warnke2016}, which follows from Freedman's
martingale inequality, applied to $F$ and $-F$ we obtain
\begin{align*}
\bP\left(
|F(X)-\E F(X)|>t
\right)
\le
2\exp\left\{
-\frac{ct^2}
{\|f\|_{C^1}^2+(np_n)^{-1/2}\|f\|_{C^1}t}
\right\},
\qquad t\ge0.
\end{align*}
For $0\le t\le4\sqrt{np_n}\|f\|_{C^1}$, the denominator
is at most $5\|f\|_{C^1}^2$.
For larger $t$, the probability on the left is zero,
as $\|F\|_\infty\le 2\sqrt{np_n}\|f\|_\infty$.
Thus, after adjusting $c>0$ if needed, we have
\begin{align}
\bP\left(
|F(X)-\E F(X)|>t
\right)
\le
2\exp\left\{-\frac{ct^2}{\|f\|_{C^1}^2}\right\},
\qquad t\ge0.
\label{eq:extended-process-subgaussian-tail}
\end{align}

Since $F(X)=T_n(f)$ on $\mathcal E_n$ and both
variables are bounded in absolute value by
$2\sqrt{np_n}\|f\|_\infty$,
\cref{lem:normalizer_codegree_estimates} gives
\begin{align*}
|\E F(X)-\E T_n(f)|
\le
4\sqrt{np_n}\|f\|_\infty\bP(\mathcal E_n^c)
\lesssim
(np_n)^{-1/2}\|f\|_{C^1}.
\end{align*}
On $\mathcal E_n$, we have
\begin{align*}
S_n(f)
=
F(X)-\E F(X)
+\E F(X)-\E T_n(f).
\end{align*}
For $s>\log2$, apply
\eqref{eq:extended-process-subgaussian-tail} with
$t=(C/2)\sqrt s\,\|f\|_{C^1}$.
Choosing $C$ sufficiently large absorbs the centering
difference and makes the exponent at least $s$.
For $0\le s\le\log2$, the desired bound is automatic
because $2e^{-s}\ge1$.
This proves \eqref{eq:holder-process-subgaussian-tail}.
\end{proof}

\begin{lemma}\label[lemma]{lem:centered-process-tightness}
Suppose that $G_n$ satisfies
\cref{ass:weights,ass:bounded-codegrees} and $np_n\to\infty$.
Assume further that $np_nq_n^r/\log(np_n)\to\infty$.
If $1/2<\alpha\le1$, then $S_n$ is tight in
$\ell^\infty(\mathcal B_{1,\alpha})$.
\end{lemma}

\begin{proof}
Equip $\mathcal B_{1,\alpha}$ with the metric
$d_{C^1}(f,g):=\|f-g\|_{C^1([0,M])}$.
Since $(f-g)/2\in\mathcal B_{1,\alpha}$ whenever
$f,g\in\mathcal B_{1,\alpha}$, linearity reduces asymptotic
$d_{C^1}$-equicontinuity to proving
\begin{align}
\E\left[
\sup_{\substack{f\in\mathcal B_{1,\alpha}\\
                 \|f\|_{C^1}\le\delta}}
|S_n(f)|\ind_{\mathcal E_n}
\right]
\lesssim \delta^{1-1/(2\alpha)},
\qquad 0<\delta\le1,
\label{eq:stopped-process-equicontinuity}
\end{align}
uniformly for all sufficiently large $n$, since
$\bP(\mathcal E_n^c)\to0$ by
\cref{lem:normalizer_codegree_estimates}. Writing $N(\varepsilon,T,d)$ for the covering number of
$(T,d)$, \citet[Theorem~2.7.1]{vdVW1996}, applied to the
derivative class and combined with a discretization of $f(0)$,
gives
\begin{align}
\log N(\varepsilon,\mathcal B_{1,\alpha},d_{C^1})
\lesssim\varepsilon^{-1/\alpha},
\qquad 0<\varepsilon\le1.
\label{eq:C1-holder-ball-entropy}
\end{align}
Thus $\mathcal B_{1,\alpha}$ is totally bounded and thus separable.
Moreover, linearity and the bound
$|S_n(f)|\le4\sqrt{np_n}\|f\|_\infty$ imply that the sample
paths are $d_{C^1}$-continuous, so the suprema below are measurable.

Since $\mathcal E_n$ does not depend on $f$, linearity and
\cref{lem:holder-process-concentration}, applied to $f-g$, give
\begin{align*}
\bP\left(
\left|
S_n(f)\ind_{\mathcal E_n}
-S_n(g)\ind_{\mathcal E_n}
\right|
>C\sqrt{s}\,d_{C^1}(f,g)
\right)
\le 2e^{-s},
\end{align*}
for every $f,g\in\mathcal B_{1,\alpha}$ and $s\ge0$,
with $C$ independent of $n,f,g,s$.
As $S_n(0)=0$, Dudley's entropy inequality
\citep[Corollary~2.2.8]{vdVW1996} yields
\begin{align*}
\E\left[
\sup_{\substack{f\in\mathcal B_{1,\alpha}\\
                 \|f\|_{C^1}\le\delta}}
|S_n(f)|\ind_{\mathcal E_n}
\right]
&\lesssim
\int_0^\delta
\sqrt{\log N(\varepsilon,\mathcal B_{1,\alpha},d_{C^1})}
\,d\varepsilon\\
&\lesssim
\int_0^\delta\varepsilon^{-1/(2\alpha)}\,d\varepsilon
\lesssim \delta^{1-1/(2\alpha)},
\end{align*}
where we used \eqref{eq:C1-holder-ball-entropy} and
$\alpha>1/2$.
This proves \eqref{eq:stopped-process-equicontinuity}
and hence asymptotic $d_{C^1}$-equicontinuity.

Finally, the variance bound in
\cref{lem:second_moment_bias_bound} gives tightness of the
finite-dimensional marginals. The conclusion follows from
\citet[Theorem~1.5.7]{vdVW1996}.
\end{proof}

\begin{proof}[Proof of \cref{prop:tightness_holder}]
Since $np_nq_n^r/\log(np_n)\to\infty$, the process $S_n$ is tight by
\cref{lem:centered-process-tightness}.
It remains to show that $T_n$ and $S_n$ are asymptotically
equivalent.

Fix $\beta\in(\alpha_*,\alpha)$ and set $m_n(f):=\E[T_n(f)]$. By linearity and \cref{lem:second_moment_bias_bound},
\begin{align*}
|m_n(f)-m_n(g)|
\lesssim
\|f-g\|_{C^{1,\beta}([0,M])}.
\end{align*}
For each fixed $f\in\mathcal B_{1,\alpha}$,
\cref{lem:fdd_convergence} gives convergence of $T_n(f)$
to a centered Gaussian law.
The variance and bias bounds in
\cref{lem:second_moment_bias_bound} also give
$\sup_n\E[|T_n(f)|^2]<\infty$, hence uniform integrability.
It follows that $m_n(f)\to0$. Since $\alpha>\beta$, the ball $\mathcal B_{1,\alpha}$
is compact in the $C^{1,\beta}$ topology. By the Arzel\`a--Ascoli theorem, the preceding equicontinuity bound and pointwise convergence imply
\begin{align*}
\|T_n-S_n\|_\infty
=
\sup_{f\in\mathcal B_{1,\alpha}}|m_n(f)|
\overset{n\to\infty}{\longrightarrow}0.
\end{align*}
This completes the proof.
\end{proof}

The finite-dimensional convergence in \cref{lem:fdd_convergence}
and the tightness in \cref{prop:tightness_holder} prove
\cref{thm:empirical-codegree-process}.

\section{Proofs of the sharpness results}
\label{sec:sharpness-proofs}

We prove the two sharpness results stated in
\cref{subsec:sharpness-results}. Both constructions concern ordinary
degrees under sampling with $p_n\equiv1/2$.

\subsection{Failure of finite-dimensional convergence below \texorpdfstring{$\alpha_*$}{alpha-star}}

\begin{proof}[Proof of \cref{prop:regularity-counterexample}]
Set
$d^{(n)}:=2\lfloor n^{1-\gamma}/2\rfloor$,
so that $d^{(n)}=(1+o(1))n^{1-\gamma}$.
Fix $\eta\in(0,1-\gamma)$. For sufficiently large $n$,
we have $n^\eta\le d^{(n)}\le n/2$.
If $G_n$ is uniformly distributed over the simple
$d^{(n)}$-regular graphs on $[n]$, then by
\citet[Theorem~A]{TikhomirovYoussef2019} there exists
$C_\eta\in(0,\infty)$ such that
\begin{align*}
    \bP\left(
        \lambda(G_n)\le C_\eta\sqrt{d^{(n)}}
    \right)
    \ge 1-\frac1n,
\end{align*}
where
$\lambda(G_n):=\max\{|\lambda_2(A_{G_n})|,
|\lambda_n(A_{G_n})|\}$,
with the adjacency eigenvalues ordered as
$d^{(n)}=\lambda_1(A_{G_n})\ge\cdots\ge\lambda_n(A_{G_n})$.
Hence, for every sufficiently large $n$, we may choose
a deterministic $d^{(n)}$-regular simple graph $G_n$
such that
\begin{align}
    \lambda(G_n)\le C_\eta\sqrt{d^{(n)}}.
    \label{eq:regularity-counterexample-spectral-gap}
\end{align}

By the expander mixing lemma
\cite[see, e.g., Section~2.4]{HooryLinialWigderson2006},
for all $S,T\subseteq[n]$,
\begin{align*}
    \left|
        \sum_{i\in S}\sum_{j\in T}
        \left(
            G_n(i,j)-\frac{d^{(n)}}{n}
        \right)
    \right|
    \le
    \lambda(G_n)\sqrt{|S||T|}.
\end{align*}
Since $q_n=d^{(n)}/n$, the inequality above and
\eqref{eq:regularity-counterexample-spectral-gap} imply
\begin{align*}
    \|W_n-1\|_\square
    &=
    \frac{1}{n^2q_n}
    \max_{S,T\subseteq[n]}
    \left|
        \sum_{i\in S}\sum_{j\in T}
        \left(
            G_n(i,j)-\frac{d^{(n)}}{n}
        \right)
    \right|\\
    &\le
    \frac{\lambda(G_n)}{d^{(n)}}
    \longrightarrow 0.
\end{align*}
Thus \cref{ass:weights,ass:star-cut-convergence,ass:bounded-codegrees}
hold with $W\equiv1$, while
\cref{ass:sampling-ratios,ass:sparsity} follow from
$p_n=1/2$, $\gamma<1/2$, and
\begin{align*}
    (np_n)^\gamma q_n
    =
    \frac{d^{(n)}}{2^\gamma n^{1-\gamma}}
    \longrightarrow 2^{-\gamma}.
\end{align*}
In particular, $\gamma_*=\gamma$.

Set
$s_n:=\sqrt{np_n}(d^{(n)})^{-(1+\alpha)/2}$,
and let $\overline T_n$ be as in \eqref{eq:proxy_Tn}.
The sampled degree
$\tld d_i=\sum_j G_n(i,j)X_j$ of every vertex
$i\in[n]$ is $\Bin(d^{(n)},p_n)$-distributed
and is independent of
$X_i=\ind\{i\in V(\tld G_n)\}$.
Taking $f(x)=|x-1|^{1+\alpha}$, we have
$f(d_i/(nq_n))=f(1)=0$ for every $i\in[n]$.
Consequently,
\begin{align}
    \E[\overline T_n(f)]
    =
    \sqrt{np_n}\,
    \E\left|
        \frac{\Bin(d^{(n)},p_n)}
             {p_n d^{(n)}}-1
    \right|^{1+\alpha}.
    \label{eq:regularity-counterexample-proxy-bias}
\end{align}
Let
\begin{align*}
    Z_n\sim
    \frac{2\Bin(d^{(n)},p_n)-d^{(n)}}
         {\sqrt{d^{(n)}}}.
\end{align*}
Since $p_n=1/2$, the binomial central limit theorem
gives $Z_n\dto Z$, where $Z$ is standard normal,
and $\E[Z_n^2]=1$.
As $1+\alpha<2$, the sequence
$(|Z_n|^{1+\alpha})_n$ is uniformly integrable.
Therefore,
\begin{align}
    \frac{\E[\overline T_n(f)]}{s_n}
    =
    \E|Z_n|^{1+\alpha}
    \longrightarrow
    b_\alpha:=\E|Z|^{1+\alpha}>0.
    \label{eq:regularity-counterexample-averaged-bias}
\end{align}

It follows that, for all sufficiently large $n$,
$\E[\overline T_n(f)]\ge (b_\alpha/2)s_n$.
By Cauchy--Schwarz,
\begin{align*}
    \E[T_n(f)]
    &=
    \E[\overline T_n(f)]
    +\E[T_n(f)-\overline T_n(f)]\\
    &\ge
    \E[\overline T_n(f)]
    -
    \left\{
        \E|T_n(f)-\overline T_n(f)|^2
    \right\}^{1/2}.
\end{align*}
By the variance bound in
\cref{lem:second_moment_bias_bound}
and \eqref{eq:deterministic-normalizer-comparison},
there exists $C_f<\infty$, independent of $n$,
such that, for sufficiently large $n$,
\begin{align}
    \E[T_n(f)]
    \ge \frac{b_\alpha}{2}s_n-C_f,
    \qquad
    \Var[T_n(f)]\le C_f.
    \label{eq:regularity-counterexample-bounds}
\end{align}

Since $d^{(n)}=(1+o(1))n^{1-\gamma}$,
\begin{align*}
    s_n
    &=
    \sqrt{np_n}(d^{(n)})^{-(1+\alpha)/2}\\
    &=
    (2^{-1/2}+o(1))
    n^{(\gamma-\alpha(1-\gamma))/2}.
\end{align*}
Moreover,
\begin{align*}
    \gamma-\alpha(1-\gamma)
    =
    (1-\gamma)(\alpha_*-\alpha)>0,
\end{align*}
because $\alpha<\alpha_*=\gamma/(1-\gamma)$.
Thus $s_n\to\infty$, and
\eqref{eq:regularity-counterexample-bounds}
gives $\E[T_n(f)]\to\infty$.
For every fixed $K>0$, Chebyshev's inequality gives
\begin{align*}
    \bP\{T_n(f)\le K\}
    &\le
    \frac{\Var[T_n(f)]}
         {(\E[T_n(f)]-K)^2}
    \longrightarrow 0.
\end{align*}
Hence $T_n(f)\pto+\infty$.
\end{proof}

\subsection{Failure of tightness below
\texorpdfstring{$1/2$}{1/2}}

\begin{proof}[Proof of \cref{prop:degree-tightness-counterexample}]
\emph{Constructing the limiting graphon and parent graphs.}
Choose $1<\beta<1/(\alpha+1/2)$, and partition $[0,1)$ into consecutive
intervals $A_j$, $j\ge1$, with
\begin{align*}
    |A_j|=:a_j
    =\frac{j^{-\beta}}{\sum_{k\ge1}k^{-\beta}}
    =\Theta(j^{-\beta}).
\end{align*}
For $j\in\bN$, let $B_j:=1/2+A_j$, and $d_j:=(\inf B_j+\sup B_j)/2$ be its midpoint. Define
\begin{align}
    W(x,y):=
    \sum_{j\ge1}\frac{d_j}{a_j}
    \ind_{A_j}(x)\ind_{A_j}(y).
    \label{eq:graphon-tightness-counterexample}
\end{align}
Then $D_W(\cdot):=\int W(\cdot,y)\ dy=d_j$ on $A_j$, so $1/2<D_W<3/2$ almost everywhere. Moreover, using that $\{A_j\}$ partition $[0,1)$ gives
\begin{align*}
    \int W(x,y)\ dxdy
    =
    \sum_{j\ge1}a_jd_j
    =
    \sum_{j\ge1} \left(\frac{a_j}{2} +\int_{A_j}x\ dx\right)
    =
    \frac12+\int_0^1x\ dx
    =
    1.
\end{align*}

Applying \cref{thm:sparse-graphon-samples} to $W$ with target densities
$\rho_n=n^{-\gamma}, r=1$ and $p_n=p_\infty=\frac{1}{2}$, it follows that the random graph $G_n$ satisfies \cref{ass:weights,ass:star-cut-convergence,ass:bounded-codegrees,ass:sampling-ratios,ass:sparsity} almost surely (see \cref{rem:check}). We may therefore fix a realization of simple unweighted parent graphs satisfying these conditions. For such realization, $q_n/\rho_n\to1$, and since $\rho_n=n^{-\gamma}$, we have $\gamma_*=\gamma$. Therefore, since $\gamma/(1-\gamma)=\alpha_*<\alpha$,
\cref{lem:fdd_convergence} gives finite-dimensional convergence of $T_n$.

\medskip
\noindent\emph{Test functions.}
Fix a $C^\infty$ function $\psi$ supported in $(-1/4,1/4)$ with $\psi(0)=0$ and $\psi'(0)=1$, and set
\begin{align*}
    f_j(x):=a_j^{1+\alpha}
    \psi\left(\frac{x-d_j}{a_j}\right),
    \qquad 0\le x\le M.
\end{align*}

Since $d_j$ is the midpoint of $B_j$, and $|A_j|=a_j$, then
\begin{align*}
\supp(f_j)\subseteq
\left(d_j-\frac{a_j}{4},d_j+\frac{a_j}{4}\right)
\subseteq B_j.
\end{align*}

We claim that the signed sums of the test functions are uniformly bounded in
$C^{1,\alpha}$:
\begin{align}
    C:=\sup_{\substack{m\in\bN\\\theta\in\{-1,1\}^m}}
    \left\|\sum_{j=1}^m\theta_jf_j\right\|_{C^{1,\alpha}}<\infty.
    \label{eq:tightness-signed-test-functions}
\end{align}
To verify this, fix $m\ge1$ and $\theta\in\{-1,1\}^m$, and write
$F=\sum_{j=1}^m\theta_jf_j$. Since the supports of the $f_j$ are
pairwise disjoint, we have
$\|F\|_\infty\le\|\psi\|_\infty$ and
$\|F'\|_\infty\le\|\psi'\|_\infty$.
It remains to bound $[F']_\alpha$ uniformly in $m$ and $\theta$.
Using $a_j\le 1$, we obtain $[f_j']_\alpha\le[\psi']_\alpha$
(see \eqref{eq:holder-norms}). Since the intervals $B_j$ are
pairwise disjoint and $\supp(f_j)\subset B_j$, we have
$F=\theta_j f_j$ on $B_j$ and $F=0$ outside
$\bigcup_{j=1}^m B_j$.

To bound $[F']_\alpha$, fix $x,y\in[0,M]$ and consider
the following cases. If $x,y\in B_j$ for some $j$, then
$F'=\theta_j f_j'$ on $B_j$, so
$|F'(x)-F'(y)|\le[f_j']_\alpha|x-y|^\alpha\le[\psi']_\alpha|x-y|^\alpha$.

If $x\in B_j$ and $y\in B_k$ for distinct indices,
say $j<k$, let $z_j:=\sup B_j$ and $z_k:=\inf B_k$ denote the facing endpoints of $B_j$ and $B_k$, respectively. Both lie between $x$ and $y$, and
$f_j'(z_j)=f_k'(z_k)=0$. Hence
\begin{align*}
    |F'(x)-F'(y)|
    &\le |f_j'(x)-f_j'(z_j)|
       +|f_k'(y)-f_k'(z_k)|\\
    &\le [\psi']_\alpha
       \bigl(|x-z_j|^\alpha+|y-z_k|^\alpha\bigr)\le 2[\psi']_\alpha|x-y|^\alpha.
\end{align*}
Finally, if exactly one of $x,y$ belongs to
$\bigcup_{j=1}^m B_j$, the same argument applies
with only one term, using the endpoint of that interval
between $x$ and $y$. If neither belongs to the union,
both derivatives vanish. Thus
$[F']_\alpha\le2[\psi']_\alpha$, proving the claim.

\medskip
\noindent\emph{Gaussian fluctuations.}
Since $D_W=d_k$ on $A_k$ and $\supp(f_j)\subset B_j$,
where the intervals $B_k$ are pairwise disjoint, we have
$f_j(d_k)=f_j'(d_k)=0$ for $k\ne j$.
On the other hand, $\psi(0)=0$ and $\psi'(0)=1$ give
$f_j(d_j)=0$ and $f_j'(d_j)=a_j^\alpha$.
Therefore,
\begin{align}
    f_j(D_W)=0,
    \qquad
    f_j'(D_W)=a_j^\alpha\ind_{A_j}
    \quad\text{a.e. on }[0,1].
    \label{eq:tightness-local-test-functions}
\end{align}

Moreover, by the definition of $W$,
\begin{align*}
\int_0^1 f_j'(D_W(y))W(y,x)\,dy
&=a_j^\alpha\int_{A_j}W(y,x)\,dy
=d_ja_j^\alpha\ind_{A_j}(x),\\
\int_0^1 D_W(y)f_j'(D_W(y))\,dy
&=d_ja_j^{1+\alpha}.
\end{align*}
For $r=1$, the first term in
\eqref{eq:limiting-influence-function} vanishes since
$f_j(D_W(x))=0$. The derivative term, however, is
$d_ja_j^\alpha\ind_{A_j}(x)$, whose integral over
$[0,1]$ is $d_ja_j^{1+\alpha}$. Since $f_j(D_W)=0$, the first term in
\eqref{eq:limiting-influence-function} vanishes.
Using the identities above, together with
$|A_j|=a_j$ and $\int_0^1D_W(x)\,dx=1$, we obtain
\begin{align*}
    \phi_{f_j}(x)
    &=d_ja_j^\alpha(\ind_{A_j}(x)-a_j)
      -2d_ja_j^{1+\alpha}(D_W(x)-1)\\
    &=d_ja_j^\alpha\ind_{A_j}(x)
      -d_ja_j^{1+\alpha}(2D_W(x)-1).
\end{align*}

Let $B$ be the Brownian motion, and put
\begin{align*}
    Z_j:=a_j^{-1/2}\int_{A_j}dB(x),
    \qquad
    V:=\int_0^1(2D_W(x)-1)\,dB(x).
\end{align*}
Then $(Z_j)_{j\ge1}$ are independent standard normal random variables.
By \cref{def:Gf_brownian_int} and the expression for $\phi_{f_j}$ above,
\begin{align}
    G(f_j)
    &=\sqrt{1-p_\infty}\,d_j
    \left(a_j^\alpha\int_{A_j}dB(x)-a_j^{\alpha+1}V\right)
    \notag\\
    &=\sqrt{1-p_\infty}\,d_j
    \left(a_j^{\alpha+1/2}Z_j-a_j^{\alpha+1}V\right).
    \label{eq:tightness-block-gaussians}
\end{align}
Set $S_m:=\sum_{j=1}^m d_ja_j^{\alpha+1/2}|Z_j|$.
Since $\E|Z_j|=\sqrt{2/\pi}$, $1/2<d_j<3/2$,
$a_j=\Theta(j^{-\beta})$, and $\beta(\alpha+1/2)<1$, we have
\begin{align*}
    \E S_m
    &=\sqrt{\frac2\pi}\sum_{j=1}^m d_ja_j^{\alpha+1/2}
    \gtrsim\sum_{j=1}^m j^{-\beta(\alpha+1/2)}
    \longrightarrow\infty.
\end{align*}
On the other hand, independence of the $Z_j$ and
$\Var(|Z_j|)=1-2/\pi$ give, for every $m\in\bN$,
\begin{align*}
    \Var(S_m)
    =\left(1-\frac2\pi\right)
      \sum_{j=1}^m d_j^2a_j^{2\alpha+1}
    \le\frac94\sum_{j=1}^m a_j
    \le\frac94,
\end{align*}
where we used $d_j<3/2$, $0<a_j\le1$, and $\alpha>0$. It follows from Chebyshev's inequality that $S_m\pto+\infty$ as
$m\to\infty$. Also, since $\beta(1+\alpha)>1$, and $a_j=\Theta(j^{-\beta})$, the random variable $R:=|V|\sum_{j\ge1}d_ja_j^{1+\alpha}$ is almost surely finite. Therefore \eqref{eq:tightness-block-gaussians} gives
\begin{align}
    \sum_{j=1}^m|G(f_j)|
    \ge\sqrt{1-p_\infty}\,(S_m-R)
    \pto+\infty.
    \label{eq:tightness-gaussian-sum-diverges}
\end{align}

For fixed $m$, with $C$ as in \eqref{eq:tightness-signed-test-functions}, define
\begin{align*}
    \mathcal C_m
    :=\left\{
    C^{-1}\sum_{j=1}^m\theta_jf_j:
    \theta\in\{-1,1\}^m
    \right\}
    \subset\mathcal B_{1,\alpha}.
\end{align*}

By linearity, $\max_{f\in\mathcal C_m}|T_n(f)|=C^{-1}\sum_{j=1}^m|T_n(f_j)|$, and since $\alpha>\alpha_*$, \cref{thm:empirical-codegree-process} and the continuous mapping theorem give, for every fixed $m$,
\begin{align*}
    C^{-1}\sum_{j=1}^m|T_n(f_j)|
    \dto C^{-1}\sum_{j=1}^m|G(f_j)|.
\end{align*}
For any $K>0$, since $\mathcal C_m\subset\mathcal B_{1,\alpha}$,
the preceding convergence and the Portmanteau theorem give
\begin{align*}
    \liminf_{n\to\infty}
    \bP\left(\sup_{f\in\mathcal B_{1,\alpha}}|T_n(f)|>K\right)
    \ge
    \bP\left(C^{-1}\sum_{j=1}^m|G(f_j)|>K\right)
\end{align*}
for every fixed $m$. By
\eqref{eq:tightness-gaussian-sum-diverges}, the right-hand side
tends to one as $m\to\infty$. Hence
\begin{align*}
    \sup_{f\in\mathcal B_{1,\alpha}}|T_n(f)|
    \pto+\infty.
\end{align*}

Finally, for every fixed $f\in C^{1,\alpha}([0,M])$,
\cref{lem:fdd_convergence,lem:second_moment_bias_bound} and uniform
integrability give $\E[T_n(f)]\overset{n\to\infty}{\longrightarrow}0$. Hence $T_n-\E T_n$ has the same finite-dimensional limits, and the preceding argument applies to the centered process as well.
\end{proof}

\begin{acks}[Acknowledgments]
The authors used generative AI tools to assist with refining the exposition and preparing figures using TikZ. All AI-assisted material was reviewed by the authors, who take full responsibility for the manuscript.
\end{acks}

\bibliographystyle{imsart-nameyear}
\bibliography{references}

@article{bdm22,
  author    = {Bhattacharya, Bhaswar B. and Das, Sayan and Mukherjee, Sumit},
  title     = {Motif estimation via subgraph sampling: The fourth-moment phenomenon},
  journal   = {The Annals of Statistics},
  year      = {2022},
  volume    = {50},
  number    = {2},
  pages     = {987--1011},
  doi       = {10.1214/21-AOS2134}
}

@article{MR458659,
  author    = {Frank, Ove},
  title     = {Estimation of graph totals},
  journal   = {Scandinavian Journal of Statistics},
  year      = {1977},
  volume    = {4},
  number    = {2},
  pages     = {81--89},
  issn      = {0303-6898},
  mrnumber  = {0458659}
}

@article{MR515656,
  author    = {Frank, Ove},
  title     = {Estimation of the number of connected components in a graph by using a sampled subgraph},
  journal   = {Scandinavian Journal of Statistics},
  year      = {1978},
  volume    = {5},
  number    = {4},
  pages     = {177--188},
  issn      = {0303-6898},
  mrnumber  = {0515656}
}

@incollection{Frank_2005,
  author    = {Frank, Ove},
  title     = {Network sampling and model fitting},
  booktitle = {Models and Methods in Social Network Analysis},
  editor    = {Carrington, Peter J. and Scott, John and Wasserman, Stanley},
  series    = {Structural Analysis in the Social Sciences},
  publisher = {Cambridge University Press},
  address   = {Cambridge},
  year      = {2005},
  pages     = {31--56},
  doi       = {10.1017/CBO9780511811395.003}
}

@inproceedings{pmlr-v75-klusowski18a,
  author    = {Klusowski, Jason M. and Wu, Yihong},
  title     = {Counting Motifs with Graph Sampling},
  booktitle = {Proceedings of the 31st Conference on Learning Theory},
  series    = {Proceedings of Machine Learning Research},
  volume    = {75},
  pages     = {1966--2011},
  year      = {2018},
  editor    = {Bubeck, S{\'e}bastien and Perchet, Vianney and Rigollet, Philippe},
  publisher = {PMLR}
}

@article{KlusowskiWu2020,
  author    = {Klusowski, Jason M. and Wu, Yihong},
  title     = {Estimating the number of connected components in a graph via subgraph sampling},
  journal   = {Bernoulli},
  year      = {2020},
  volume    = {26},
  number    = {3},
  pages     = {1635--1664},
  doi       = {10.3150/19-BEJ1147}
}

@article{ZhangKolaczykSpencer2015,
  author    = {Zhang, Yaonan and Kolaczyk, Eric D. and Spencer, Bruce D.},
  title     = {Estimating network degree distributions under sampling: An inverse problem, with applications to monitoring social media networks},
  journal   = {The Annals of Applied Statistics},
  year      = {2015},
  volume    = {9},
  number    = {1},
  pages     = {166--199},
  doi       = {10.1214/14-AOAS800}
}

@article{10.1214/aoms/1177705148,
  author    = {Goodman, Leo A.},
  title     = {Snowball Sampling},
  journal   = {The Annals of Mathematical Statistics},
  year      = {1961},
  volume    = {32},
  number    = {1},
  pages     = {148--170},
  doi       = {10.1214/aoms/1177705148}
}

@article{BCLSV08,
  author    = {Borgs, Christian and Chayes, Jennifer T. and Lov{\'a}sz, L{\'a}szl{\'o} and S{\'o}s, Vera T. and Vesztergombi, Katalin},
  title     = {Convergent sequences of dense graphs {I}: Subgraph frequencies, metric properties and testing},
  journal   = {Advances in Mathematics},
  year      = {2008},
  volume    = {219},
  number    = {6},
  pages     = {1801--1851},
  doi       = {10.1016/j.aim.2008.07.008}
}

@book{Lovasz2012,
  author    = {Lov{\'a}sz, L{\'a}szl{\'o}},
  title     = {Large Networks and Graph Limits},
  series    = {American Mathematical Society Colloquium Publications},
  volume    = {60},
  publisher = {American Mathematical Society},
  address   = {Providence, RI},
  year      = {2012},
  isbn      = {978-0-8218-9085-1},
  doi       = {10.1090/coll/060}
}

@book{vdVW1996,
  author    = {van der Vaart, Aad W. and Wellner, Jon A.},
  title     = {Weak Convergence and Empirical Processes: With Applications to Statistics},
  series    = {Springer Series in Statistics},
  publisher = {Springer},
  address   = {New York},
  year      = {1996},
  doi       = {10.1007/978-1-4757-2545-2}
}

@article{Warnke2016,
  author  = {Warnke, Lutz},
  title   = {On the Method of Typical Bounded Differences},
  journal = {Combin. Probab. Comput.},
  year    = {2016},
  volume  = {25},
  number  = {2},
  pages   = {269--299},
  doi     = {10.1017/S0963548315000103}
}

@article{BCCZ2019LpI,
  author  = {Borgs, Christian and Chayes, Jennifer T. and Cohn, Henry and Zhao, Yufei},
  title   = {An {$L^p$} theory of sparse graph convergence {I}: Limits, sparse random graph models, and power law distributions},
  journal = {Transactions of the American Mathematical Society},
  year    = {2019},
  volume  = {372},
  number  = {5},
  pages   = {3019--3062},
  doi     = {10.1090/tran/7543}
}

@article{BCCZ2018LpII,
  author  = {Borgs, Christian and Chayes, Jennifer T. and Cohn, Henry and Zhao, Yufei},
  title   = {An {$L^p$} theory of sparse graph convergence {II}: {LD} convergence, quotients and right convergence},
  journal = {The Annals of Probability},
  year    = {2018},
  volume  = {46},
  number  = {1},
  pages   = {337--396},
  doi     = {10.1214/17-AOP1187}
}

@article{HladkyPelekisSileikis2021,
  author  = {Hladk{\'y}, Jan and Pelekis, Christos and {\v{S}}ileikis, Matas},
  title   = {A limit theorem for small cliques in inhomogeneous random graphs},
  journal = {Journal of Graph Theory},
  year    = {2021},
  volume  = {97},
  number  = {4},
  pages   = {578--599},
  doi     = {10.1002/jgt.22673}
}

@article{BhattacharyaChatterjeeJanson2023,
  author  = {Bhattacharya, Bhaswar B. and Chatterjee, Anirban and Janson, Svante},
  title   = {Fluctuations of subgraph counts in graphon based random graphs},
  journal = {Combinatorics, Probability and Computing},
  year    = {2023},
  volume  = {32},
  number  = {3},
  pages   = {428--464},
  doi     = {10.1017/S0963548322000335}
}

@article{ChatterjeeDanBhattacharya2024,
  author  = {Chatterjee, Anirban and Dan, Soham and Bhattacharya, Bhaswar B.},
  title   = {Higher-order graphon theory: Fluctuations, degeneracies and inference},
  journal = {The Annals of Statistics},
  year    = {2026},
  volume  = {54},
  number  = {4},
  pages   = {2136--2163},
  doi     = {10.1214/26-AOS2644}
}

@incollection{Meckes2009,
  author = {Meckes, Elizabeth S.},
  title = {On Stein's method for multivariate normal approximation},
  booktitle = {High Dimensional Probability V: The Luminy Volume},
  series = {Institute of Mathematical Statistics Collections},
  volume = {5},
  pages = {153--178},
  publisher = {Institute of Mathematical Statistics},
  year = {2009},
  doi       = {10.1214/09-IMSCOLL511}
}

@article{Raic2004Decomposable,
  author  = {Rai{\v c}, Martin},
  title   = {A Multivariate {CLT} for Decomposable Random Vectors with Finite Second Moments},
  journal = {Journal of Theoretical Probability},
  volume  = {17},
  number  = {3},
  pages   = {573--603},
  year    = {2004},
  doi     = {10.1023/B:JOTP.0000040290.44087.68}
}

@book{DeVoreLorentz1993,
  author    = {Ronald A. DeVore and George G. Lorentz},
  title     = {Constructive Approximation},
  series    = {Grundlehren der mathematischen Wissenschaften},
  volume    = {303},
  publisher = {Springer},
  address   = {Berlin},
  year      = {1993}
}

@book{Janson2013Graphons,
  author    = {Svante Janson},
  title     = {Graphons, Cut Norm and Distance, Couplings and Rearrangements},
  series    = {New York Journal of Mathematics Monographs},
  volume    = {4},
  publisher = {State University of New York, University at Albany},
  year      = {2013}
}

@book{Serfling1980,
  author    = {Serfling, Robert J.},
  title     = {Approximation Theorems of Mathematical Statistics},
  series    = {Wiley Series in Probability and Mathematical Statistics},
  publisher = {John Wiley \& Sons},
  address   = {New York},
  year      = {1980},
  doi       = {10.1002/9780470316481}
}

@article{DelmasDhersinSciauveau2021,
  author  = {Delmas, Jean-Fran{\c{c}}ois and Dhersin, Jean-St{\'e}phane and Sciauveau, Marion},
  title   = {Asymptotic for the cumulative distribution function of the degrees and homomorphism densities for random graphs sampled from a graphon},
  journal = {Random Structures \& Algorithms},
  year    = {2021},
  volume  = {58},
  number  = {1},
  pages   = {94--149},
  doi     = {10.1002/rsa.20965}
}

@article{Maugis2024,
  author  = {Maugis, Pierre-Andr{\'e} G.},
  title   = {Central limit theorems for local network statistics},
  journal = {Biometrika},
  year    = {2024},
  volume  = {111},
  number  = {3},
  pages   = {743--754},
  doi     = {10.1093/biomet/asad080}
}

@article{KloppVerzelen2019,
  author  = {Klopp, Olga and Verzelen, Nicolas},
  title   = {Optimal graphon estimation in cut distance},
  journal = {Probability Theory and Related Fields},
  year    = {2019},
  volume  = {174},
  number  = {3--4},
  pages   = {1033--1090},
  doi     = {10.1007/s00440-018-0878-1}
}

@misc{LyuMukherjee2024,
    author       = {Lyu, Hanbaek and Mukherjee, Sumit},
    title        = {Large random matrices with given margins},
    year         = {2024},
    eprint       = {2407.14942},
    archivePrefix = {arXiv},
    primaryClass = {math.PR},
     note         = {arXiv:2407.14942}
}

@article{TikhomirovYoussef2019,
  author  = {Tikhomirov, Konstantin and Youssef, Pierre},
  title   = {The spectral gap of dense random regular graphs},
  journal = {The Annals of Probability},
  volume  = {47},
  number  = {1},
  pages   = {362--419},
  year    = {2019},
  doi     = {10.1214/18-AOP1263}
}

@article{HooryLinialWigderson2006,
  author  = {Hoory, Shlomo and Linial, Nathan and Wigderson, Avi},
  title   = {Expander graphs and their applications},
  journal = {Bulletin of the American Mathematical Society},
  volume  = {43},
  number  = {4},
  pages   = {439--561},
  year    = {2006},
  doi     = {10.1090/S0273-0979-06-01126-8}
}

@article{Frank1980DegreeFrequencies,
  author  = {Frank, Ove},
  title   = {Estimation of the number of vertices of different degrees in a graph},
  journal = {Journal of Statistical Planning and Inference},
  year    = {1980},
  volume  = {4},
  number  = {1},
  pages   = {45--50},
  doi     = {10.1016/0378-3758(80)90032-4}
}

@article{HunterGoodreauHandcock2008,
  author  = {Hunter, David R. and Goodreau, Steven M. and Handcock, Mark S.},
  title   = {Goodness of fit of social network models},
  journal = {Journal of the American Statistical Association},
  year    = {2008},
  volume  = {103},
  number  = {481},
  pages   = {248--258},
  doi     = {10.1198/016214507000000446}
}

@article{Stewart2025EmpiricalDistributions,
  author  = {Stewart, Jonathan R.},
  title   = {Consistency of empirical distributions of sequences of graph statistics in networks with dependent edges},
  journal = {Journal of Multivariate Analysis},
  year    = {2025},
  volume  = {207},
  pages   = {105420},
  doi     = {10.1016/j.jmva.2025.105420}
}

@article{Janson2026UncoveringGraph,
  author  = {Janson, Svante},
  title   = {Uncovering a graph},
  journal = {Probability Theory and Related Fields},
  year    = {2026},
  volume  = {195},
  pages   = {689--740},
  doi     = {10.1007/s00440-025-01401-9}
}

@article{ChenKato2020,
  author  = {Chen, Xiaohui and Kato, Kengo},
  title   = {Jackknife multiplier bootstrap: finite sample approximations to the {$U$}-process supremum with applications},
  journal = {Probability Theory and Related Fields},
  year    = {2020},
  volume  = {176},
  number  = {3--4},
  pages   = {1097--1163},
  doi     = {10.1007/s00440-019-00936-y}
}

@article{BoistardLopuhaaRuizGazen2017,
  author  = {Boistard, H{\'e}l{\`e}ne and Lopuha{\"a}, Hendrik P. and Ruiz-Gazen, Anne},
  title   = {Functional central limit theorems for single-stage sampling designs},
  journal = {The Annals of Statistics},
  year    = {2017},
  volume  = {45},
  number  = {4},
  pages   = {1728--1758},
  doi     = {10.1214/16-AOS1507}
}

@article{HsiehHsuKoKovarikLogan2024,
  author  = {Hsieh, Chih-Sheng and
             Hsu, Yu-Chin and
             Ko, Stanley I. M. and
             Kov{\'a}{\v r}{\'i}k, Jarom{\'i}r and
             Logan, Trevon D.},
  title   = {Non-representative sampled networks:
             Estimation of network structural properties by weighting},
  journal = {Journal of Econometrics},
  year    = {2024},
  volume  = {240},
  number  = {1},
  pages   = {105689},
  doi     = {10.1016/j.jeconom.2024.105689}
}

@article{ConlonFoxZhao2014,
  author  = {Conlon, David and Fox, Jacob and Zhao, Yufei},
  title   = {Extremal Results in Sparse Pseudorandom Graphs},
  journal = {Advances in Mathematics},
  volume  = {256},
  pages   = {206--290},
  year    = {2014},
  doi     = {10.1016/j.aim.2013.12.004}
}

\appendix

\section{Verification of assumptions for random graphon samples}
\label{app:sparse-graphon-samples}

This appendix proves \cref{thm:sparse-graphon-samples} for the model in
\cref{subsec:graphon-examples}. Since $K_W^{(1)}=W$, the case $t=1$
of \cref{ass:star-cut-convergence} is the usual cut-norm convergence
of normalized graphon samples. For sparse samples, this is proved in
\citet[Theorem~2.14(b)]{BCCZ2019LpI} with target-density normalization.
We extend the argument to higher-order star kernels and establish the
uniform codegree bounds.

Use $I_i$, $\Delta_{n,t}$, and $\sigma_n$ as defined in the introduction.
For each $n$, suppress the superscript in $A_{ij}^{(n)}$ and write
$A_{ij}=A_{ij}^{(n)}$. We begin with an approximation lemma for the
latent variables.

For a positive integer $d$ and $F\in L^1([0,1]^d)$, let $\Pi_nF$ be the step function given by, for $s\in [n]^d$, and $\mathbf x\in I_{s(1)}\times\cdots\times I_{s(d)}$
\begin{align*}
    (\Pi_nF)(\mathbf x)
    :=
    \begin{cases}
        F(U_{\sigma_n(s(1))},\ldots,U_{\sigma_n(s(d))}),
        &s\in[n]_d,\\
        0,&s\notin[n]_d.
    \end{cases}
\end{align*}

\begin{lemma}\label[lemma]{lem:latent-empirical-approximation}
Let $d\in \bN$ be a fixed integer, $F\in L^1([0,1]^d)$, and let $F_n$ be a sequence of measurable functions on $[0,1]^d$ such that $\sup_n |F_n|\le \varphi$ for some $\varphi\in L^1([0,1]^d)$, and $F_n\to F$ almost everywhere, as $n\to\infty$. Then
\begin{align*}
    \Pi_n F_n\ \overset{n\to \infty}{\longrightarrow} F\text{ in }L^1([0,1]^d)
    \qquad\text{almost surely}.
\end{align*}
\end{lemma}

\begin{proof}
First assume $F_n=F$ for all $n$. Glivenko--Cantelli's theorem gives
\begin{align*}
    \delta_n:=\max_{1\le i\le n} \left|U_{\sigma_n(i)}-\frac{i}{n}\right|\overset{n\to\infty}{\longrightarrow} 0
    \qquad\text{almost surely}.
\end{align*}
Let $G\in C([0,1]^d)$, and denote its modulus of continuity by
$\omega_G$. If $x=(x_1,\dots,x_d)\in
I_{s(1)}\times\cdots\times I_{s(d)}$, then
$|x_a-s(a)/n|\le 1/n$ for $1\le a\le d$. Hence
\begin{align*}
    \|\Pi_n G-G\|_1\le \omega_G\left(\delta_n+\frac{1}{n}\right)+\|G\|_\infty\left(1-\frac{(n)_d}{n^d}\right)\overset{n\to\infty}{\longrightarrow} 0
    \qquad\text{almost surely}.
\end{align*}
Given $F\in L^1([0,1]^d),\ \varepsilon>0$, and a continuous function $G$ with $\|F-G\|_1\le \varepsilon$, Hoeffding's strong law for $U$-statistics \citep[Theorem~5.4A]{Serfling1980} yields
\begin{align*}
    \|\Pi_n(F-G)\|_1
    &=\frac1{n^d}\sum_{s\in[n]_d}
       |(F-G)(U_{s(1)},\ldots,U_{s(d)})|\\
    &\longrightarrow\|F-G\|_1\le\varepsilon
       \qquad\text{almost surely}.
\end{align*}
This extends the conclusion to every $F\in L^1([0,1]^d)$.
For the varying sequence, set
$\varphi_N:=\sup_{\ell\ge N}|F_\ell-F|$. Then
$0\le\varphi_N\le2\varphi$ almost everywhere, and for $n\ge N$,
\begin{align*}
    |\Pi_nF_n-\Pi_nF|\le\Pi_n|F_n-F|\le\Pi_n\varphi_N.
\end{align*}
The fixed-function result, applied to each $\varphi_N$, gives almost surely
\begin{align*}
    \limsup_{n\to\infty}\|\Pi_nF_n-\Pi_nF\|_1\le\|\varphi_N\|_1.
\end{align*}
Since $\varphi_N\to0$ almost everywhere, dominated convergence gives
$\|\varphi_N\|_1\to0$. Combining this with
$\|\Pi_nF-F\|_1\to0$ completes the proof.
\end{proof}

\begin{proof}[Proof of \cref{thm:sparse-graphon-samples}]
Set $p_{ij}:=1\wedge\rho_nW(U_i,U_j)$ for $i\ne j$ and
$p_{ii}:=0$. For $(x,y)\in I_i\times I_j$, define the graphon
$\overline W_n(x,y):=\rho_n^{-1}
A_{\sigma_n(i)\sigma_n(j)}$. Whenever $q_n>0$, we have
$W_n=(\rho_n/q_n)\overline W_n$.

Fix $1\le t\le r$ and $s=(s(1),\ldots,s(t))\in[n]_t$. Conditional on
$U_{s(1)},\ldots,U_{s(t)}$, the variables
$\{\prod_{a=1}^t p_{s(a),j}\}_{j\in[n]\setminus\bar s}$ are independent, take values in $[0,1]$, and
\begin{align*}
    \Var\left(\prod_{a=1}^t p_{s(a),j}\middle | U_{s(1)},\dots,U_{s(t)}\right)\le \E\left[\prod_{a=1}^t p_{s(a),j}\middle| U_{s(1)},\dots,U_{s(t)}\right]\le \|D_W^{(t)}\|_\infty\rho_n^t,
\end{align*}
for all $j\in [n]\setminus \bar s$. Hence, by
Bernstein's inequality and the law of total probability, there exists
$c=c(t,W)>0$ such that
\begin{align*}
    &\bP\left(
        \sum_{j\in [n]\setminus\bar s}\prod_{a=1}^t p_{s(a),j}>\left(1+\|D_W^{(t)}\|_\infty\right)n\rho_n^t
    \right)
    \le 2\exp\{-cn\rho_n^t\}.
\end{align*}
Since $n\rho_n^t\ge n\rho_n^r\gg\log n$, these probabilities remain
summable after taking a union bound over $s\in[n]_t$. Thus, by
Borel--Cantelli, almost surely, for all sufficiently large $n$,
\begin{align}
    \max_{s\in[n]_t}\sum_{j\in [n]\setminus\bar s}
    \prod_{a=1}^t p_{s(a),j}
    &\le
    \left(1+\|D_W^{(t)}\|_\infty\right)n\rho_n^t.
    \label{eq:uniform-expected-codegree-bound}
\end{align}
Conditionally on all the latent variables $U_1,U_2,\dots$, the
variables $\left\{\prod_{a=1}^tA_{s(a),j}\right\}_{j\in [n]\setminus\bar s}$ are
independent Bernoulli random variables. Their total conditional
variance is at most their total conditional mean, which equals the
corresponding sum in \eqref{eq:uniform-expected-codegree-bound}. The
same Bernstein and Borel--Cantelli argument shows that, almost surely,
for all sufficiently large $n$,
\begin{align}
    \max_{s\in[n]_t}\sum_{j\in [n]\setminus\bar s}
    \prod_{a=1}^tA_{s(a),j}
    &\le 2\left(1+\|D_W^{(t)}\|_\infty\right)n\rho_n^t.
    \label{eq:random-graphon-codegree}
\end{align}

\textbf{Claim 1.} We have
\begin{align}
    \left\|
        \E\left[
            K_{\overline W_n}^{(t)}\ind_{\Delta_{n,t}}
            \,\middle|\,U_1,U_2,\ldots
        \right]
        -K_W^{(t)}
    \right\|_1
    \xrightarrow{\mathrm{a.s.}}0.
    \label{eq:cond-mean-L1conv}
\end{align}
Indeed, by \eqref{eq:codegree-kernel-definition}, $\|K_W^{(t)}\|_1=\|D_W^{(t)}\|_1\le\|D_W^{(t)}\|_\infty<\infty$. To identify the conditional expectation, let $\mathbf x\in I_{s(1)}\times\cdots\times I_{s(t)}$ for some $s=(s(1),\dots,s(t))\in [n]^t$, and let $y\in I_j=\left(\frac{j-1}{n},\frac{j}{n}\right]$.
If $(s(1),\ldots,s(t),j)\in[n]_{t+1}$, conditional independence gives
\begin{align*}
    \E\left[
        K_{\overline W_n}^{(t)}(\mathbf x;y)
        \,\middle|\,U_1,U_2,\ldots
    \right]
    =
    \rho_n^{-t}
    \prod_{a=1}^t p_{\sigma_n(s(a)),\sigma_n(j)}
    =
    \prod_{a=1}^t (W\wedge\rho_n^{-1})\bigl(U_{\sigma_n(s(a))},U_{\sigma_n(j)}\bigr).
\end{align*}
When the indices are not all distinct,
$K_{\overline W_n}^{(t)}\ind_{\Delta_{n,t}}$ and
$\Pi_nK_{W\wedge\rho_n^{-1}}^{(t)}$ both vanish. Hence
\begin{align*}
    \E\left[
        K_{\overline W_n}^{(t)}\ind_{\Delta_{n,t}}
        \,\middle|\,U_1,U_2,\ldots
    \right]
    =
    \Pi_nK_{W\wedge\rho_n^{-1}}^{(t)}.
\end{align*}

Under either condition in the theorem,
$W\wedge\rho_n^{-1}\to W$ almost everywhere. Under condition~(i) the convergence is uniform, since
\begin{align*}
    \|W-W\wedge\rho_n^{-1}\|_\infty
    =
    \bigl(\|W\|_\infty-\rho_n^{-1}\bigr)_+
    \longrightarrow0,
\end{align*}
whereas under condition~(ii) it follows from $\rho_n^{-1}\to\infty$, which gives
\begin{align*}
    K_{W\wedge\rho_n^{-1}}^{(t)}
    \longrightarrow K_W^{(t)},
    \quad
    0\le K_{W\wedge\rho_n^{-1}}^{(t)}
    \le K_W^{(t)}.
\end{align*}
Applying \Cref{lem:latent-empirical-approximation} with $d=t+1$ and $\varphi=K_W^{(t)}$
 gives
\begin{align*}
    \Pi_nK_{W\wedge\rho_n^{-1}}^{(t)}
    \longrightarrow K_W^{(t)}
    \quad\text{in }L^1
    \quad\text{almost surely},
\end{align*}
which proves Claim 1.

\textbf{Claim 2.} We have
\begin{align}
    \left\|
        K_{\overline W_n}^{(t)}\ind_{\Delta_{n,t}}
        -
        \E\left[
            K_{\overline W_n}^{(t)}\ind_{\Delta_{n,t}}
            \ \middle|\,U_1,U_2,\ldots
        \right]
    \right\|_{\square,t}\xrightarrow{\mathrm{a.s.}}0.
    \label{eq:centererd-ker-conv}
\end{align}
We first prove the claim for $t=1$. By \eqref{eq:cond-mean-L1conv},
\begin{align*}
    \frac{1}{n^2\rho_n}\sum_{i,j=1}^np_{ij}
    =\|\Pi_n(W\wedge\rho_n^{-1})\|_1
    \xrightarrow{\mathrm{a.s.}}\|W\|_1=1.
\end{align*}
Conditionally on the latent variables,
\citet[Lemma~7.2]{BCCZ2019LpI} gives, for every $\varepsilon>0$,
\begin{align*}
    &\bP\left(
        \|\overline W_n-\Pi_n(W\wedge\rho_n^{-1})\|_\square
        >
        \frac{\varepsilon}{n^2\rho_n}
        \sum_{i,j=1}^np_{ij}
        \,\middle|\,U_1,U_2,\ldots
    \right)\\ &\qquad\le 2^{n+1}
    \exp\left\{
        -\frac{\min\{\varepsilon,\varepsilon^2\}}{24}
        \sum_{i,j=1}^np_{ij}
    \right\}.
\end{align*}
For sufficiently large $n$, the right-hand side is at most
$2^{n+1}\exp\{-c_\varepsilon n^2\rho_n\}$ and is summable because
$n\rho_n\to\infty$. Hence
$\|\overline W_n-\Pi_n(W\wedge\rho_n^{-1})\|_\square
\xrightarrow{\mathrm{a.s.}}0$ by Borel--Cantelli, proving the claim
for $t=1$ under either condition (i) or (ii). 

It remains to prove Claim~2 for $2\le t\le r$. Fix such a $t$. The standard matrix representation of the cut norm for step kernels (see, e.g., \citealp[Sections~3--4]{Janson2013Graphons};
\citealp[Appendix~A.1 and the proof of Lemma~4]{KloppVerzelen2019};
and \citealp[Proposition~9.1]{LyuMukherjee2024}) gives
\begin{align*}
    &\left\|
        K_{\overline W_n}^{(t)}\ind_{\Delta_{n,t}}
        -
        \E\left[
            K_{\overline W_n}^{(t)}\ind_{\Delta_{n,t}}
            \,\middle|\,U_1,U_2,\ldots
        \right]
    \right\|_{\square,t}\\
    &\quad=
    \frac{1}{n^{t+1}\rho_n^t}
    \max_{\substack{\mathcal A\subseteq[n]_t\\ B\subseteq[n]}}
    \left|
        \sum_{\substack{s\in \mathcal A\\ j\in B}}
        \left[
            \prod_{a=1}^t
            A_{s(a),j}-\prod_{a=1}^tp_{s(a),j}
        \right]
    \right|.
\end{align*}
Write $\eta_{ij}:=A_{ij}-p_{ij}$, so that
\begin{align*}
    \prod_{a=1}^t A_{s(a),j} - \prod_{a=1}^t p_{s(a),j}
    =
    \sum_{\varnothing\ne J\subseteq[t]}
    \prod_{a\notin J}p_{s(a),j}
    \prod_{a\in J}\eta_{s(a),j}.
\end{align*}
Also, for $j,k\in[n]$, define
\begin{align*}
    R_{jk}
    &:=\sum_{i=1}^n p_{ij}p_{ik},\qquad
    S_{jk}:=\sum_{i=1}^n \eta_{ij}\eta_{ik}.
\end{align*}
Then
\begin{align}
    &\left\|
        K_{\overline W_n}^{(t)}\ind_{\Delta_{n,t}}
        -
        \E\left[
            K_{\overline W_n}^{(t)}\ind_{\Delta_{n,t}}
            \,\middle|\,U_1,U_2,\ldots
        \right]
    \right\|_{\square,t} \notag\\
    &\quad\le
    \frac{1}{n^{t+1}\rho_n^t}
    \sum_{\varnothing\ne J\subseteq[t]}
    \max_{\substack{\mathcal A\subseteq[n]_t\\ B\subseteq[n]}}
    \left|
        \sum_{s\in\mathcal A}\sum_{j\in B}
        \prod_{a\notin J}p_{s(a),j}
        \prod_{a\in J}\eta_{s(a),j}
    \right| \notag\\
    &\quad\le
    \frac{1}{n^{t+1}\rho_n^t}
    \sum_{\varnothing\ne J\subseteq[t]}
    \max_{\substack{\mathcal A\subseteq[n]_t\\ B\subseteq[n]}}
    |\mathcal A|^{1/2}
    \left[
        \sum_{s\in\mathcal A}
        \left(
            \sum_{j\in B}
            \prod_{a\notin J}p_{s(a),j}
            \prod_{a\in J}\eta_{s(a),j}
        \right)^2
    \right]^{1/2} \notag\\
    &\quad\le
    \frac{1}{n^{t+1}\rho_n^t}
    \sum_{\varnothing\ne J\subseteq[t]}
    n^{t/2}
    \max_{B\subseteq[n]}
    \left[
        \sum_{s\in[n]^t}
        \left(
            \sum_{j\in B}
            \prod_{a\notin J}p_{s(a),j}
            \prod_{a\in J}\eta_{s(a),j}
        \right)^2
    \right]^{1/2} \notag\\
    &\quad=
    \frac{1}{n^{t+1}\rho_n^t}
    \sum_{\varnothing\ne J\subseteq[t]}
    n^{t/2}
    \max_{B\subseteq[n]}
    \left(
        \sum_{j,k\in B}
        R_{jk}^{t-|J|}S_{jk}^{|J|}
    \right)^{1/2} \notag\\
    &\quad\le
    \frac{1}{n^{t+1}\rho_n^t}
    \sum_{\varnothing\ne J\subseteq[t]}
    n^{t/2}
    \left(
        \sum_{j,k=1}^n
        R_{jk}^{t-|J|}|S_{jk}|^{|J|}
    \right)^{1/2}.
    \label{eq:RS-bound}
\end{align}
By \eqref{eq:uniform-expected-codegree-bound}, almost surely, for all
sufficiently large $n$,
\begin{align}
    \max_{j\ne k} R_{jk}
    \le \left(1+\|D_W^{(2)}\|_\infty\right)n\rho_n^2,\quad
    \max_j R_{jj}
    \le 
    \max_j \sum_i p_{ij}
    \le \left(1+\|D_W\|_\infty\right)n\rho_n.
    \label{eq:two-star-prob-bounds}
\end{align}
For $j\ne k$, conditionally on all the latent variables, the summands in
$S_{jk}$ are independent and centered, and
\begin{align*}
    \sum_{i=1}^n
    \operatorname{Var}\left(
        \eta_{ij}\eta_{ik}
        \,\middle|\,U_1,U_2,\ldots
    \right)
    \le \sum_{i=1}^np_{ij}p_{ik}=R_{jk}.
\end{align*}
Since $n\rho_n^2\ge n\rho_n^r\gg\log n$, Bernstein's inequality and
\eqref{eq:two-star-prob-bounds} give, almost surely, for all sufficiently
large $n$ and every $j\ne k$,
\begin{align*}
    \bP\left(
        |S_{jk}|>
        4\sqrt{\left(1+\|D_W^{(2)}\|_\infty\right)
            n\rho_n^2\log n}
        \,\middle|\,U_1,U_2,\ldots
    \right)
    \le 2\exp(-4\log n)=2n^{-4}.
\end{align*}
A union bound over $j\ne k$, and Borel--Cantelli give, almost surely, for sufficiently large $n$,
\begin{align}
    \max_{j\ne k}|S_{jk}|
    &\le4\sqrt{\left(1+\|D_W^{(2)}\|_\infty\right)
        n\rho_n^2\log n},
    &\max_j S_{jj}
    &\le3\left(1+\|D_W\|_\infty\right)n\rho_n.
    \label{eq:centered-two-star-bounds}
\end{align}
Moreover, by \eqref{eq:random-graphon-codegree} and the bound on $\max_j\sum_i p_{ij}$ in \eqref{eq:two-star-prob-bounds}, we have, almost surely for sufficiently large $n$,
\begin{align*}
    \max_j S_{jj}\le\max_j\sum_i(A_{ij}+p_{ij})\le 3\left(1+\|D_W\|_\infty\right)n\rho_n.
\end{align*}
Therefore, it follows that almost surely, for $1\le \ell \le t$ and sufficiently large $n$,
\begin{align*}
    \max_{j\neq k} R_{jk}^{t-\ell}|S_{jk}|^\ell\lesssim(n\rho_n^2)^{t-\ell/2}(\log n)^{\ell/2},\quad
    \max_{j} R_{jj}^{t-\ell}|S_{jj}|^\ell\lesssim (n\rho_n)^t.
\end{align*}
Substitution into \eqref{eq:RS-bound} therefore gives
\begin{align*}
    &\left\|
        K_{\overline W_n}^{(t)}\ind_{\Delta_{n,t}}
        -
        \E\left[
            K_{\overline W_n}^{(t)}\ind_{\Delta_{n,t}}
            \,\middle|\,U_1,U_2,\ldots
        \right]
    \right\|_{\square,t}\lesssim_{t,W}
    \left(\frac{\log n}{n\rho_n^2}\right)^{1/4}
    +\frac{1}{\sqrt{n\rho_n^t}}
    \xrightarrow{\mathrm{a.s.}}0.
\end{align*}
This verifies {\bf Claim 2} for $2\le t\le r$.

Now, for every realization of $G_n$, using the bound $\|F\|_{\square,t}\le \|F\|_1$ for every integrable kernel $F$ gives
\begin{align*}
    \left\|
        K_{\overline W_n}^{(t)}
        -K_{\overline W_n}^{(t)}\ind_{\Delta_{n,t}}
    \right\|_{\square,t}
    \le
    \rho_n^{-t}\bigl|[0,1]^t\setminus\Delta_{n,t}\bigr|\lesssim_t(n\rho_n^t)^{-1}\overset{n\to\infty}{\longrightarrow}0.
\end{align*}
The above display, along with \eqref{eq:cond-mean-L1conv} and
\eqref{eq:centererd-ker-conv}
 yield
\begin{align*}
    \|K_{\overline W_n}^{(t)}-K_W^{(t)}\|_{\square,t}
    \xrightarrow{\mathrm{a.s.}}0,\qquad 1\le t\le r.
\end{align*}
In particular for $t=1$, using the definitions of $q_n$ and $W_n$ (see \eqref{eq:W_n-q-definition}) and $\overline{W}_n=\frac{q_n}{\rho_n}W_n$, we get
\begin{align*}
    \left|\frac{q_n}{\rho_n} - 1\right|
    =\left|\int_{[0,1]^2}
        \left(\overline W_n(x,y)-W(x,y)\right)\,dx\,dy\right|
    \le \|\overline W_n-W\|_{\square}
    \xrightarrow{\mathrm{a.s.}}0.
\end{align*}
This gives $q_n>0$ a.s. for all $n$ large enough, and consequently, using again the fact that $\overline W_n=(q_n/\rho_n)W_n$ gives
\begin{align*}
    \left\|K_{W_n}^{(t)}-K_W^{(t)}\right\|_{\square,t}
    \le
    \left(\frac{\rho_n}{q_n}\right)^t
    \left\|K_{\overline W_n}^{(t)}-K_W^{(t)}\right\|_{\square,t}
    +\left|\left(\frac{\rho_n}{q_n}\right)^t-1\right|
    \left\|K_W^{(t)}\right\|_{\square,t}
    \xrightarrow{\mathrm{a.s.}}0.
\end{align*}
Finally, \eqref{eq:random-graphon-codegree} shows that, almost surely,
for all sufficiently large $n$,
\begin{align*}
    \left\|D_{W_n}^{(t)}\ind_{\Delta_{n,t}}\right\|_\infty
    \le
    2\left(1+\|D_W^{(t)}\|_\infty\right)
    \left(\frac{\rho_n}{q_n}\right)^t
    \le4\left(1+\|D_W^{(t)}\|_\infty\right),
\end{align*}
This proves the asserted bound. Choosing $M$ strictly larger than the
right-hand sides verifies \cref{ass:bounded-codegrees} after discarding
finitely many initial terms, and completes the proof.
\end{proof}

\section{Proof of Lemma \ref{lem:multijoin_bound}}
\label{sec:join_bounds}

\begin{definition}\label[definition]{def:graph-degeneracy}
Let $d$ be a nonnegative integer.
A finite simple graph $H$ is called \emph{$d$-degenerate}
if its vertices admit an ordering $v_1,\ldots,v_{|V(H)|}$
such that, for each $j$, at most $d$ indices $i>j$ satisfy
$v_jv_i\in E(H)$.
Such an ordering is called a \emph{$d$-degeneracy ordering}.
\end{definition}

To prove \cref{lem:multijoin_bound}, we will use the following auxiliary results.

\begin{lemma}\label[lemma]{lem:degen-count-bound}
Let $J$ be a fixed $r$-degenerate finite simple graph, and suppose that the graph sequence
$(G_n)$ satisfies \cref{ass:weights,ass:bounded-codegrees}. Then
\begin{align}
    \InjHom(J,G_n)\le M^{|V(J)|} n^{|V(J)|}q_n^{|E(J)|}.
    \label{eq:degenerate-count-bound}
\end{align}
\end{lemma}

\begin{proof}
Let $v_1,\ldots,v_m$ be an $r$-degeneracy ordering of $J$, where
$m=|V(J)|$, and set
\begin{align*}
    d_j
    :=
    \bigl|\{\ell>j:v_jv_\ell\in E(J)\}\bigr|
    \le r.
\end{align*}
Then
\begin{align*}
    \InjHom(J,G_n)
    &=
    \sum_{i_m=1}^n
    \sum_{i_{m-1}\in[n]\setminus\{i_m\}}
    \cdots
    \sum_{i_1\in[n]\setminus\{i_2,\ldots,i_m\}}
    \prod_{j=1}^m
    \prod_{\substack{\ell>j:\\v_jv_\ell\in E(J)}}
    G_n(i_j,i_\ell).
\end{align*}
For fixed distinct $i_{j+1},\ldots,i_m$, \cref{ass:bounded-codegrees} gives
\begin{align*}
    \sum_{i_j=1}^n
    \prod_{\substack{\ell>j:\\v_jv_\ell\in E(J)}}
    G_n(i_j,i_\ell)
    \le M n q_n^{d_j},
\end{align*}
where the same bound holds for $d_j=0$ because $M\ge1$. Applying this
estimate successively from the innermost sum outward and using
$\sum_{j=1}^m d_j=|E(J)|$ completes the proof.
\end{proof}

\begin{lemma}\label[lemma]{lem:biclique-spanning-subgraph}
Let $r\in\bN$, $H\in\{K_{r,k}:k\ge1\}\cup\{K_1,K_{1,1}\}$, and let $\emptyset\ne S\subseteq V(H)$. Then there exists a spanning subgraph
$H_S$ of $H$ (that is, with $V(H_S)=V(H)$ and $E(H_S)\subseteq E(H)$), such that
\begin{align*}
    E(H_S[S])=\emptyset,
    \qquad
    |E(H_S)|\ge |E(H)|-r(|S|-1).
\end{align*}
Moreover, $H_S$ has an $r$-degeneracy ordering in which every vertex of
$V(H)\setminus S$ precedes every vertex of $S$.
\end{lemma}

\begin{proof}
For $H=K_1$, take $H_S=H$. For $H=K_{1,1}$,
keep its edge if $|S|=1$ and delete it if $|S|=2$;
the conclusions follow because $r\ge1$.
Now suppose $H=K_{r,k}$ for some $k\in\bN$, and write $L(H)$ and $R(H)$
for its bipartition classes of sizes $r$ and $k$, respectively.
Set
\begin{align*}
\begin{gathered}
    S_L:=S\cap L(H),\quad S_R:=S\cap R(H),\quad
    a:=|S_L|,\quad b:=|S_R|,\quad s:=a+b=|S|.
\end{gathered}
\end{align*}
Choose $S_R'\subseteq S_R$ with $|S_R'|=\min\{b,r\}$.
Define two edge sets by
\begin{align*}
\begin{aligned}
    E_1:=\{uv:u\in L(H),\
        v\in R(H)\setminus S_R\},\quad
    E_2:=\{uv:u\in L(H)\setminus S_L,\
        v\in S_R'\}.
\end{aligned}
\end{align*}
Let $H_S$ be the spanning subgraph of $H$ with $V(H_S):=V(H)$ and
$E(H_S):=E_1\cup E_2$. Since every edge of $E_1$ has its right endpoint outside $S$ and every edge of $E_2$ has its left endpoint outside $S$, we have $E(H_S[S])=\emptyset$. Moreover, $E_1\cap E_2=\emptyset$ because their right endpoints belong to the disjoint sets $R(H)\setminus S_R$ and $S_R'$, respectively, giving
\begin{align*}
    |E(H_S)|=|E_1|+|E_2|
    =r(k-b)+(r-a)\min\{b,r\}.
\end{align*}
If $b=0$, then $|E(H_S)|=rk$, which gives the desired bound
because $s\ge1$. If $1\le b\le r$, then
\begin{align*}
    |E(H_S)|=rk-ab,
    \qquad ab\le ar\le r(a+b-1)=r(s-1).
\end{align*}
If $b>r$, then, since $r\ge1$,
\begin{align*}
    |E(H_S)|=r(k-b)+(r-a)r
    =rk-r(s-r)\ge rk-r(s-1).
\end{align*}
Thus, in every case,
\begin{align*}
    |E(H_S)|\ge rk-r(s-1)=|E(H)|-r(|S|-1).
\end{align*}

Finally, order the vertices by listing $R(H)\setminus S_R$ first,
then $L(H)\setminus S_L$, and then $S$, with arbitrary order within
each group. A vertex of $R(H)\setminus S_R$ has at most the $r$
vertices of $L(H)$ as later neighbors. For a vertex of
$L(H)\setminus S_L$, every incident edge in $E_1$ leads to an
earlier vertex, while every incident edge in $E_2$ leads to a vertex
of $S_R'$. It therefore has at most $|S_R'|=\min\{b,r\}\le r$
later neighbors. Finally, a vertex of $S$ has no later neighbors,
since $E(H_S[S])=\emptyset$. This is the required
$r$-degeneracy ordering.
\end{proof}

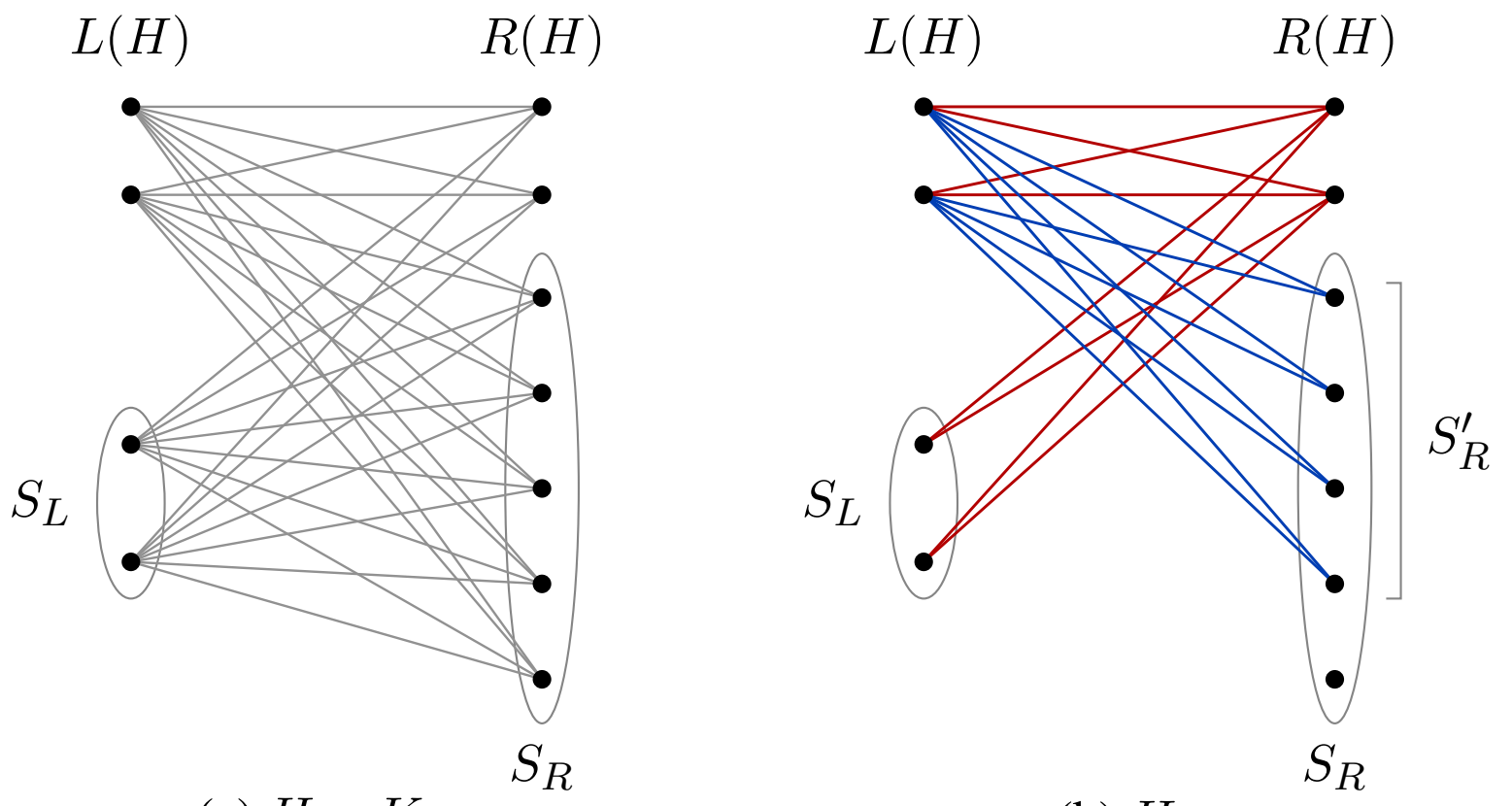
\begin{figure}[H]
\centering
\begin{tikzpicture}[
    x=1cm,y=1cm,font=\small,
    every node/.style={inner sep=1.5pt},
    vertex/.style={circle,fill=black,draw=none,
        inner sep=0pt,minimum size=3.6pt},
    original/.style={draw=black!43,line width=0.45pt},
    first/.style={draw=red!70!black,line width=0.55pt},
    second/.style={draw=blue!65!cyan!70!black,line width=0.55pt},
    boundary/.style={draw=black!48,line width=0.4pt}
]
\path[use as bounding box] (-0.80,-1.25) rectangle (9.30,4.90);

\foreach \p/\dx in {a/0,b/5.40} {
    \begin{scope}[xshift=\dx cm]
        \coordinate (\p-l1) at (0,4.20);
        \coordinate (\p-l2) at (0,3.60);
        \coordinate (\p-l3) at (0,1.90);
        \coordinate (\p-l4) at (0,1.10);
        \coordinate (\p-r1) at (2.80,4.20);
        \coordinate (\p-r2) at (2.80,3.60);
        \coordinate (\p-r3) at (2.80,2.90);
        \coordinate (\p-r4) at (2.80,2.25);
        \coordinate (\p-r5) at (2.80,1.60);
        \coordinate (\p-r6) at (2.80,0.95);
        \coordinate (\p-r7) at (2.80,0.30);

        \node at (0,4.65) {$L(H)$};
        \node at (2.80,4.65) {$R(H)$};

        \draw[boundary] (0,1.50)
            ellipse [x radius=0.23cm,y radius=0.65cm];
        \draw[boundary] (2.80,1.60)
            ellipse [x radius=0.25cm,y radius=1.60cm];

        \node[anchor=east] at (-0.34,1.50) {$S_L$};
        \node at (2.80,-0.30) {$S_R$};
    \end{scope}
}

\foreach \i in {1,2,3,4} {
    \foreach \j in {1,2,3,4,5,6,7} {
        \draw[original] (a-l\i) -- (a-r\j);
    }
}

\foreach \i in {1,2,3,4} {
    \foreach \j in {1,2} {
        \draw[first] (b-l\i) -- (b-r\j);
    }
}

\foreach \i in {1,2} {
    \foreach \j in {3,4,5,6} {
        \draw[second] (b-l\i) -- (b-r\j);
    }
}

\foreach \p in {a,b} {
    \foreach \v in {l1,l2,l3,l4,r1,r2,r3,r4,r5,r6,r7} {
        \node[vertex] at (\p-\v) {};
    }
}

\begin{scope}[xshift=5.40cm]
    \draw[boundary] (3.15,3.00) -- (3.25,3.00)
        -- (3.25,0.85) -- (3.15,0.85);
    \node[anchor=west] at (3.36,1.925) {$S_R'$};
\end{scope}

\node at (1.40,-0.68) {(a) $H=K_{4,7}$};
\node at (6.80,-0.68) {(b) $H_S$};
\end{tikzpicture}

\caption{Construction of $H_S$ for $H=K_{4,7}$, with $|S_L|=2$
and $|S_R|=5$. The outlined sets have union $S=S_L\cup S_R$.
The retained edges are partitioned into $E_1$ (red), joining $L(H)$
to $R(H)\setminus S_R$, and $E_2$ (blue), joining
$L(H)\setminus S_L$ to $S_R'$. Here $|S_R'|=4$.
The bottom-right vertex belongs to $S_R\setminus S_R'$ and is
isolated in $H_S$.}
\label{fig:biclique-spanning-subgraph}
\end{figure}

\begin{lemma}\label[lemma]{lem:union-spanning-subgraph}
Let $r,m,n\in\bN$, and let $F_1,\ldots,F_m$ be finite simple graphs, each isomorphic to a graph in
$\{K_{r,k}:k\ge1\}\cup\{K_1,K_{1,1}\}$, with $\bigcup_{j=1}^m V(F_j)\subseteq [n]$. Let $\mathcal{G}=\mathcal{G}(F_1,\cdots,F_m)$ be a graph on vertex set $[m]$, with an edge $ij$ whenever $V(F_i)\cap V(F_j)\ne\emptyset$, and suppose that $\mathcal{G}$ is connected. Let $J:=\bigcup_{j=1}^m F_j$, set $v:=\sum_{j=1}^m|V(F_j)|$ and $e:=\sum_{j=1}^m|E(F_j)|$,
and write $h:=v-|V(J)|-(m-1)$.
Then $h\ge0$, and $J$ has an $r$-degenerate spanning subgraph
$J'$ with $|E(J')|\ge e-rh$.
\end{lemma}

\begin{proof}
Set $v_h:=|V(J)|=v-(m-1)-h$.
By connectedness of the intersection graph $\mathcal{G}$, we may order
$F_1,\ldots,F_m$ so that $S_j:=V(F_j)\cap\left(\bigcup_{i<j}V(F_i)\right)\ne\emptyset$
for every $2\le j\le m$. Since $V(F_j)\setminus\bigcup_{i<j}V(F_i)=V(F_j)\setminus S_j$ for $2\le j\le m$, we have $v_h=|V(F_1)|+\sum_{j=2}^m(|V(F_j)|-|S_j|)$, and hence
\begin{align}\label{eq:sum-Sa-deficit}
    \sum_{j=2}^m|S_j|
    &=v-v_h=m-1+h,
    \qquad
    \sum_{j=2}^m(|S_j|-1)=h.
\end{align}

Since each $S_j$ is nonempty, \eqref{eq:sum-Sa-deficit} also gives $h\ge0$.
Set $F_1':=F_1$. For each $j\ge2$,
\cref{lem:biclique-spanning-subgraph} with $H=F_j$ and $S=S_j$ gives a spanning subgraph
$F_j'\subseteq F_j$ satisfying
\begin{align*}
    E(F_j'[S_j])=\emptyset,
    \qquad
    |E(F_j')|\ge |E(F_j)|-r(|S_j|-1),
\end{align*}
with an $r$-degeneracy ordering in which $V(F_j)\setminus S_j$ precedes $S_j$. Define the spanning subgraph $J'\subseteq J$ by $V(J'):=V(J),\ E(J'):=\bigcup_{j=1}^m E(F_j')$.

For $i<j$, if $uv\in E(F_i')\cap E(F_j')$, then
$u,v\in V(F_i)\cap V(F_j)\subseteq S_j$, so $uv\in E(F_j'[S_j])$, contradicting $E(F_j'[S_j])=\emptyset$.
Thus the sets $E(F_1'),\ldots,E(F_m')$ are pairwise disjoint.
Using $F_1'=F_1$, the bounds on $|E(F_j')|$, and
\eqref{eq:sum-Sa-deficit}, we obtain
\begin{align}
    |E(J')|
    =\sum_{j=1}^m|E(F_j')|
    \ge e-r\sum_{j=2}^m(|S_j|-1)
    =e-rh.
    \label{eq:spanning-degenerate-subgraph-bound}
\end{align}

Construct an ordering of $V(J')$ as follows.
\begin{enumerate}
    \item Choose an $r$-degeneracy ordering
    $v_1,\ldots,v_{|V(F_1')|}$ of $F_1'=F_1$ (such an ordering exists since $F_1$ is a biclique).
    \item For $j=2,\ldots,m$, let $v_1,\ldots,v_\ell$
    be the ordering of $\bigcup_{i<j}V(F_i')$ constructed
    in the previous step.
    Write the chosen $r$-degeneracy ordering of $F_j'$ as
    $u_1,\ldots,u_t,u_{t+1},\ldots,u_{|V(F_j')|}$, where
    $\{u_1,\ldots,u_t\}=V(F_j')\setminus S_j$ and
    $\{u_{t+1},\ldots,u_{|V(F_j')|}\}=S_j$.
    Since the vertices of $S_j$ already occur among
    $v_1,\ldots,v_\ell$, replace the preceding ordering with
    $u_1,\ldots,u_t,v_1,\ldots,v_\ell$,
    which lists each vertex of $\bigcup_{i=1}^jV(F_i')$
    exactly once.
\end{enumerate}
We verify inductively that after step $j$, the constructed
ordering is an $r$-degeneracy ordering of $\bigcup_{i=1}^jF_i'$. The case $j=1$ holds by construction.
Fix $2\le j\le m$ and assume that $v_1,\ldots,v_\ell$
is an $r$-degeneracy ordering of $\bigcup_{i<j}F_i'$.
We check the ordering $u_1,\ldots,u_t,v_1,\ldots,v_\ell$
of $\bigcup_{i=1}^jF_i'$ in two cases.

\smallskip
\noindent
\emph{Vertices $u_a$, $1\le a\le t$.}
Every edge incident to $u_a$ in this union belongs to
$F_j'$, since $u_a\notin\bigcup_{i<j}V(F_i')$.
The vertices of $F_j'$ occurring after $u_a$ are exactly
$u_{a+1},\ldots,u_{|V(F_j')|}$.
Thus the required bound follows from the chosen
$r$-degeneracy ordering of $F_j'$.

\smallskip
\noindent
\emph{Vertices $v_a$, $1\le a\le\ell$.}
The vertices occurring after $v_a$ are
$v_{a+1},\ldots,v_\ell$.
No edge $v_av_b$ belongs to $F_j'$, since this would imply
$v_a,v_b\in S_j$, contrary to $E(F_j'[S_j])=\emptyset$.
Hence the edges from $v_a$ to these vertices are unchanged
from $\bigcup_{i<j}F_i'$, and the required bound follows
from the induction hypothesis.
\end{proof}

\begin{proof}[Proof of \cref{lem:multijoin_bound}]
Write $\mathcal K_{n,m}^{(h)}$ for
$\mathcal K_{n,m}^{(h)}(H_1,\ldots,H_m)$ and set
$v_h:=v-(m-1)-h$, so that for every $\mathbf s=(s_1,\dots,s_m)\in\mathcal K_{n,m}^{(h)}$, we have $|\cup_{j=1}^m \bar s_j|=v_h$ (see \eqref{eq:v}). We may assume that
$\mathcal K_{n,m}^{(h)}\ne\emptyset$, so $1\le v_h\le v$.

For $\mathbf s\in\mathcal K_{n,m}^{(h)}$, let
$J_{\mathbf s}$ denote the graph $\bigcup_{j=1}^m s_j(H_j)$, which has  vertex set contained in $[n]$ and has size $v_h$.
Since $\mathcal G(\mathbf s)$ is connected,
\cref{lem:union-spanning-subgraph}, applied to
$s_1(H_1),\ldots,s_m(H_m)$, gives an $r$-degenerate
spanning subgraph $J_{\mathbf s}'\subseteq J_{\mathbf s}$
with $|E(J_{\mathbf s}')|\ge e-rh$.
By \cref{ass:weights}, deleting edges can only increase
the weight of each injective map. Thus
\cref{lem:degen-count-bound} and $0<q_n\le1$ give
\begin{align*}
    \InjHom(J_{\mathbf s},G_n)
    \le \InjHom(J_{\mathbf s}',G_n)
    \le M^{v_h}n^{v_h}q_n^{|E(J_{\mathbf s}')|}
    \le M^{v_h}n^{v_h}q_n^{e-rh}.
\end{align*}

Let $\mathcal C_h$ be the family of graphs on $[v_h]$
isomorphic to some $J_{\mathbf s}$ with
$\mathbf s\in\mathcal K_{n,m}^{(h)}$.
The same bound holds for every $F\in\mathcal C_h$.
For any fixed labeled union $F$ and an injective placement of its vertices
in $[n]$, there are at most $v_h^v\le v^v$ choices of the maps from
the disjoint vertex sets of $H_1,\ldots,H_m$ into $V(F)$.
Moreover, $|\mathcal C_h|\le2^{\binom{v_h}{2}}\le2^{\binom v2}$.
Since each edge weight lies in $[0,1]$, repeated occurrences of an edge
can only decrease the product weight. Therefore
\begin{align*}
    \sum_{\mathbf s\in\mathcal K_{n,m}^{(h)}}
        \prod_{j=1}^mG_n(H_j,s_j)
    &\le v^v\sum_{F\in\mathcal C_h}\InjHom(F,G_n)\\
    &\le2^{\binom v2}v^vM^{v_h}n^{v_h}q_n^{e-rh}\\
    &\le2^{\binom v2}(vM)^v
       n^{v-(m-1)}q_n^e(nq_n^r)^{-h},
\end{align*}
where $M>1$ and $v_h\le v$. This is the asserted bound.
\end{proof}

\section{Proof of Theorem \ref{thm:raic-decomposable-bound}}\label{app:stein-bounds}

We first record the Gaussian Stein equation and the bounds on its
solution used below.

\begin{lemma}[{\cite[Lemmas 1, 2]{Meckes2009}}]\label[lemma]{lem:meckes-solution}
Let $\Sigma\in\mathbb R^{m\times m}$ be symmetric and nonnegative
definite, and let $Z_\Sigma\sim\mathcal N(0,\Sigma)$. Suppose that
$\varphi\in C^2(\mathbb R^m)$ and that
\begin{align*}
    \sup_{x\in\mathbb R^m}\|\nabla\varphi(x)\|
    +
    \sup_{x\in\mathbb R^m}
    \|\nabla^2\varphi(x)\|_{\mathrm{op}}
    +
    \operatorname{Lip}(\nabla^2\varphi)
    <\infty.
\end{align*}
Define
\begin{align}
    g_\varphi(x)
    &:=
    -\int_0^1\frac{1}{2t}
    \left[
        \E\varphi\left(
            \sqrt t\,x+\sqrt{1-t}\,Z_\Sigma
        \right)
        -
        \E\varphi(Z_\Sigma)
    \right]\,dt.
    \label{eq:semidefinite-stein-solution}
\end{align}
Then $g_\varphi\in C^2(\mathbb R^m)$,
\begin{align}
    \left\langle
        \Sigma,\nabla^2g_\varphi(x)
    \right\rangle_{\mathrm{HS}}
    -
    \left\langle
        x,\nabla g_\varphi(x)
    \right\rangle
    &=
    \varphi(x)-\E\varphi(Z_\Sigma).
    \label{eq:semidefinite-stein-equation}
\end{align}
Moreover,
\begin{align}
    \sup_x\|\nabla g_\varphi(x)\|
    &\le\sup_x\|\nabla\varphi(x)\|,\notag\\
    \sup_x\|\nabla^2g_\varphi(x)\|_{\mathrm{op}}
    &\le\frac12\sup_x\|\nabla^2\varphi(x)\|_{\mathrm{op}},\notag\\
    \operatorname{Lip}(\nabla^2g_\varphi)
    &\le\frac13\operatorname{Lip}(\nabla^2\varphi).
    \label{eq:semidefinite-bound}
\end{align}
\end{lemma}

\begin{proof}[Proof of \cref{thm:raic-decomposable-bound}]
We follow the algebraic decomposition in the proof of
\cite[Theorem~2]{Raic2004Decomposable}. The preceding lemma supplies the required Stein solution bounds for the covariance matrix $\Sigma$, without requiring $\Sigma$ to be invertible. Let
\begin{align*}
    \mathcal I
    &:=
    \left\{
        (H,s):
        H\in\mathcal H,\ s\in[n]_{V(H)}
    \right\},
    \qquad
    \xi_\alpha:=\xi_n(H,s),
    \quad \alpha=(H,s).
\end{align*}
Define on $\mathcal I$ the graph in which $\alpha=(H,s)$ and $\beta=(H',s')$ are adjacent whenever
$\bar s\cap\bar s'\neq\varnothing$. This is a dependency graph for $\{\xi_\alpha\}_{\alpha\in \mathcal I}$. If the right-hand side of the claimed bound is infinite, there is nothing to prove. We may therefore suppose that it is finite. In
particular, all the expectations appearing below are finite. Let $\varphi$ belong to the class of functions in the definition of $d_3$, and let $Z_\Sigma\sim\mathcal N(0,\Sigma)$, and $g=g_\varphi$ be as in \cref{lem:meckes-solution} and thus satisfying \eqref{eq:semidefinite-stein-equation}, i.e., such that
\begin{align}
    \E\varphi(S)-\E\varphi(Z_\Sigma)
    &=
    \E\left[
        \left\langle
            \Sigma,\nabla^2g(S)
        \right\rangle_{\mathrm{HS}}
        -
        \left\langle
            S,\nabla g(S)
        \right\rangle
    \right].
    \label{eq:stein-discrepancy-S}
\end{align}
Moreover, because $\operatorname{Lip}(\nabla^2\varphi)\le1$, then $\operatorname{Lip}(\nabla^2g)\le\frac13$. For $\alpha\in\mathcal I$, let $K_\alpha$ be the closed neighborhood of $\alpha$ in the dependency graph and set
\begin{align}
    Z_\alpha
    &:=
    \sum_{\beta\in K_\alpha}\xi_\beta,
    \quad
    S_\alpha:=S-Z_\alpha.
    \label{eq:Za-Sa}
\end{align}
The random vector $S_\alpha$ is independent of $\xi_\alpha$.
Since $\E\xi_\alpha=0$, it follows that
\begin{align}
    \Sigma
    =
    \E SS^\top
    =
    \sum_{\alpha\in\mathcal I}
    \E\left[\xi_\alpha S^\top\right]
    =
    \sum_{\alpha\in\mathcal I}
    \E\left[\xi_\alpha Z_\alpha^\top\right]
    =
    \sum_{\alpha\in\mathcal I}
    \sum_{\beta\in K_\alpha}
    \E\left[\xi_\alpha\xi_\beta^\top\right].
    \label{eq:sigma-dependency-decomposition}
\end{align}
Similarly,
\begin{align*}
    \E\left\langle S,\nabla g(S)\right\rangle
    &=
    \sum_{\alpha\in\mathcal I}
    \E\left\langle
        \xi_\alpha,
        \nabla g(S)-\nabla g(S_\alpha)
    \right\rangle,
\end{align*}
because $\E\left\langle\xi_\alpha,\nabla g(S_\alpha)\right\rangle = \left\langle \E\xi_\alpha,\E\nabla g(S_\alpha)\right\rangle = 0$. Also, the fundamental theorem of calculus gives
\begin{align*}
    \nabla g(S)-\nabla g(S_\alpha)
    &=
    \int_0^1
    \nabla^2g(S_\alpha+tZ_\alpha)Z_\alpha\,dt.
\end{align*}
Since $\nabla^2g$ is bounded and the summands are square-integrable, the integrand is integrable. Fubini's theorem
therefore yields
\begin{align}
    \E\left\langle S,\nabla g(S)\right\rangle
    &=
    \sum_{\alpha\in\mathcal I}
    \sum_{\beta\in K_\alpha}
    \int_0^1
    \E\left\langle
        \xi_\alpha,
        \nabla^2g(S_\alpha+tZ_\alpha)\xi_\beta
    \right\rangle\,dt.
    \label{eq:gradient-dependency-decomposition}
\end{align}

For $\alpha\in\mathcal I$ and $\beta\in K_\alpha$, define
\begin{align}
    V_{\alpha\beta}
    &:=
    \sum_{\gamma\in K_\beta\setminus K_\alpha}\xi_\gamma,
    \qquad
    S_{\alpha\beta}
    :=S-Z_\alpha-V_{\alpha\beta}=S - \sum_{\gamma\in K_\alpha\cup K_\beta}\xi_\gamma,
    \label{eq:Vab}
\end{align}
and in particular $S_{\alpha\beta}$ is independent of the pair $(\xi_\alpha,\xi_\beta)$. Let $(\widehat\xi_\alpha)_{\alpha\in \mathcal I}$ be an independent copy of $\{\xi_\alpha\}_{\alpha\in \mathcal I}$.
Equation \eqref{eq:sigma-dependency-decomposition} gives
\begin{align}
    \E\left\langle
        \Sigma,\nabla^2g(S)
    \right\rangle_{\mathrm{HS}}
    &=
    \sum_{\alpha\in\mathcal I}
    \sum_{\beta\in K_\alpha}
    \E\left\langle
        \widehat\xi_\alpha,
        \nabla^2g(S)\widehat\xi_\beta
    \right\rangle.
    \label{eq:covariance-independent-copy}
\end{align}
Because $S_{\alpha\beta}$ is independent of
$(\xi_\alpha,\xi_\beta)$ and the copied pair has the same joint
distribution,
\begin{align*}
    \E\left\langle
        \widehat\xi_\alpha,
        \nabla^2g(S_{\alpha\beta})
        \widehat\xi_\beta
    \right\rangle
    &=
    \E\left\langle
        \xi_\alpha,
        \nabla^2g(S_{\alpha\beta})\xi_\beta
    \right\rangle.
\end{align*}
Thus, combining
\eqref{eq:stein-discrepancy-S},
\eqref{eq:gradient-dependency-decomposition}, and
\eqref{eq:covariance-independent-copy}, and adding and subtracting the two equal quantities in the preceding display, we obtain
\begin{align}
    \E\varphi(S)-\E\varphi(Z_\Sigma)
    &=
    \sum_{\alpha\in\mathcal I}
    \sum_{\beta\in K_\alpha}
    \E\left\langle
        \widehat\xi_\alpha,
        \left[
            \nabla^2g(S)-\nabla^2g(S_{\alpha\beta})
        \right]
        \widehat\xi_\beta
    \right\rangle\notag\\
    &\quad-
    \sum_{\alpha\in\mathcal I}
    \sum_{\beta\in K_\alpha}
    \int_0^1
    \E\left\langle
        \xi_\alpha,
        \left[
            \nabla^2g(S_\alpha+tZ_\alpha)
            -
            \nabla^2g(S_{\alpha\beta})
        \right]
        \xi_\beta
    \right\rangle\,dt.
    \label{eq:raic-decomposition-semidefinite}
\end{align}
Note that
\begin{align*}
    S-S_{\alpha\beta}=Z_\alpha+V_{\alpha\beta},\quad
    S_\alpha+tZ_\alpha-S_{\alpha\beta}=
    V_{\alpha\beta}+tZ_\alpha,
\end{align*}
and recall that for a matrix $A$ and vectors $u,v\in\mathbb R^m$, we have $\left|\langle u,Av\rangle\right|\le \|A\|_{\mathrm{op}}\|u\|\|v\|$. Then
\begin{align*}
    \left|
        \E\left\langle
            \widehat\xi_\alpha,
            \bigl(\nabla^2g(S)-\nabla^2g(S_{\alpha\beta})\bigr)
            \widehat\xi_\beta
        \right\rangle
    \right|& \le \E\left[\|\nabla^2g(S)-\nabla^2g(S_{\alpha\beta})\|_{\mathrm{op}}\|\widehat\xi_\alpha\|
        \|\widehat\xi_\beta\|\right]\\
    &\overset{\eqref{eq:semidefinite-bound}}{\le}
    \frac13
    \E\left[
        \|Z_\alpha+V_{\alpha\beta}\|
        \|\widehat\xi_\alpha\|
        \|\widehat\xi_\beta\|
    \right]\\
    &=
    \frac13
    \E\left[\|\xi_\alpha\|\cdot\|\xi_\beta\|\right]
    \E[\|Z_\alpha+V_{\alpha\beta}\|].
\end{align*}
Similarly,
\begin{align*}
    &\left|
        \E\left\langle
            \xi_\alpha,
            \bigl(\nabla^2g(S_\alpha+tZ_\alpha)
                -\nabla^2g(S_{\alpha\beta})\bigr)\xi_\beta
        \right\rangle
    \right|
    \le \frac13
    \E\left[
        \|\xi_\alpha\|\cdot\|\xi_\beta\|\cdot\|V_{\alpha\beta}+tZ_\alpha\|
    \right].
\end{align*}
Applying the two estimates above to
\eqref{eq:raic-decomposition-semidefinite} gives
\begin{align*}
    &\left|\E\varphi(S)-\E\varphi(Z_\Sigma)\right|\\
    &\quad\le\frac13
    \sum_{\alpha\in\mathcal I}\sum_{\beta\in K_\alpha}
    \Biggl\{
        \E\bigl[\|\xi_\alpha\|\,\|\xi_\beta\|\bigr]
        \E\|Z_\alpha+V_{\alpha\beta}\|\\
    &\hspace{4.2cm}
        +\int_0^1\E\bigl[
            \|\xi_\alpha\|\,\|\xi_\beta\|\,
            \|V_{\alpha\beta}+tZ_\alpha\|
        \bigr]\,dt
    \Biggr\}\\
    &\quad\le\frac13
    \sum_{\alpha\in\mathcal I}\sum_{\beta\in K_\alpha}
    \sum_{\gamma\in K_\alpha\cup K_\beta}
    \Bigl\{
        \E\bigl[\|\xi_\alpha\|\,\|\xi_\beta\|\bigr]
        \E\|\xi_\gamma\|
        +\E\bigl[\|\xi_\alpha\|\,\|\xi_\beta\|\,\|\xi_\gamma\|\bigr]
    \Bigr\}.
\end{align*}
The second inequality follows from
$Z_\alpha+V_{\alpha\beta}
=\sum_{\gamma\in K_\alpha\cup K_\beta}\xi_\gamma$ and the triangle
inequality. Reintroducing the copy indices yields
\begin{align*}
    &\left|\E\varphi(S)-\E\varphi(Z_\Sigma)\right|\\
    &\quad\le\frac13\sum_{H_1,H_2,H_3\in\mathcal H}
       \sum_{(s_1,s_2,s_3)\in\mathcal K_{n,3}(H_1,H_2,H_3)}
       \Biggl\{\E\prod_{a=1}^3\|\xi_n(H_a,s_a)\|\\
    &\hspace{4.5cm}
       +\E\left[\prod_{a=1}^2\|\xi_n(H_a,s_a)\|\right]
          \E\|\xi_n(H_3,s_3)\|\Biggr\}.
\end{align*}
Here every triple with $\beta\in K_\alpha$ and
$\gamma\in K_\alpha\cup K_\beta$ is connected; enlarging the sum to
all connected triples is valid because the summands are nonnegative.
Taking the supremum over $\varphi$ then completes the proof.
\end{proof}

\end{document}